\documentclass[reqno, 11pt]{amsart}
\pdfoutput=1

\usepackage{eucal}
\usepackage{outlines}
\usepackage{bbm}
\usepackage{amssymb,amsfonts}
\usepackage{stmaryrd}

\usepackage{autobreak}
\usepackage{enumerate}
\usepackage{enumitem}
\usepackage{mathrsfs}
\usepackage{tikz-cd}
\usepackage{mathtools}
\usepackage{bm}
\usepackage{quiver} 
\usepackage{scalerel}
\usepackage{stackengine,wasysym}
\usepackage{csquotes}
\usepackage{amsmath,calligra,mathrsfs}
\usepackage[style=alphabetic,
minalphanames=1,
maxalphanames=7,
maxbibnames=99]{biblatex}
\usepackage{hyperref}
\usepackage[capitalise]{cleveref}

\newtheorem{thm}{Theorem}[section]

\newtheorem{cor}[thm]{Corollary}
\newtheorem{prop}[thm]{Proposition}
\newtheorem{lem}[thm]{Lemma}
\newtheorem{conj}[thm]{Conjecture}

\newtheorem{claim}[thm]{Claim}
\newtheorem{mainthm}{Main Theorem}

\theoremstyle{definition}

\newtheorem{defn}[thm]{Definition}

\newtheorem{con}[thm]{Construction}
\newtheorem{exmp}[thm]{Example}

\newtheorem{notn}[thm]{Notation}

\newtheorem{disc}[thm]{Discussion}

\newtheorem{warn}[thm]{Warning}

\newtheorem{conv}[thm]{Convention}

\newtheorem*{mainidea*}{Main Idea}
\newtheorem*{disc*}{Discussion}

\theoremstyle{remark}

\newtheorem{warning}[thm]{Warning}
\newtheorem{sch}[thm]{Scholium}

\newtheorem{rem}[thm]{Remark}

\newtheorem{notation}[thm]{Notation}
\newtheorem{convention}[thm]{Convention}
\newtheorem{construction}[thm]{Construction}

\newtheorem{conventions}[thm]{Conventions}

\makeatletter
\let\c@equation\c@thm
\makeatother
\numberwithin{equation}{section}

\newcommand{\beqn}{\begin{equation}}
	\newcommand{\eeqn}{\end{equation}}
\newcommand{\bclaim}{\begin{claim}}
	\newcommand{\eclaim}{\end{claim}}
\newcommand{\blem}{\begin{lem}}
	\newcommand{\elem}{\end{lem}}
\newcommand{\bproof}{\begin{proof}}
	\newcommand{\eproof}{\end{proof}}
\newcommand{\bdef}{\begin{defn}}
	\newcommand{\edefn}{\end{defn}}
\newcommand{\bprop}{\begin{prop}}
	\newcommand{\eprop}{\end{prop}}
\newcommand{\bthm}{\begin{thm}}
	\newcommand{\ethm}{\end{thm}}
\newcommand{\brem}{\begin{rem}}
	\newcommand{\erem}{\end{rem}}
\newcommand{\bcor}{\begin{cor}}
	\newcommand{\ecor}{\end{cor}}

\newcommand{\bbD}{\mathbb{D}}
\newcommand{\bbE}{\mathbb{E}}

\newcommand{\bbG}{\mathbb{G}}

\newcommand{\bbL}{\mathbb{L}}

\newcommand{\bbN}{\mathbb{N}}

\newcommand{\bbS}{\mathbb{S}}

\newcommand{\bbZ}{\mathbb{Z}}

\newcommand{\calA}{\mathcal{A}}

\newcommand{\calD}{\mathcal{D}}

\newcommand{\calG}{\mathcal{G}}
\newcommand{\calH}{\mathcal{H}}

\newcommand{\calL}{\mathcal{L}}
\newcommand{\calM}{\mathcal{M}}

\newcommand{\calO}{\mathcal{O}}

\newcommand{\calR}{\mathcal{R}}
\newcommand{\calS}{\mathcal{S}}

\newcommand{\calU}{\mathcal{U}}

\newcommand{\calX}{\mathcal{X}}
\newcommand{\calY}{\mathcal{Y}}
\newcommand{\calZ}{\mathcal{Z}}

\newcommand{\scrC}{\mathscr{C}}
\newcommand{\scrD}{\mathscr{D}}
\newcommand{\scrE}{\mathscr{E}}

\newcommand{\scrK}{\mathscr{K}}
\newcommand{\scrL}{\mathscr{L}}

\newcommand{\scrR}{\mathscr{R}}

\newcommand{\sfA}{\mathsf{A}}

\newcommand{\frakg}{\mathfrak{g}}

\renewcommand{\a}{\alpha}
\renewcommand{\b}{\beta}

\newcommand{\on}{\operatorname}
\newcommand{\wit}{\widetilde}

\newcommand{\un}{\underline}
\newcommand{\wih}{\widehat}

\newcommand{\G}{\mathbb{G}}
\newcommand{\C}{\mathbb{C}}

\newcommand{\Aone}{\mathbb{A}^1}

\newcommand{\hook}{\hookrightarrow}
\newcommand{\F}{\mathcal{F}}

\newcommand{\Ze}{\mathbb{Z}}

\newcommand{\kb}{k[\b^{\pm 1}]}

\newcommand{\shear}{{\mathbin{\mkern-6mu\fatslash}}}
\newcommand{\unshear}{{\mathbin{\mkern-6mu\fatbslash}}}

\newcommand{\mathendash}{\text{\textendash}}
\renewcommand{\mod}{\mathendash\mathrm{mod}}

\newcommand{\llb}{\llbracket}
\newcommand{\rrb}{\rrbracket}
\newcommand{\llp}{(\!(}
\newcommand{\rrp}{)\!)}
\mathchardef\mhyphen="2D

\renewcommand{\bmod}{\mathendash\mathbf{mod}}

\newcommand{\heart}{\heartsuit}

\newcommand{\oblv}{\mathrm{oblv}}

\newcommand{\rh}{\mathrm{rh}}

\newcommand{\an}{\mathrm{an}}
\newcommand{\Cat}{\mathcal{C}\mathrm{at}}
\newcommand{\id}{\mathrm{id}}
\newcommand{\Map}{\calM\mathrm{ap}}
\newcommand{\End}{\mathrm{End}}
\newcommand{\cl}{\mathrm{cl}}

\newcommand{\sm}{\mathrm{sm}}

\newcommand{\Alg}{\mathrm{Alg}}
\newcommand{\CAlg}{\mathrm{CAlg}}
\newcommand{\Sch}{\mathrm{Sch}}

\newcommand{\dgcat}{\mathbf{dgcat}}

\newcommand{\QCoh}{\on{\mathsf{QCoh}}}
\newcommand{\IndCoh}{\on{\mathsf{IndCoh}}}
\newcommand{\Vect}{\mathsf{Vect}}
\newcommand{\Sing}{\on{\mathsf{Sing}}}
\newcommand{\Coh}{\on{\mathsf{Coh}}}
\newcommand{\Perf}{\on{\mathsf{Perf}}}
\newcommand{\ShvCat}{\on{\mathsf{ShvCat}}}
\newcommand{\act}{\mathsf{act}}
\newcommand{\Ind}{\on{\mathsf{Ind}}}

\newcommand{\Rep}{\on{Rep}}

\newcommand{\Shv}{\on{Shv}}

\newcommand{\pr}{\on{pr}}
\newcommand{\crit}{\on{crit}}
\newcommand{\Fun}{\on{Fun}}

\newcommand{\op}{\on{op}}
\newcommand{\pt}{\on{pt}}

\newcommand{\res}{\on{res}}

\newcommand{\HH}{\on{HH}}

\newcommand{\MF}{\on{MF}}

\newcommand{\HP}{\on{HP}}
\newcommand{\QC}{\on{QC}}

\newcommand{\N}{\on{N}}
\newcommand{\PreMF}{\on{PreMF}}
\newcommand{\aalg}{\mathendash\textnormal{alg}}

\newcommand{\Orl}{\on{Orl}}
\newcommand{\LG}{\on{LG}}

\newcommand{\coker}{\on{coker}}

\newcommand{\unMF}{\un{\MF}}

\newcommand{\unHH}{\un{\HH}}
\newcommand{\unPreMF}{\un{\PreMF}}
\newcommand{\gr}{\on{gr}}

\newcommand{\DR}{\on{DR}}

\newcommand{\et}{\textrm{\'et}}

\DeclareMathOperator{\Mod}{Mod}

\DeclareMathOperator{\Sym}{Sym}
\DeclareMathOperator{\fib}{fib}

\DeclareMathOperator{\Spec}{Spec}
\DeclareMathOperator{\Hom}{Hom}
\DeclareMathOperator{\Sp}{Sp}

\DeclareMathOperator{\sHom}{\mathscr{H}\text{\kern -3pt {\calligra\large om}}\,}

\renewcommand{\over}{\overline}

\newcommand{\Ch}{\on{Ch}}

\newcommand{\proj}{\mathrm{proj}}

\newcommand{\ex}{\mathrm{ex}}

\newcommand{\Grp}{\mathrm{Grp}}

\renewcommand{\O}{\calO}

\newcommand{\wE}{\wih{\scrE}}

\renewcommand{\G}{\bbG}

\newcommand{\Zar}{\mathrm{Zar}}

\newcommand{\mult}{\mathrm{mult}}

\newcommand{\Temp}{\mathcal{RH}}
\newcommand{\red}{\mathrm{red}}

\renewcommand{\sHom}{\calH\mathrm{om}}

\newcommand{\lLoc}{\mathrm{Loc}}
\newcommand{\muRH}{\Temp^{\mu}}
\newcommand{\poly}{\,\mathrm{alg}}
\renewcommand{\res}{\mathrm{res}}
\newcommand{\super}{\mathrm{super}}

\newcommand{\ku}{k[u^{\pm 1}]}
\newcommand{\ka}{k[\a^{\pm 1}]}
\newcommand{\kz}{k[z^{\pm 1}]}
\newcommand{\calMF}{\mathcal{MF}}

\newcommand{\idm}{\mathrm{idm}}
\newcommand{\coAlg}{\mathrm{coAlg}}

\newcommand{\enh}{\mathrm{enh}}
\newcommand{\unOmega}{\underline{\Omega}}
\newcommand{\cmod}{\mathendash\hspace{.05em}\mathrm{cmod}}
\newcommand{\Cb}{\C[\b^{\pm 1}]}
\newcommand{\Stk}{\mathsf{Stk}}

\newcommand{\PreStk}{\mathsf{PreStk}}
\newcommand{\Mon}{\mathrm{Mon}}
\newcommand{\CMon}{\mathrm{CMon}}

\newcommand{\Grpd}{\mathsf{Grpd}}
\newcommand{\aft}{\mathrm{aft}}
\newcommand{\laft}{\mathrm{laft}}
\newcommand{\Noeth}{\mathrm{Noeth}}

\newcommand{\Corr}{\mathrm{Corr}}
\newcommand{\open}{\mathrm{open}}
\newcommand{\Tor}{\mathrm{Tor}}
\renewcommand{\sch}{\mathrm{sch}}
\newcommand{\Artin}{\mathrm{Artin}}
\newcommand{\singsupp}{\on{singsupp}}
\newcommand{\Moduli}{\mathrm{Moduli}}
\newcommand{\Lie}{\mathrm{Lie}}
\newcommand{\LieAlg}{\mathsf{LieAlg}}
\newcommand{\ft}{\mathrm{ft}}
\newcommand{\aff}{\mathrm{aff}}
\newcommand{\BM}{\mathrm{BM}}
\newcommand{\per}{\mathrm{per.}}
\newcommand{\dR}{\mathrm{dR}}
\newcommand{\dual}{\mathrm{dual}}
\newcommand{\Swap}{\calS\mathrm{wap}}
\newcommand{\unHP}{\un{\HP}}
\newcommand{\tr}{\on{tr}}

\newcommand{\nil}{\mathrm{nil}}
\newcommand{\unit}{\mathrm{unit}}
\newcommand{\GStk}{\mathsf{GStk}}

\newcommand{\Z}{\bm{Z}}

\newcommand{\X}{\bm{X}}
\newcommand{\K}{\scrK}
\newcommand{\T}{\bm{T}}

\newcommand{\Y}{\bm{Y}}
\newcommand{\bmN}{\bm{N}}

\DeclareFontFamily{U}{mathx}{\hyphenchar\font45}
\DeclareFontShape{U}{mathx}{m}{n}{
      <5> <6> <7> <8> <9> <10>
      <10.95> <12> <14.4> <17.28> <20.74> <24.88>
      mathx10
      }{}
\DeclareSymbolFont{mathx}{U}{mathx}{m}{n}
\DeclareFontSubstitution{U}{mathx}{m}{n}
\DeclareMathAccent{\widecheck}{0}{mathx}{"71}
\DeclareMathAccent{\wideparen}{0}{mathx}{"75}

\makeatletter
\newcommand*\bigcdot{\mathpalette\bigcdot@{.5}}
\newcommand*\bigcdot@[2]{\mathbin{\vcenter{\hbox{\scalebox{#2}{$\m@th#1\bullet$}}}}}
\makeatother

\usepackage{geometry}
\title{A microlocal Feigin--Tsygan--Preygel theorem}

\author{Kendric Schefers}

\begin{document}

\begin{abstract}
Let $\Z$ be a derived global complete intersection over $\C$. 
We compute the periodic cyclic homology of the category of
ind-coherent sheaves with prescribed singular support on $\Z$ 
in terms of the microlocal homology, a family of chain theories 
living between cohomology and Borel--Moore homology.
Our result is a microlocal generalization of both the Feigin--Tsygan theorem
identifying $\HP_{\bullet}^{\C}(\QCoh(\Z))$ with
the $2$-periodized cohomology of $\Z$ and A. Preygel's 
theorem identifying $\HP_{\bullet}^{\C}(\IndCoh(\Z))$ with
the $2$-periodized Borel--Moore homology of $\Z$.

Our proof strategy makes extensive use of categories of
matrix factorizations, which we treat using Preygel's formalism. 
This paper contains generalizations
of several known results in the subject
which we prove in this formalism, and which may be of independent interest
to the reader.
\end{abstract}

\maketitle

\tableofcontents


\section{Introduction}

\subsection{Overview}

Let $k$ be an algebraically closed field.
The philosophy of noncommutative algebraic geometry
proposes to study $k$-schemes by
studying their dg-categories of quasicoherent sheaves. In doing so,
it is useful to consider not just categories of sheaves, but \emph{all}
complete cocomplete dg-categories. Such a category is
called a noncommutative space, and is thought
of as the proxy for an underlying fictional space $X$, on which it is
the category of quasicoherent sheaves. The prototypical example
of a noncommutative space is $A\mod$ for a noncommutative ring
$A$, which one pretends is the category of sheaves on
the noncommutative scheme ``$\Spec A$."

Since we imagine a noncommutative space as
standing in for its possibly nonexistent underlying 
geometric space,
the typical maneuver in the subject
is to define operations on noncommutative spaces by analogy with the
corresponding operation on $\QCoh(X)$ for an actual scheme or
stack $X$. For example, the famous theorem of Hochschild--Kostant--Rosenberg
identifies the Hochschild homology of $\QCoh(X)$ with the 
space of differential forms on $X$, so we define the space of differential forms
on $\scrC$ to be its $k$-linear Hochschild homology $\HH_{\bullet}^k(\scrC)$.

The Feigin--Tsygan theorem identifies the $k$-linear periodic cyclic homology
of $\QCoh(X)$ with the $2$-periodized de Rham complex of $X$,
which suggests that the appropriate definition for the
de Rham cohomology of a noncommutative space $\scrC$
is the periodic cyclic homology $\HP_{\bullet}^k(\scrC)$,
or $\HP_{\bullet}^{\kb}(\scrC)$, when $\scrC$ is a 
$2$-periodic (i.e. $\kb$-linear) dg-category.
The periodic cyclic homology has several features
which indicate that it indeed behaves like de Rham cohomology; 
for example, there seems to be a kind of noncommutative
Hodge structure on $\HP_{\bullet}^{\kb}(\scrC)$ (\cite{KKP})
just as there is a usual Hodge structure on cohomology.

There is, of course, another category of sheaves other than
$\QCoh$ which one might consider on a scheme\footnotemark
\footnotetext{Schemes and stacks for us will always be derived,
even when the descriptor ``derived" is omitted.} 
$X$: the category of \emph{ind-coherent sheaves} $\IndCoh(X)$, 
a perfectly good noncommutative space itself which is
in some sense dual to $\QCoh(X)$.
The ind-coherent counterpart to the Feigin--Tsygan theorem
is the result of T. Preygel (\cite{Preygel})
identifying the periodic cyclic homology of $\IndCoh(X)$
with the $2$-periodized Borel--Moore complex
on $X$.

For any reasonable 
scheme $X$, there is a fully faithful inclusion $\QCoh(X)
\hook \IndCoh(X)$ which is an equivalence
if and only if $X$ is regular. When $X$ is quasi-smooth,
the difference between $\QCoh(X)$ and $\IndCoh(X)$ is
mediated by the singular support theory of Arinkin--Gaitsgory. 
To each object $\F \in \IndCoh(X)$ is assigned a closed, conic subset 
of $T^*[-1]X$, called the singular support of $\F$, 
which, roughly speaking, measures the (shifted) (co-)directions obstructing
$\F$ from belonging to $\QCoh(X)$.

One may consider
the full subcategories $\IndCoh_{\Lambda}(X) \subset 
\IndCoh(X)$ spanned by objects
whose singular support is contained in a prescribed
subset $\Lambda \subset T^*[-1]X$. The theorems
of Feigin--Tsygan and Preygel may therefore be seen
as identifying the de Rham cohomologies of the extremal
members of a directed system of noncommutative spaces:
$\IndCoh_{X}(X)$ and $\IndCoh_{T^*[-1]X}(X)$,
respectively.

Let $\Z$ be a derived global complete intersection.
The goal of this paper, simply put, is to identify in familiar
terms the de Rham cohomology (i.e. the periodic cyclic homology)
of the family of noncommutative spaces given by the categories,
\[\IndCoh_{\Lambda}(\Z),\]
of ind-coherent sheaves on $\Z$
with singular support contained in a fixed closed conic subset
$\Lambda \subset T^*[-1]\Z$.

We find that the periodic cyclic homology
of the intermediate categories $\IndCoh_{\Lambda}(\Z)$ 
may be computed in terms of the
\emph{microlocal homology} of $\Z$, a family of chain theories
sitting between de Rham cohomology and Borel--Moore
homology of the underlying classical scheme of $\Z$. 
In this way, our theorem places the computations of both Feigin--Tsygan
and Preygel in a context where
each is derived from a computation
performed on the entire (shifted) cotangent bundle
by specializing to the cases of, respectively, the minimal (zero section $\Z$) and
maximal ($T^*[-1]\Z$) closed conic subsets.

\subsection{Background}
Before stating our main results in the next section, 
we briefly recall the definition
of microlocal homology and some results from \cite{Kenny22}
on the singular support of Borel--Moore chains.

We will assume familiarity on the part of the reader with
the theory of singular support for coherent
sheaves introduced in \cite{AG15}. 

In what follows, and in the rest of the paper, all
functors will be implicitly derived.

\subsubsection{Microlocal homology}
A useful feature of oriented, smooth manifolds is
Poincar\'e duality, under which singular cohomology
and Borel--Moore homology are isomorphic. 
When $Z$ is a manifold with singularities, 
Poincar\'e duality fails, and
in general there is only a map
$H^i(Z) \to H^{\BM}_{\dim Z - i}(Z)$,
for each $i$, given by cap product with the fundamental
class $[Z] \in H^{\BM}_{\dim Z}(Z)$.
The microlocal homology is a family of chain
theories introduced by D. Nadler as a tentative
microlocal measure for the failure of Poincar\'e
duality on certain tame singular spaces. The precise
definition is as follows.

\begin{defn}(c.f. \cite{NadlerNotes})
Let $X$ be a smooth
complex algebraic variety of dimension $m$,
and let $V$ be a complex vector space of finite dimension $n$, 
viewed as a scheme. Suppose that $f: X \to V$, and let $\Lambda \subset V^{\vee}$
be a closed conic subset. Then the \emph{$\Lambda$-microlocal 
homology of $f$} is defined to be
\[H_{\bullet}^{\Lambda}(f) := \Gamma_{\Lambda}\mu_0(\un{\C}_X[m]),\]
where $\Gamma_{\Lambda}$ denotes the functor of
global sections with support contained in $\Lambda$,
and $\mu_0$ denotes the Kashiwara--Schapira functor
of microlocalization along $\{0\} \subset V$ (see \cite{KS90}).
\end{defn}

The microlocal homology derives its name from the
following important fact. Using the standard yoga of the six formalism,
one easily sees that
\begin{itemize}
\item $H_{\bullet}^{\{0\}}(f)$ recovers the singular cochain complex of $Z = f^{-1}(0)^{\cl}$ with a shift
by $m-n$;\footnotemark
\footnotetext{Note that $m-n$ is the expected dimension of the zero fiber $Z$.}
\item $H_{\bullet}^{V^{\vee}}(f)$ recovers the Borel--Moore complex of $Z$;
\item and that, for conical $\Lambda' \subset \Lambda$, there 
exists a map $H_{\bullet}^{\Lambda'}(f) \to H_{\bullet}^{\Lambda}(f)$.
\end{itemize}
where $Z$ denotes the classical zero fiber of $f$.
Thus, the $\Lambda$-microlocal homology interpolates
between the cohomology and Borel--Moore homology of $Z$
according to the microlocal parameter $\Lambda$.

The microlocal homology a priori depends on
the map $f$---on the presentation of the variety
$Z$. In \cite{Kenny22}, we show that $H_{\bullet}^{\Lambda}(f)$
depends only on the \emph{derived zero fiber} of $f$,
which we denote by $\Z$, rather than on the
full data of the map itself.
This is achieved by introducing a derived algebro-geometric
variant of the Kashiwara--Schapira microlocalization functor
to define a sheaf on $T^*[-1]\Z$, the underlying variety
of the $-1$-shifted cotangent
bundle of $\Z$, whose sections
with support recover are easily seen to
recover the microlocal homology spaces when $\Z$
is proper.

The main theorem of \cite{Kenny22} shows that
this sheaf is equivalent to the canonical perverse
sheaf of twisted vanishing cycles $\varphi_{\T^*[-1]\Z}$ on $\T^*[-1]\Z$,
considered with its canonically oriented $-1$-shifted
symplectic structure. Since the latter is intrinsic
to $\Z$, we see that the various microlocal homology
spaces $H_{\bullet}^{\Lambda}(f)$ depend only
on the derived structure on $\Z$. We reproduce 
the precise statement below.

\begin{thm}
\label{thmm: intro theoremm}
Let $\Lambda \subset V^{\vee}$ be
a closed conic subset, and suppose $f: X \to V$
is a proper map from a smooth scheme $X$
to a finite dimensional vector space $V$.
Let $Z$ denote the classical zero fiber of $f$,
and let $\Z$ denote derived zero fiber. Let
$\wit{\Lambda}$ denote the intersection of
$Z \times \Lambda \subset Z \times V^{\vee}
\simeq \N_{\Z/X}$ with $T^*[-1]\Z \subset
\N_{\Z/X}$.
Then there exists an equivalence-up-to-shifts,
\[\Gamma_{\wit{\Lambda}}\left(\varphi_{\T^*[-1]\Z}\right) 
	\simeq H_{\bullet}^{\Lambda}(f),\]
where $\Gamma_{\wit{\Lambda}}(-)$ denotes
global sections with support contained in $\wit{\Lambda}$.
\end{thm}

\begin{notn}
In light of the above theorem, we prefer to work with
$\Gamma_{\wit{\Lambda}}(\varphi_{\T^*[-1]\Z})$ rather than the
$\Lambda$-microlocal homology $H_{\bullet}^{\Lambda}(f)$, even
when $f$ is not proper. We denote the former by
$H_{\bullet}^{\wit{\Lambda}}(\Z)$, emphasizing its
dependence on $\Z$, and lack of dependence on $f$.
\end{notn}

\subsubsection{Singular support of Borel--Moore chains}
\Cref{thmm: intro theoremm} allows one to
sensibly talk about a singular support theory for classes
in Borel--Moore homology taking values in closed
conical subsets of $T^*[-1]\Z$.
Intuitively, one should think of the space 
$H_{\bullet}^{\Lambda}(\Z)$---or, more properly, its image in 
$H_{\bullet}^{T^*[-1]\Z}(\Z) = H_{\bullet}^{\BM}(Z)$---as 
comprising the classes in the Borel--Moore homology of $Z$ with singular
support contained in $\Lambda \subset T^*[-1]\Z$.
The points of $\Lambda$ should be 
thought of as the ($-1$-shifted cotangent) 
directions obstructing the classes in $H_{\bullet}^{\Lambda}(\Z)$ 
from having a representative in singular cochains.

\subsection{Main results}
The following are the precise statements
of our main results described earlier
in the introduction.

\begin{mainthm}[``Microlocal Feigin--Tsygan--Preygel theorem", \cref{cor: interpolating}]
\label{mainthm: intro theorem 1}
Let $\Z$ be a proper derived global complete intersection
over $\C$, presented as the derived zero fiber of
a map $f: X \to V$, where $X$ is smooth of dimension $m$
and $V$ is a vector space of dimension $n$. Let $\Lambda \subset V^{\vee}$
be a closed conic subset, and let $\wit{\Lambda} \subset T^*[-1]\Z$ denote
the set-theoretic intersection $Z \times \Lambda \cap T^*[-1]\Z$ taken in
inside the conormal bundle, $N_{\Z/X} = Z \times V^{\vee}$. Then
\[\HP_{\bullet}^{\C}(\IndCoh_{\wit{\Lambda}}(\Z))[\b^{\pm 1}] \simeq 
	(H_{\bullet}^{\Lambda}(f)\llp \xi \rrp [\b^{\pm 1}])[m],\]
where $\b$ and $\xi$ are formal variables in cohomological degree
$2$ and $0$, respectively.
\end{mainthm}

\brem
\Cref{mainthm: intro theorem 1} can be summarized as saying that
the microlocal homology decategorifies the Arinkin--Gaitsgory
singular support of (ind-)coherent sheaves.
\erem

The above theorem is deduced, in combination with 
the results of \cite{Kenny22}, from the following stronger theorem
exhibiting an equivalence between the periodic cyclic homology
of $\IndCoh(\Z)$ as a Zariski sheaf on $T^*[-1]\Z$ and the restriction
of the perverse sheaf $\varphi_{\T^*[-1]\Z}$ to the Zariski topology.

\begin{mainthm}[\cref{thm: main theorem}]
\label{mainthm: intro theorem 2}
Let $\Z$ be a derived global complete
intersection over $\C$, presented as
the derived zero fiber of a map $f: X \to V$,
where $X$ is smooth of dimension $m$ and $V$ is a complex vector
space of dimension
$n$. Let $\tau$ denote the canonical map from the analytic topology 
to the $\G_m$-invariant Zariski topology on $T^*[-1]\Z$.
Then there is an equivalence,
\[\unHP_{\bullet}^{\C}(\IndCoh(\Z))[\b^{\pm 1}][m+n] \simeq \tau_*\varphi_{\T^*[-1]\Z}\llp \xi \rrp[\b^{\pm 1}],\]
of $\C\llp \xi \rrp[\b^{\pm 1}]\mod$-valued sheaves on the
$\G_m$-invariant Zariski topology on $T^*[-1]\Z$.
\end{mainthm}

\begin{disc}
Let us briefly explain the appearance of the degree
$0$ and degree $2$ variables, $\xi$ and $\b$, in the
above main theorems, starting with a reminder of the periodic
cyclic homology.

Given a dualizable object $c$ in a symmetric monoidal category
$\scrC$, the periodic cyclic homology $\HP_{\bullet}^{\scrC}(c)$ is defined to be
the Tate fixed points of the categorical trace $\HH_{\bullet}^{\scrC}(c)$
with respect to its canonical $S^1$-action. What kind of object
$\HP_{\bullet}^{\scrC}(c)$ is depends on the structure of $\scrC$.

For example, when $\scrC$ is the category of $k$-linear stable $\infty$-categories
(so that $c$ is a dualizable dg-category), the trace $\HH_{\bullet}^k(c)$
is object of $k\mod$, and the topological $S^1$-action on trace
upgrades to a $k[S^1]$-module by adjunction.
The invariants of any $k[S^1]$-module is
naturally a module over invariants of the trivial $k[S^1]$-module,
$k^{S^1} = C^*(BS^1;k) \simeq k[u]$, where $u$ has
cohomological degree $2$. Taking the Tate fixed points
corresponds to inverting the generator $u$, and obtains
a module over $k[u^{\pm 1}]$. One sees directly the familiar fact
that $\HP_{\bullet}^k(\scrC)$ is $2$-periodic (i.e. a $k[u^{\pm 1}]$-module) in this way.
On the other hand, when $\scrC$ is the category of 
$2$-periodic dg-categories (i.e. $\kb$-linear
stable $\infty$-categories, for $\b$ a second copy of the generator of $C^*(BS^1;k)$), 
$\HP_{\bullet}^{\kb}(c)$ naturally has the structure of 
a \emph{doubly $2$-periodic} $k$-module,
that is, an object of $k[\b^{\pm 1}, u^{\pm 1}]$.

Now let $k=\C$. As explained in the next section, the proof of \cref{mainthm: intro theorem 2}
proceeds by comparing the $k$-linear periodic cyclic homology 
of $\IndCoh(\Z)$, naturally a plain dg-category, to the $\kb$-linear periodic cyclic
homology of a certain category of matrix
factorizations, naturally a $2$-periodic dg-category. 
Thus, in order to compare the two invariants, one must
$2$-periodize (tensor with another copy of $C^*(BS^1;k)$)
the former. The degree $0$ variable $\xi$ appears because the 
$k[\b^{\pm 1}, u^{\pm 1}]$-module
so-obtained is $u$-complete, and therefore equivalent to
a complete $k\llp \xi \rrp[\b^{\pm 1}]$-module. 
\end{disc}

\brem
It is the our feeling that there should be a
way of relating $\HP_{\bullet}^{\C}(\IndCoh(\Z))$
to the microlocal homology
without passing through a category of matrix factorizations,
and, therefore, without needing to $2$-periodize both sides
of the equivalence in \cref{mainthm: intro theorem 2}.
See the discussion in \cref{ssec: some philosophical
points} below.
\erem

\brem
Decategorification results like \cref{mainthm: intro theorem 1} and
\cref{mainthm: intro theorem 2} have found
applications in the subject of geometric representation theory, such
as in the recent paper \cite{GammageHilburn} on $2$-categories $\calO$,
where recent results in the subject (in the case of \textit{op. cit.}, 
Koszul duality for categories $\calO$) can be 
recovered via decategorification from particular instances of the 
$2$-categorical 3d mirror symmetry correspondence.
\erem

\subsection{Methods}
The passage from ind-coherent sheaves to vanishing cycles
is mediated by the category of matrix factorizations, which
is well-known to categorify vanishing cycles in some
precise sense (see \cite{BRTV} and \cite{Efimov}).

Our results require us to leverage several
results in the subject of matrix factorizations
which, before the writing of this paper, were 
not sufficiently natural in their known formulations 
for our purposes.
In the task of proving sufficiently natural
and sheafified versions of such statements, we
found the formalism of matrix factorizations
via formal group actions on categories
developed in Preygel's thesis \cite{PreygelT} 
to be the most helpful. Since this formalism
is not widely known, and since we believe
it offers some distinct advantages over the
usual formalisms of matrix factorizations (such as
curved dg-modules over a certain curved dg-algebra), 
we briefly recall the idea.

\subsubsection{Preygel's matrix factorizations}
A Landau--Ginzburg pair for us is a pair $(M,w)$
where $M$ is a smooth scheme over an algebraically
closed field $k$ of characteristic $0$, and 
$w: M \to \Aone$ is a regular function. The function
$w$ is typically called the ``superpotential" of the Landau--Ginzburg
pair in application to string theory, where the category
of matrix factorizations is known as the D-branes in 
a Landau--Ginzburg model obtained from the given pair $(M,w)$.

A conceptually straightforward computation involving
the derived center (i.e. the Hochschild cohomology) of $\Coh(M)$
shows that the following data are equivalent,
\[\{\text{$k$-linear $B\wih{\G}_a$-actions on $\Coh(M)$}\} \longleftrightarrow \{\text{regular functions on $M$}\}.\]

\begin{defn}
\label{intro def: matrix factorizations}
The small category of matrix factorizations of the
Landau--Ginzburg pair $(M,w)$ is defined to the
Tate fixed points,
\[\MF(M,w) := \Coh(M)^{tB\wih{\G}_a},\] 
of the $B\wih{\G}_a$-action on $\Coh(M)$
determined by the regular function $w: M \to \Aone$.
The large category of matrix factorizations $\MF^{\infty}(M,w)$ is
defined to be the ind-completion of $\MF(M,w)$.
\end{defn}

The idea behind the above definition is as follows.
Heuristically, a $B\wih{\G}_a$-action
on $\Coh(M)$ is a compatible family of $B\wih{\G}_a$-actions
on the Hom-spaces of $\Coh(M)$, which is the same
as a compatible family of automorphisms of all the 
Hom-spaces. The $B\wih{\G}_a$-action corresponding to $w$
is given by multiplying the Hom-spaces by $w$. 
Taking $B\wih{\G}_a$-invariants amounts to trivializing
these automorphisms, which in this case amounts
to a nullhomotopy of $w$, so it should not
be surprising the $\Coh(M)^{B\wih{\G}_a} \simeq \Coh(M_0)$.
On the other hand, coinvariants of the $B\wih{\G}_a$-action given by $w$
can be shown to be $\Perf(M_0)$, with the natural map
$\Coh(M)_{B\wih{\G}_a} \to \Coh(M)^{B\wih{\G}_a}$
corresponding to the inclusion $\Perf(M_0) \to \Coh(M_0)$.
The Tate fixed points 
\[\Coh(M_0)^{tB\wih{\G}_a} = \coker\left(\Coh(M)_{B\wih{\G}_a} \to \Coh(M)^{B\wih{\G}_a}\right)\]
therefore corresponds to the quotient
$\Coh(M_0)/\Perf(M_0)$. The connection
between $\MF(M,w)$ and the usual
category of matrix factorizations\footnotemark is made
through Orlov's theorem, which famously shows that the singularity
category of $M$, i.e. the dg-quotient $\Coh(M_0)/\Perf(M_0)$,
is equivalent to the usual dg-category of matrix factorizations.
\footnotetext{There are several common variations of
the usual category of matrix factorizations, but here we mean,
roughly, a definition in which the objects are pair of sheaves
(either locally free or coherent) $\F_0, \F_1$ on $M$, 
with morphisms $d_0: \F_0 \to \F_1$ and $d_1: \F_1 \to \F_0$
such that $d_0d_1 = d_1d_0 = w \cdot \id$ and morphisms
are what one would expect. See e.g. \cite{Efimov}, \cite{Orlov},
\cite{Toda}, etc. for precise examples of such a definition.}

\subsubsection{Graded matrix factorizations}
Suppose that $M$ were equipped with $\G_m$-action and
that $w$ were a $\G_m$-equivariant function. Then there is
variant of \cref{intro def: matrix factorizations} corresponding
to the classical category of graded matrix factorizations
which we introduce in the main body of this paper, and whose properties
and relationship to $\MF(M,w)$ we develop. This definition will end up playing
an important role in our proof of the main theorems, since the
category of ind-coherent sheaves on the derived zero fiber of
a map $f: X \to V$ is equivalent to a certain category 
of graded matrix factorizations $\MF(X \times V^{\vee}/\G_m, \wit{f})$---this
is a variant of the main theorem of \cite{Isik} in our formalism.

\subsection{Some philosophical points}
\label{ssec: some philosophical points}
In the theory of $\scrD$-modules and microlocal sheaf theory, 
the vanishing cycles functor and the functors of microlocalization 
and Fourier transform are intimately related. For example, if $f : X \to \C$ 
is a holomorphic function on a complex manifold $X$ and the zero locus of $f$ 
is non-singular, $\varphi_f$ is entirely determined by the microlocalization 
along the submanifold $f^{-1}(0)$ (\cite[Proposition
8.6.3]{KS90}). On the other hand, the stalk of the Fourier transform 
at a generic covector is given by vanishing cycles with 
respect to the functional determined by that covector (\cite{V83a}).
More recently, in \cite{Kenny22}, we have provided a description of
the entire derived microlocalization functor defined therein 
in terms of vanishing cycles. Microlocalization and vanishing
cycles are, in this way, two sides of the same coin. They dovetail
to compute much the same data, and in many circumstances,
can even be identified with one another.

A basic tension present in comparing the Hochschild invariants
of $\MF^{\infty}$ and $\IndCoh$ is that they are each defined using
slightly different input data: the input to define matrix factorizations
is a regular function $w$, whereas the input to define
$\IndCoh$ is a stand-alone scheme (in our case, the zero locus
of a map $f: X \to V$). This tension manifests itself in the tortured formulations
of our main theorems, in which we are required by our
method of proof to $2$-periodize both objects under 
consideration in order to obtain an equivalence.

We speculate that the tension we describe
indicates that, while $\MF^{\infty}(M,w)$
naturally categorifies vanishing cycles $\varphi_w$, 
$\IndCoh(\Z(f))$, where $\Z(f)$ denotes
the derived zero fiber of $f$, more naturally 
categorifies the (quasi-smooth) microlocalization $\mu_{\Z(f)/X}$.
Indeed, while $\MF^{\infty}(M,w)$ and $\varphi_w$
both take the function $w$ as their input,
the $\QCoh(X)$-module $\IndCoh(\Z(f))$
and $\mu_{\Z(f)/X}$ both instead take
the closed immersion $\Z(f) \hook X$
as their input. This point of view is expressed
in the following table:

\begin{center}
\begin{tabular}{ c | c }
 ``microlocalization" & ``vanishing cycles" \\
\hline \\[-6pt]
 $\IndCoh(\Z(f))$ & $\MF^{\infty}(M,w)$ \\  
 $\mu_{\Z(f)/X}$ & $\varphi_w$    
\end{tabular}
\end{center}
in which the items are arranged top-to-bottom
in decreasing categoricity.

From this point of view, our proof of \cref{mainthm: intro theorem 1}
is rather circuitous, first exhibiting a connection
between the Hochschild invariants of $\IndCoh(\Z)$ and 
$\MF^{\infty}(X \times V^{\vee}, \wit{f})$---crossing over
from ``microlocalization" to ``vanishing cycles"---
then exploiting the connection between matrix factorizations
and vanishing cycles, and, finally, concluding with
the literal equivalence of vanishing cycles with quasi-smooth microlocalization
established in \cite{Kenny22}. Diagrammatically, our proof looks as
follows.

\[\begin{tikzcd}
	{\IndCoh(\Z)} && {\MF^{\infty}(X \times V^{\vee}, \wit{f})} \\
	{\mu_{\Z/X}} && {\varphi_{\wit{f}}}
	\arrow["{\text{decategorification}}", squiggly, from=1-3, to=2-3]
	\arrow["{\text{comparing} \HP_{\bullet}}", squiggly, from=1-1, to=1-3]
	\arrow["{\text{equivalence}}", squiggly, from=2-3, to=2-1].
\end{tikzcd}\]

In principle, there should exist a more direct proof that $\IndCoh_{\Lambda}(\Z)$
categorifies $\mu_{\Z/X}$ and, by extension, the microlocal
homology $H_{\bullet}^{\Lambda}(\Z)$. That is, in the diagram
of the structure of our proof above, there should be a squiggly
arrow labeled ``decategorification" going from $\IndCoh(\Z)$ to
$\mu_{\Z/X}$ which does not pass through the ``vanishing cycles"
side of things. Such a proof would, we imagine,
eliminate the need to $2$-periodize the objects in
the equivalence of \cref{mainthm: intro theorem 1}.

More generally, we wonder whether working with $\IndCoh$
as a categorification of derived microlocalization might prove fruitful
in tackling such problems as proving the Joyce conjecture, 
where the relevance of categories of matrix factorizations
in a potential proof of the conjecture is already recognized.

\subsection{Structure of the paper}
\Cref{sec: conventions and notation} establishes the main conventions and
notation used throughout the paper. It also contains a glossary of categories
used in this work that the reader might find helpful.

In \cref{sec: relative matrix factorizations} we recall the notion of
category of relative matrix factorizations over a base
introduced in \cite{PreygelT}, and proves some general facts thereof.

In \cref{sec: matrix factorizations over Spec k and BG_m}
we discuss the two special cases of relative matrix factorizations
of interest to us: relative matrix factorizations over
the base $\Spec k$ and over the base $B\G_m$. 

\Cref{sec: descent} contains various results about the descent
properties of categories and invariants, in particular matrix
factorizations and their periodic cyclic homology, which we use in the proof
of our main result.

In \cref{sec: ind-coherent sheaves and graded matrix factorizations}
we prove an equivalence between the dg-categories of ind-coherent sheaves
on a derived zero fiber $\Z(f)$ of a map $f: X \to V$, on the one hand,
and graded matrix factorizations for a graded
Landau--Ginzburg pair built from the data of $f$, on the other.
We show that our equivalence is actually one of module categories over 
$\QCoh(X \times V^{\vee}/\G_m)$, and therefore sends the singular
support of ind-coherent sheaves to the set-theoretic support of objects
in graded matrix factorizations. Moreover, from our proof it is clear
that this equivalence satisfies Zariski descent. These results have been proven
or asserted at various points in the literature, but not in
our present formalism of matrix factorizations, 
nor with the requisite amount of naturality
needed to show descent easily. 

In \cref{sec: comparing Hochschild invariants} compares the $k$-linear periodic cyclic
homology of graded matrix factorizations for a fixed graded Landau--Ginzburg pair
to the $\kb$-linear periodic cyclic homology of the category 
of matrix factorizations determined by the same pair viewed
as an ungraded Landau--Ginzburg pair. We find that when the domain of
our superpotential is evenly graded, shearing allows us to identify the latter
with the $2$-periodization of the former. Moreover, it is easily
seen that this identification sheafifies.

In \cref{sec: periodic cyclic homology of matrix factorizations}
we compute the sheafified $\kb$-linear periodic cyclic homology of $\MF^{\infty}(M,w)$.
We find, as expected, that this is equivalent to the constructible sheaf of vanishing
cycles with respect to $w$, restricted to the Zariski topology. Our results in
this section are a derived, sheafified version of the original result in \cite{Efimov}
relating the periodic cyclic homology of matrix factorizations to vanishing
cycles.

In \cref{sec: the main theorem and applications}, the last of the main sections,
we finally state and prove \cref{mainthm: intro theorem 1}
and  \cref{mainthm: intro theorem 2}.

There are also several appendices.

\Cref{sec: primer on IndCoh} is a review of the theory of ind-coherent
sheaves. This appendix contains more than is strictly needed to read
the main body of this text.

In \cref{sec: very good stacks} we briefly recall the notion of a very good
stack, and then generalize some results of T. Preygel in order to show that
a certain class of stacks that includes $B\G_m$ is very good.

\Cref{sec: singular support} is a review of the Arinkin--Gaitsgory 
theory of singular support of ind-coherent sheaves. This topic
receives treatment in a section separate from the review of
$\IndCoh$ because of the prominent role
this theory plays throughout this paper, and for ease of access
for readers who are familiar with $\IndCoh$ but whose memories
of the details of the theory of singular support may be hazy.

\Cref{sec: formal groups} collects the necessary facts
about actions of formal groups (both graded and
ungraded) on categories needed to work with
our formalism of matrix factorizations.

\Cref{sec: shearing} is a review of the shearing functors
defined on graded categories. This material is really
only needed for the treatment of evenly graded Landau--Ginzburg
pairs in \cref{sec: comparing Hochschild invariants}.

In \cref{sec: Hochschild invariants}, we recall the definitions of
Hochschild homology and periodic cyclic homology, taking care to
specify the structure on such objects for use in the main body of the paper.
This section also contains a very brief review of the Tate fixed points construction
of Nikolaus--Scholze.

\subsection{Acknowledgments}
I am deeply grateful to David Ben-Zvi, who first introduced
me to the microlocal homology and suggested its connection
to categories of ind-coherent sheaves with prescribed singular support.
I have benefited enormously from his suggestions; in particular, the idea
of using matrix factorizations in this context is due to him.

I am equally grateful
to Sam Raskin for countless helpful discussions about this material,
as well as for lending his immense help with a number of technical
details. I owe to him a great deal of my understanding of the
categorical details in this paper.

I would also like to thank Leonid Positselski
and Yukinobu Toda for answering questions
about their works. 

Finally, I thank Harrison Chen, Tom Gannon, Rok Gregoric, Sam Gunningham,
Peter Haine, Justin Hilburn, Joshua Mundinger, Pavel Safronov, and Germ\'an Stefanich
for their contributions to my understanding of the subject of this work.


\section{Conventions and notation}
\label{sec: conventions and notation}

\subsubsection{}
Fix an algebraically closed field $k$ of characteristic zero, once and for all.

\subsubsection{}
Given a map $f: X \to V$,
$\Z(f)$ will denote the
derived zero fiber of $f$
while $Z(f)$ will denote
the classical zero fiber.
Note that $Z(f)$ is the
underlying classical scheme of
$\Z(f)$.

\subsubsection{Hom-objects and mapping spaces}
The reader will encounter $\infty$-categories in
this paper equipped with several different structures,
each of which we use at different times. Here we
establish notation for the various hom-objects and
mapping spaces associated with these structures.

Let $\scrC$ be an $\infty$-category, and let $a,b \in \scrC$
be two objects. We adopt the following notation.
\begin{enumerate}
\item \label{hom notation 1}
$\Map_{\scrC}(a,b) \in \Grpd_{\infty}$ denotes the mapping
space from $a$ to $b$. We drop the subscript ``$\scrC$" 
when the category is clear from context.

\item \label{hom notation 2}
If $(\scrC, \otimes)$ is closed symmetric monoidal, then $\sHom_{\scrC}(a,b) \in \scrC$ denotes
the internal hom-object from $a$ to $b$, i.e. the object determined by,
\[\Map_{\scrC}(c, \sHom_{\scrC}(a,b)) \simeq \Map_{\scrC}(c \otimes a, b),\]
for $c \in \scrC$. As above, we drop the subscript ``$\scrC$" when the
category is clear from context.

\item \label{hom notation 3} 
If $\scrC$ is a module category over a
rigid symmetric monoidal category $\scrR$, 
then $\Hom_{\scrR}(a,b) \in \scrR$
is the object determined by
\[\Map_{\scrR}(V, \Hom_{\scrR}(a,b)) \simeq \Map_{\scrC}(V \otimes a, b),\]
for $V \in \scrR$, $a, b \in \scrC$.

When $\scrR$ is $R\mod$ for a ring spectrum $R$, we alternatively write
$\Hom_R(a,b)$. For example, if $\scrC$ is an object of $\dgcat_k$, 
then $\Hom_k(a,b)$ is the usual $k$-module spectrum of maps between $a$ and $b$.

\end{enumerate}

We will mainly be interested in categories of sheaves in this paper,
for which we have the following alternative notation.

Let $f: X \to S$ be a derived stack over $S$. Pullback along $f$ gives 
a natural symmetric monoidal functor $f^*:\QCoh(S) \to \QCoh(X)$,
making $\QCoh(X)$ into an algebra over $\QCoh(S)$. Composing
$f^*$ with the canonical symmetric monoidal embedding $\Upsilon: \QCoh(X) \to \IndCoh(X)$
obtains a symmetric monoidal functor $\QCoh(S) \to \IndCoh(X)$,
making $\IndCoh(X)$ into an algebra over $\QCoh(S)$, as well.

\begin{enumerate}

\item We use $\sHom_X(-,-)$ to denote the hom-object
(\labelcref{hom notation 2}) for $(\scrC, \otimes)$ taken to be either
$(\QCoh(X), \otimes)$ or $(\IndCoh(X), \overset{!}{\otimes})$.
When there is potential for confusion, we write $\sHom_{\QCoh(X)}(-,-)$
or $\sHom_{\IndCoh(X)}(-,-)$.

\item We use $\Hom_S(-,-)$ to denote the hom-object of 
(\labelcref{hom notation 3}) above for $\scrR = \QCoh(S)$
and $\scrC$ taken to be either $\QCoh(X)$ or $\IndCoh(X)$.
\end{enumerate}

\brem
There is a useful relationship between $\sHom_{\QCoh(X)}(-,-)$
and $\Hom_S(-,-)$:
\[f_*\sHom_{\QCoh(X)}(-,-) \simeq \Hom_S(-,-).\]
When $S$ is smooth, and $f$ is schematic, then
we similarly have
\[f_*\sHom_{\IndCoh(X)}(-,-) \simeq \Hom_S(-,-).\]

Note that if $X$ is taken to be a stack over itself, i.e.
$S = X$ and $f = \id_X$, then $\Hom_S(-,-) = \sHom_X(-,-)$.
\erem

\brem
In the case when the base stack $S$ is an affine scheme $\Spec R$, 
we also write $\Hom_R(-,-)$ in place of $\Hom_S(-,-)$, by slight abuse of
notation. Indeed note that an action of $\QCoh(\Spec R)$ is an action of $R\mod$,
and that $\Gamma(\Hom_{\Spec R}(-,-)) \simeq \Hom_R(-,-)$.
\erem

\subsubsection{Inclusions}
In order to minimize verbiage, we will use
$i_{?}$ as standing, all-purpose notation 
to denote the inclusion of some space (``$?$")
into another, whenever said inclusion is clear
from context.

\subsection{Glossary of categories}
${}$ \\
\begin{outline}
	\1 $\dgcat_k$ denotes the $(\infty,1)$-category of modules in $\Cat_{\infty}$ over
	$\Perf_k$ which are stable, and for which the action of $\Perf_k$ is exact in
	both variables.
		\2 $\dgcat_k^{\idm}$ denotes the full subcategory of $\dgcat_k$
		on the collection of small stable idempotent-complete $k$-linear $\infty$-categories.
		\2 $\dgcat_k^{\infty}$ denotes the (non-full!) $\infty$-subcategory on 
		the collection of stable cocomplete $k$-linear categories with colimit-preserving
		$k$-linear functors.
	\1 $\Ch_k$ denote the dg-category of strict chain complexes over the field $k$. There is
	a localization map $\Ch_k \to \Vect_k$ inverting quasi-isomorphisms.
	\1 $\QCoh(X)$ denotes the $k$-linear (stable cocomplete) 
	$\infty$-category of quasicoherent complexes on a derived stack $X$. 
	It is equipped with a natural t-structure, whose heart $\QCoh(X)^{\heart}$ 
	is equivalent to the (ordinary) category of quasicoherent complexes on $\pi_0(X)$.
		\2 $\Perf(X) \subset \QCoh(X)$ is the full subcategory spanned
		by perfect complexes.If $X$ is perfect (e.g. a quasi-compact, quasi-separated
		derived scheme) in the sense of \cite{BZFN}, then $\QCoh(X) = \Ind(\Perf(X))$.
		\2 Let $X$ be a coherent derived stack. Then $\Coh(X) \subset \QCoh(X)$ is the full 
		subcategory spanned by objects with locally bounded, coherent cohomology sheaves. 
	\1 Given a derived stack $X$, $\Sing(X)$ denotes the the quotient dg-category
	$\Coh(X)/\Perf(X)$.
		\2 $\Sing^{\infty}(X)$ denotes the quotient of large dg-categories $\IndCoh(X)/\QCoh(X)$. 
		It is readily seen that $\Ind(\Sing(X)) \simeq \Sing^{\infty}(X)$.
	\1 $\IndCoh(X)$ is defined in detail in \cref{sec: primer on IndCoh}, but for sufficiently
	nice $X$, it is simply the ind-completion of $\Coh(X)$ defined above.
	\1 For $A$ a graded $k$-algebra, $A\mod_{\gr}$ denotes the 
	$\Rep(\G_m)$-module category of graded modules over $A$.
		\2 $\Perf_{\gr} A$ denotes the full subcategory of $A\mod_{\gr}$
		spanned by compact objects.
	\1 $\Grpd_{\infty}$ denotes the $\infty$-category of $\infty$-groupoids, whose
	objects are alternatively called \emph{spaces} or \emph{animae}.
	\1 Given an $\infty$-category $\scrC$ with a t-structure, $\scrC^+$ will denote $\cup_{n \geq 0} \scrC^{\geq -n}$,
	and $\scrC^-$ will denote $\cup_{n \geq 0} \scrC^{\leq n}$.
\end{outline}


\section{Relative matrix factorizations}
\label{sec: relative matrix factorizations}
In this paper we adopt the
approach to matrix factorizations envisioned
in \cite{PreygelT}. 


\subsection{Relative Landau--Ginzburg pairs}
\label{ssec: relative Landau--Ginzburg pairs}
Suppose that $S$ is a smooth, very good stack,\footnotemark
\footnotetext{See \cref{def: very good stack} for the definition
of a very good stack.}
and let $\scrL$ denote a fixed line bundle on $S$.

\bdef
\label{def: relative LG pair}
A \textit{relative Landau--Ginzburg pair} is a
pair $(M/S, f)$ consisting of a smooth
$S$-scheme $M$ of finite type
over $S$ and a map of $S$-schemes,
\[f: M \to \scrL.\]
\edefn

\brem
By abuse of notation, we refer to both
the sheaf $\scrL$ and its associated $S$-scheme,
$\Spec_{S}(\Sym_{\O_S}(\scrL))$, by $\scrL$.
\erem

\subsubsection{}
Let $(M/S, f)$ be a chosen relative
Landau--Ginzburg pair over $S$, with
the choice of $\scrL$ implicitly taken when
we choose $S$.

\bdef
We let $\G := \wih{(S \times_{\scrL} S)}_S$ denote
the formal completion of the derived self-intersection of
the zero section of $\scrL$ along the diagonal $\Delta: S
\to S \times_{\scrL} S$.
\edefn

Alternatively, we could have defined $\G$ to be $S \times_{\wih{\scrL}_S} S$,
by the following lemma shown to us by Rok Gregoric.

\blem
\label{lem: alternate defn of G}
Let $X \to Y$ be an arbitrary map of derived prestacks.
Then $\wih{(X \times_{Y} X)}_X \simeq X \times_{\wih{Y}_X} X$.
\elem

\bproof
By definition, one has
\beqn
\label{eqn: eh2}
\wih{(X \times_{Y} X)}_X = (X \times_{Y} X) \times_{(X \times_{Y} X)_{\dR}} X_{\dR},
\eeqn
where the subscript ``$\dR$" denotes the de Rham prestack (see \cite[\S 1.1]{GR17II}).
The functor sending a derived prestack to its de Rham prestack is
a right adjoint, so preserves limits. In particular, we have $(X \times_{Y} X)_{\dR}
\simeq X_{\dR} \times_{Y_{\dR}} X_{\dR}$. Tautologically,
we can write $X_{\dR} \simeq X_{\dR} \times_{X_{\dR}} X_{\dR}$,
so the right-hand side of \labelcref{eqn: eh2} may be rewritten as
\[(X \times_{Y} X) \times_{X_{\dR} \times_{Y_{\dR}} X_{\dR}} (X_{\dR} \times_{X_{\dR}} X_{\dR}),\]
which is shown by some formal nonsense to be
\[(X \times_{X_{\dR}} X_{\dR}) \times_{Y \times_{Y_{\dR}} X_{\dR}} (X \times_{X_{\dR}} X_{\dR})
\simeq X \times_{\wih{Y}_X} X.\]
\eproof

The formal derived scheme $\G$ has the structure of a group 
over $S$ by ``composition of loops." 
It also compatibly has a commutative group 
structure over $S$ by ``pointwise addition" (using 
the additive structure on the fibers of $\scrL$). 
For the sake of clarity we give an explicit
but partial description of each of these structures below.

\begin{con}
\label{con: group operations on G} ${}$ \\
\begin{itemize}
\item
Using the alternative description of $\G$ offered by
\cref{lem: alternate defn of G}, we may identify
the loop composition product $\circ_{\mu} : \G \times_S \G \to \G$ and identity $\id : S \to \G$ 
explicitly as
\[\circ_{\mu} : \G \times_S \G \simeq (S \times_{\wih{\scrL}_S} S) \times_S 
	(S \times_{\wih{\scrL}_S} S) \simeq S \times_{\wih{\scrL}_S} S 
		\times_{\wih{\scrL}_S} S \xrightarrow{\pr_{13}} S \times_{\wih{\scrL}_S} S \simeq \G\]
\[\id : S \xrightarrow{\Delta} S \times_{\wih{\scrL}_S} S \simeq \G\]
The rest of the Segal-monoid structure admits a similar description via projections and diagonals.
A homotopy inverse is given by the explicit anti-isomorphism $i : \G \simeq \G^{\op}$ 
which on underlying spaces can be identified with
\[i : \G = S \times_{\wih{\scrL}_S} S \xrightarrow{\text{switch}} S \times_{\wih{\scrL}_S} S = \G.\]

\item
The pointwise addition product $\circ_{+} : \G \times_S \G \to \G$, identity $0: S \to \G$, and 
inverse $- : \G \simeq \G^{\op}$ may be explicitly identified as follows. The
additive structure on $\scrL$ as an $S$-scheme induces an additive structure
on $\wih{\scrL}_S$ as a formal scheme over $S$ in the standard way.
The commutative diagram
\[\begin{tikzcd}
	{S \times_S S} & {S} \\
	{\wih{\scrL}_S \times_{S} \wih{\scrL}_S} & \wih{\scrL}_S
	\arrow["{\circ_{+}}", from=2-1, to=2-2]
	\arrow[from=1-2, to=2-2]
	\arrow["{=}", from=1-1, to=1-2]
	\arrow["{0_{\wih{\scrL}_S} \times 0_{\wih{\scrL}_S}}"', from=1-1, to=2-1]
\end{tikzcd}\]
gives rise to a map
\[\circ_{+} : \G \times_S \G \simeq (S \times_S S) \times_{\wih{\scrL}_S \times_S \wih{\scrL}_S} (S \times_S S) \to S \times_{\wih{\scrL}_S} S = \G.\]
Analogously, base changing the identity $0_{\scrL}: S \to \wih{\scrL}_S$ and the inverse map $- : \wih{\scrL}_S \to \wih{\scrL}_S$ one obtains maps
\[0 : S = S \times_S S \to S \times_{\wih{\scrL}_S} S = \G\]
\[- : \G = S \times_{\wih{\scrL}_S} S \to S \times_{\wih{\scrL}_S} S = \G\]
and an anti-isomorphism $- : \G \simeq \G^{\op}$.
\end{itemize}
\end{con}

\begin{notn}
We call $\G$ the space of \emph{formal
loops in $\scrL$ based at $0$}.
\end{notn}

The category $\IndCoh(\G)$
obtains two different convolution monoidal structures 
from each the loop composition group structure
and the pointwise addition group structure on $\G$.
More precisely, \cref{con: group operations on G}
provides a lift of $\G$ from $\Sch_{/S}$ to 
$\Mon(\CMon(\Sch_{/S}))$, the category of monoids
in the category of commutative monoids in $\Sch_{/S}$,
which we denote by $(\G, \mu, +)$. Composing
with the lax symmetric monoidal functor,
\[X \mapsto \IndCoh(X), \hspace{5mm} f \mapsto f_*^{\IndCoh},\footnotemark\]
\footnotetext{Note that $f_*^{\IndCoh}$ is naturally $\QCoh(S)$-linear
for $f: X \to Y \in \Sch_{/S}$. Indeed,
$f_*^{\IndCoh}$ has a unique structure of map of $\QCoh(Y)$-modules,
where the $\QCoh(Y)$-action on $\IndCoh(X)$ is obtained via
the symmetric monoidal functor $f^*: \QCoh(Y) \to \QCoh(X)$.
Restricting along the map $\QCoh(S) \to \QCoh(Y)$, we
obtain a natural $\QCoh(S)$-linear structure on $f_*^{\IndCoh}$ for
the canonical $\QCoh(S)$-actions on source and target.
Though we do not show it, we claim that a relative
version over $S$ of the functor $\IndCoh$ appearing
in \S\labelcref{sssec: upgrading to a functor} in
\cref{sec: primer on IndCoh} exists
with all the expected properties. We will denote
this functor by
\[\IndCoh(-/S): \Sch_{/S}^{\aft} \to \QCoh(S)\mod.\]}
one obtains $(\IndCoh(\G), \circ_{\mu}, \circ_+) \in \Alg(\CAlg(\QCoh(S)\mod))$.

\bdef
$(\IndCoh(\G), \circ)$ is the image of $(\IndCoh(\G), \circ_{\mu}, \circ_+)$
under the forgetful functor 
$\oblv: \Alg(\CAlg(\QCoh(S)\mod)) \to \CAlg(\QCoh(S)\mod)$.
\edefn

\subsubsection{Koszul duality}

\blem
\label{lem: general Koszul duality}
There is a symmetric monoidal equivalence,
\[(\IndCoh(\G), \circ) = \left(\IndCoh\left(S \times_{\scrL} S\right)_S, \circ\right) 
	\simeq \left(\Sym_{\O_S}\left(\scrL^{\vee}[-2]\right)\mod(\QCoh(S)), \otimes\right).\]
In particular, we obtain an equivalence of 
the symmetric monoidal subcategories of compact objects,
\[(\Coh(\G), \circ) = \left(\Coh\left(S \times_{\scrL} S\right)_S, \circ\right) 
	\simeq \left(\left(\Sym_{\O_S}\left(\scrL^{\vee}[-2]\right)\mod(\QCoh(S))\right)^c, \otimes\right).\]
\elem

\bproof
See \cite[\S 5.4.6]{PreygelT}.
\eproof

\subsection{Relative matrix factorizations}
\label{ssec: relative matrix factorizations}
Let $(M/S, f)$ be a fixed relative
Landau--Ginzburg pair over $S$, as in the
previous section.

Let $M_{0_{\scrL}} := M \times_{\scrL} S$ denote the
fiber product of $M$ along $f$ with $S$ along 
the zero section of $\scrL$.

\subsubsection{}
There is a right action of $\G$ with its $\mu$-product
(``loop composition") structure on $M_{0_{\scrL}}$, as well as a
compatible action of $\G$ with its $+$-product 
(``pointwise addition") structure on $M_{0_{\scrL}}$.
Their definitions could be given rigorously
in the manner of \cref{con: group operations on G},
though we will not do so.

Applying $\IndCoh(-)$ to the
group actions equips $\IndCoh(M_{0_{\scrL}})$ with
the mutually compatible structures of
right $(\IndCoh(\G), \circ_{\mu})$-module
(under convolution along loop composition)
and $(\IndCoh(\G), \circ_+)$-module
(under convolution along addition) in $\QCoh(S)\mod$.
These two module structures are the same
in the precise sense of the Eckmann--Hilton
argument, whose formulation we recall below.

\blem[{\cite[Lemma 3.1.6]{PreygelT}} ``Eckmann--Hilton"]
Let $\scrC^{\otimes}$ be a symmetric 
monoidal $\infty$-category, and suppose 
$A \in \Alg(\CAlg(\scrC^{\otimes}))$. Let $\overline{A} \in \CAlg(\scrC^{\otimes})$
and $\wit{A} \in \Alg(\scrC^{\otimes})$ be its images under the obvious forgetful functors. 
Set $\scrD^{\otimes} = \overline{A}\mod(\scrC^{\otimes})$, and note that there is a (lax symmetric
monoidal) forgetful functor $\scrD^{\otimes} \to \scrC^{\otimes}$. Then:
\begin{enumerate}[label=(\alph*), ref=\alph*]
\item The other commuting product on $A$ gives rise to a lift of $\wit{A}$ to an object $A' \in \Alg(\scrD^{\otimes})$.
\item There are equivalences of $\infty$-categories
\[\scrD \xleftarrow{\simeq} A'\mod(\scrD) \xrightarrow{\simeq} \wit{A}\mod(\scrC).\]
\end{enumerate}
\elem

\brem
In particular, letting $\scrC^{\otimes} = \QCoh(S)\mod$ with its
$\otimes_S$ symmetric monoidal structure, we obtain from
the Eckmann--Hilton argument an equivalence,
\[(\IndCoh(\G), \circ_{\mu})\mod(\QCoh(S)\mod) \simeq (\IndCoh(\G), \circ)\mod(\QCoh(S)\mod).\]
We therefore will refer to ``the" $\IndCoh(\G)$-module structure on $\IndCoh(M_{0_{\scrL}})$
from now on.
\erem

\subsubsection{}
The $\circ$-monoidal structure on $\IndCoh(\G)$ preserves
compact objects; that is, $(\Coh(\G), \circ)$ is a symmetric
monoidal subcategory of $(\IndCoh(\G), \circ)$. 

\blem
\label{lem: small action}
There exists a $\Perf S$-linear action of $(\Coh(\G), \circ)$ on $\Coh(M_{0_{\scrL}})$
such that passing to ind-objects recovers the $\QCoh(S)$-linear action of
$(\IndCoh(\G), \circ)$ on $\IndCoh(M_{0_{\scrL}})$ defined above.
\elem

\bproof
The argument is sketched in the last
bullet point of \cite[Construction 3.1.5]{PreygelT},
which we now reproduce. 

Observe that the structure
maps of the $\IndCoh(\G)$-action on $\IndCoh(M_{0_{\scrL}})$
(i.e. the action map $\IndCoh(\G) 
\otimes_{\QCoh(S)} \IndCoh(M_{0_{\scrL}}) \to \IndCoh(\G)$
along with the maps encoding the
higher compatibility of this action map
with the monoid structure of $\IndCoh(\G)$)
are affine and finite over $S$ on $\pi_0$, so
the $\QCoh(S)$-linear $\IndCoh$-pushforward
along each of them preserves 
compact objects (i.e. $\Coh$).
Thus, $\Coh(M_{0_{\scrL}})$ is a $\Coh(\G)$-module
in $\Perf S\mod(\dgcat_k^{\idm})$
under convolution which recovers 
the action of $\IndCoh(\G)$
on $\IndCoh(M_{0_{\scrL}})$ by passing to the ind-completions.
\eproof

\bdef
\label{def: matrix pre-factorizations}
We define the category of 
\textit{matrix pre-factorizations}
of the relative Landau--Ginzburg pair 
$(M/S, \scrL, f)$ to be the category,
\[\PreMF(M/S, \scrL, f) := \Coh(M_{0_{\scrL}}),\]
viewed as a module
over $(\Coh(\G), \circ)$
via the action
described above.
\edefn

\bdef
\label{def: matrix factorizations}
We define the category of
\textit{matrix factorizations} of
the relative Landau--Ginzburg pair
$(M/S, \scrL, f)$ to be
\[\MF(M/S, \scrL, f) := \PreMF(M/S, \scrL, f) 
	\otimes_{\Coh(\G)} \Sing(\G).\]
\edefn

\begin{notation}
We denote the categories of relative
matrix (pre-)factorizations by $\PreMF(M/S,f)$
and $\MF(M/S,f)$ whenever the line 
bundle $\scrL \to S$ is clear from context.
\end{notation}

\subsubsection{Large or ``infinite rank" matrix factorizations}
Of course, $\IndCoh(M_{0_{\scrL}})$ also has its
$\IndCoh(\G)$-module structure coming from convolution,
so one might just as well define the \emph{large} category
variants of \cref{def: matrix pre-factorizations}
and \cref{def: matrix factorizations} above.

\bdef
\label{def: PreMF infty}
Let $(M/S, \scrL, f)$ be a relative Landau--Ginzburg
pair over $S$. Then $\PreMF^{\infty}(M/S, \scrL, f)$ is
defined to be the $\IndCoh(\G)$-module $\IndCoh(M_{0_{\scrL}})$.
\edefn

\bdef
\label{def: MF infty}
$\MF^{\infty}(M, \scrL, f) := \PreMF^{\infty}(M, \scrL, f) \otimes_{\IndCoh(\G)} \Sing^{\infty}(\G)$.
\edefn

\brem
These definitions for the categories of matrix 
(pre-)factorizations of a relative Landau--Ginzburg pair
 over $S$ were proposed in \cite[\S 5.4]{PreygelT}.
\erem

\brem
\label{lem: Ind(MF) = MF infty}
An upshot of \cref{lem: small action}
is that we may recover the large category of
matrix (pre-)factorizations from the small
one by passing to its ind-completion:
\begin{align*}
\PreMF^{\infty}(M/S, \scrL, f) \simeq& \Ind(\PreMF(M/S, \scrL, f)) \\
\MF^{\infty}(M/S, \scrL, f) \simeq& \Ind(\MF(M/S, \scrL, f))
\end{align*}
\erem

\brem
Heuristically, one may think of
the objects in $\MF^{\infty}$ as
matrix factorizations of infinite rank,
in contrast with objects in $\MF$ which are finite rank.
\erem


\subsection{Support conditions}
\label{ssec: support conditions}
In this subsection, we discuss
the notion of relative matrix (pre-)factorizations
with support conditions. Throughout this
section, unless otherwise specified,
$S$ will be a very good
stack, $\scrL$ a line bundle on $S$,
and $(M,\scrL, f)$ a relative Landau--Ginzburg
pair over $S$. To ease the notation, 
we also use $M_0$ to denote $M_{0_{\scrL}}$.

\subsubsection{}
Suppose that $Z \subset M_{0_{\scrL}}$ is
a closed $S$-subscheme. Then we may
consider the category 
\[\Coh(M_{0_{\scrL}})_Z := \ker\left(\Coh(M_{0_{\scrL}}) \xrightarrow{j^*} \Coh(M_{0_{\scrL}}\setminus Z)\right),\]
where $j: M_{0_{\scrL}} \setminus Z \to M_{0_{\scrL}}$
is the open inclusion. It is easy to see
that $\Coh(M_{0_{\scrL}})_Z$ inherits
a $\Coh(\G)$-module structure from
$\Coh(M_{0_{\scrL}})$.

\blem
\label{lem: MF with supports}
There is a natural $\Coh(\G)$-action 
on $\Coh(M_{0_{\scrL}})_Z$ such that the
kernel map,
\[\Coh(M_{0_{\scrL}})_Z \to \Coh(M_{0_{\scrL}}),\]
lifts to a $\Coh(\G)$-linear
map for the $\Coh(\G)$-action
on $\Coh(M_{0_{\scrL}})$.
\elem

Because $M \setminus Z$ is
an open substack of $M$, it is
itself a smooth $S$-scheme,
and $(M \setminus Z, f|_{M \setminus Z})$
is also a relative Landau--Ginzburg pair
over $S$. As such, $(M \setminus Z)_{0_{\scrL}}
\simeq M_{0_{\scrL}} \setminus Z$ obtains a
natural action by the formal groupoid
$\G$, and the open inclusion, 
$j: M_{0_{\scrL}} \setminus Z \to M_{0_{\scrL}}$,
is naturally $\G$-equivariant structure
for this action. The $\G$-equivariant
structure of $j$ obtains a $\Coh(\G)$-linear
structure on $j^*$ after applying $\Coh(-)$. Taking the
kernel of $j^*$ in $\Coh(\G)\mod(\Perf S\mod(\dgcat_k^{\idm}))$
yields a $\Coh(\G)$-module in $\Perf S\mod(\dgcat_k^{\idm})$
whose underlying
$\Perf S$-module is $\Coh(M_{0_{\scrL}})_Z$.

In light of the above lemma, we make
the following definitions.

\bdef
\label{def: PreMF with supports}
Given a relative Landau--Ginzburg pair
$(M/S, \scrL, f)$, define the category of \emph{relative matrix
pre-factorizations with support in $Z$}
to be the $\Coh(\G)$-module, 
\[\PreMF_Z(M/S, \scrL, f) := \Coh(M_{0_{\scrL}})_Z,\]
with the module structure described
in \cref{lem: MF with supports}.
\edefn

\bdef
\label{def: MF with supports}
Given a relative Landau--Ginzburg pair
$(M/S, \scrL, f)$, define the category of
\emph{relative matrix factorizations with
support in $Z$} to be the $\Sing(\G)$-module,
\[\MF_Z(M/S, \scrL, f) := \PreMF_Z(M/S,\scrL, f) \otimes_{\Coh(\G)} \Sing(\G).\]
\edefn

\brem
The definitions of relative matrix 
(pre-)factorizations with support
are clearly compatible with the definitions of
relative matrix (pre-)factorizations
in the sense that $\PreMF_{M_{0_{\scrL}}}(M/S,\scrL,f)
\simeq \PreMF(M/S,\scrL, f)$ and $\MF_{M_{0_{\scrL}}}(M/S,\scrL, f)
\simeq \MF(M/S,\scrL, f)$.
\erem

\subsubsection{}
Let $X$ be an $S$-scheme, and let $Z \subset X$
be a closed $S$-subscheme.
By \cite[Proposition 4.1.7(b)]{GaitsgoryIndCoh},
the subcategory $\IndCoh(X)_Z \subset \IndCoh(X)$
identifies with the ind-completion of
\[\Coh(X)_Z := \ker\left(\Coh(X) \xrightarrow{j^*} \Coh(X \setminus Z)\right).\]
The preceding discussion therefore applies \textit{mutatis
mutandis} to large categories after passing to the
ind-completions of the categories under discussion.
In particular, $\PreMF^{\infty}_Z(M/S,\scrL,f)$ and
$\MF^{\infty}_Z(M/S,\scrL,f)$ are defined, and
\begin{align*}
\PreMF^{\infty}_Z(M/S,\scrL,f) &\simeq \Ind(\PreMF_Z(M/S,\scrL,f)) \\
\MF^{\infty}_Z(M/S,\scrL,f) &\simeq \Ind(\MF_Z(M/S,\scrL,f))
\end{align*}
as one might expect.

\section{Matrix factorizations over $\Spec k$ and $B\G_m$}
\label{sec: matrix factorizations over Spec k and BG_m}

Though we have defined the category of relative
matrix factorizations over $S$ in great generality, 
we will only need two special cases of its definition:
$S = \Spec k$ and $S = B\G_m$. The categories
obtained in these two cases we will call, respectively,
\emph{$2$-periodic} and \emph{graded} matrix factorizations.

In this section, we write down some details of their construction;
how they compare to the ``classical" categories of graded 
and ungraded matrix factorizations; and put forth an alternative
definition for these particular categories of matrix factorizations
which will prove very useful in computing Hochschild invariants
in a later section.


\subsection{$2$-periodic matrix factorizations}

An important special case of the set-up
in the previous sections is the one in 
which $S = \Spec k = \ast$
(and so $\scrL = \Aone$ necessarily). In this case,
the definition of matrix factorizations
given in \cref{def: matrix factorizations}
recovers the classical category of
matrix factorizations found in the literature,
as shown in \cite[Proposition 3.4.1]{PreygelT}.
We work out the the details of this special case below.

\begin{convention}
For remainder of this and the following section,
we take $\Spec k$ and $\Aone_k$ to be our choices
of base stack, $S$, and fixed line bundle $\scrL$, respectively.
Furthermore, we denote $\Aone_k$ by $\Aone$,
eliding the subscript.
\end{convention}

\subsubsection{}
Observe that
\[\wih{(0 \times_{\Aone} 0)}_0 \simeq 0 \times_{\Aone} 0 = \Omega_0 \Aone\]
because the diagonal morphism $\Delta: 0 \to 0 \times_{\Aone} 0$
is a surjective closed immersion\footnotemark.
\footnotetext{See \cite[Remark 2.1.2]{Preygel} for
a justification of this assertion.}
Thus, $\G$ (in the notation of the previous section) 
is simply the derived scheme of based loops in $\Aone$ at
the origin, and the convolution structure on
$\IndCoh(\G)$ under this equivalence 
is the same as the one induced by the 
group operation of loop composition
in $\Omega_0 \Aone$.

\subsubsection{}
The Koszul duality equivalence described
in \cref{lem: general Koszul duality} specializes to
the following equivalence in the current set-up:
\beqn
\label{eqn: Koszul duality for S=pt}
(\Coh(\G), \circ) \simeq 
	\left((\Sym_k(k[x][-2])\mod(\Vect_k)\right)^c, \otimes) 
		\simeq (\Perf k[\b], \otimes),
\eeqn
where $k[x] = \O_{{\Aone}^{\vee}}$ and
$\b$ denotes a formal variable in cohomological
degree $2$ (i.e.\ $k[\b]$ is the free
$\bbE_{\infty}$-algebra object in $\Vect_k$
on a single generator in 
cohomological degree $2$).

The following lemma characterizes
the essential image of $\Perf(\G)$ under the above 
equivalence.

\blem
\label{lem: image of QC under KD}
The essential image of $\Perf(\G)$, 
viewed as a subcategory of $\Coh(\G)$
in the usual way, under the Koszul
duality equivalence \labelcref{eqn: 
Koszul duality for S=pt} is the full subcategory
of $\Perf k[\b]$ on locally $\b$-torsion
perfect $k[\b]$-modules. Thus, \labelcref{eqn:
Koszul duality for S=pt} induces a symmetric
monoidal equivalence,
\[(\Sing(\G), \circ) \simeq (\Perf \kb, \otimes).\]
\elem

\bproof
This lemma is a part of \cite[Lemma 3.1.9]{PreygelT},
and its proof may be found in \textit{loc. cit.}
\eproof

\subsubsection{}
Thus, we see that the general definition
of $\PreMF$ specializes in
this case to a $k[\b]$-linear category,
\[\PreMF(M/\ast, \Aone, f) := \Coh(M_0),\]
and the definition of $\MF$
specializes to a $\kb$-linear, i.e. $2$-periodic, category,
\[\MF(M/\ast, \Aone, f) := \Coh(M_0) \otimes_{\Perf k[\b]} \Perf \kb.\]

\begin{notation}
In the sequel, when the base stack $S$
of a relative Landau--Ginzburg pair is $\Spec k$, 
we denote matrix pre-factorizations
by $\PreMF(M,f)$ and matrix factorizations
by $\MF(M,f)$.
\end{notation}





\subsection{Graded matrix factorizations}
\label{ssec: graded matrix factorizations}
The definition of $\MF(M/S,\scrL,f)$ in its
full generality was not strictly necessary for our discussion
of $2$-periodic (i.e. ``usual") matrix factorizations
above. In particular, the choice of a line bundle
over $S$ is a vacuous one in the case when $S$
is a point. Rather, the motivating
example for the definition of 
relative matrix factorizations is
the case when
\[(S,\scrL) = (B\G_m, \O_{B\G_m}(d)),\]
where $B\G_m = (\Spec k)/\G_m$ 
and $d$ is an integer. Note that
$B\G_m$ is a very good stack by
\cref{cor: BG_m is a very good stack}.
The category obtained for Landau--Ginzburg pairs over
$(B\G_m, \O_{B\G_m}(d))$ is the category
of graded matrix factorizations.

\subsubsection{}
The special case when $d=2$ will be very
important for us later on; the following lemma
computes $\Sing(\G)$ explicitly under
Koszul duality.

\blem
\label{lem: Sing is Vect}
Let $\G := B\G_m \times_{\O_{B\G_m}(2)} B\G_m$.
Then
\[(\Sing(\G), \circ) \simeq (\Perf_{\gr} \kb, \otimes)\]
under the equivalence of
\cref{lem: general Koszul duality}. 
\elem

\bproof
The invariants of a $\G_m$-category is given by
the totalization of a certain cosimplicial category, 
so the formation of $\G_m$-invariants commutes
with coequalizers. In particular, we have
\[\Sing(\G_{\Aone}/\G_m) \simeq \Coh(\G_{\Aone}/\G_m)/\Perf(\G_{\Aone}/\G_m)
\simeq \Coh(\G_{\Aone})^{\G_m}/\Perf(\G_{\Aone})^{\G_m} \simeq \Sing(\G_{\Aone})^{\G_m}.\]
Hence, the assertion follows from the fact that
the equivalence $(\Sing(\G_{\Aone}), \circ) \simeq (\Perf \kb, \otimes)$
is $\G_m$-equivariant.
\eproof

\subsubsection{Graded Landau--Ginzburg pairs}
\label{sssec: graded Landau--Ginzburg pairs}
Suppose that $(M,f,d)$ is a graded 
Landau--Ginzburg pair of weight $d$ over $\Spec k$
in the following sense: $\G_m$ acts
on $M$, and $f$ is $\G_m$-equivariant when
we equip $\Aone_k$ with the weight $d$ 
$\G_m$-scaling action for some integer $d$.
Such an object may be regarded
as a relative Landau--Ginzburg pair
over $(B\G_m, \O_{B\G_m}(d))$
as described in \cite[Example 5.4.4]{PreygelT};
that is,
\[\begin{tikzcd}
	{M/\G_m} & {\scrL = \Aone/\G_m} \\
	{S = B\G_m}
	\arrow[from=1-1, to=2-1]
	\arrow[from=1-2, to=2-1]
	\arrow["f", from=1-1, to=1-2]
\end{tikzcd}\]
is a relative Landau--Ginzburg pair
over $B\G_m$.

\begin{notation}
We will typically denote a
relative Landau--Ginzburg pair over $B\G_m$
obtained from a graded
Landau--Ginzburg pair $(M,f,d)$ in the above fashion by
$(M/\G_m, f)$ rather than \\ $(M/\G_m/B\G_m,\O_{B\G_m}(d), f)$.
\end{notation}

\begin{notation}
Similarly, we denote the categories of relative
matrix (pre-)factorizations of \\ $(M/\G_m/B\G_m, \O_{B\G_m}(d), f)$ 
by 
\[\PreMF(M/\G_m,\O_{B\G_m}(d),f) \hspace{1em} \textrm{and} \hspace{1em} \MF(M/\G_m, \O_{B\G_m}(d),f)\] 
rather than 
\[\PreMF(M/\G_m/B\G_m, \O_{B\G_m}(d),f) \hspace{1em} \textrm{and} \hspace{1em} \MF(M/\G_m/B\G_m, \O_{B\G_m}(d),f).\] 
or by $\PreMF(M/\G_m,f)$ and $\MF(M/\G_m,f)$ when the
line bundle $\O_{B\G_m}(d)$ is clear from context.
We refer to these categories as the category of
\emph{graded} matrix (pre-)factorizations.
\end{notation}

\brem
\label{irem}
We record the following simple
but important observations. A
relative Landau--Ginzburg pair $(M/\G_m,f)$ over
$(B\G_m, \O_{B\G_m}(d))$ obtained
from a graded Landau--Ginzburg pair
is related to $(M,f)$ (i.e. the underlying
Landau--Ginzburg pair over a point of
the graded pair $(M,f,d)$) by
the following pullback diagrams:

\[\begin{tikzcd}
	M & {M/\G_m} && {\Aone_k} & {\O_{B\G_m}(d)} && M & {M/\G_m}\\
	\ast & {B\G_m} && \ast & {B\G_m} && {\Aone_k} & {\O_{B\G_m}(d)}
	\arrow[from=1-5, to=2-5]
	\arrow[from=1-4, to=2-4]
	\arrow[from=2-4, to=2-5]
	\arrow[from=1-4, to=1-5]
	\arrow["\lrcorner"{anchor=center, pos=0.125}, draw=none, from=1-4, to=2-5]
	\arrow[from=1-1, to=2-1]
	\arrow[from=2-1, to=2-2]
	\arrow[from=1-2, to=2-2]
	\arrow[from=1-1, to=1-2]
	\arrow["\lrcorner"{anchor=center, pos=0.125}, draw=none, from=1-7, to=2-8]
	\arrow["{[f]}", from=1-8, to=2-8]
	\arrow["f"', from=1-7, to=2-7]
	\arrow[from=2-7, to=2-8]
	\arrow[from=1-7, to=1-8]
	\arrow["\lrcorner"{anchor=center, pos=0.125}, draw=none, from=1-7, to=2-8],
\end{tikzcd}\]
where we have denoted the induced map
$M/\G_m \xrightarrow{f} \O_{B\G_m(d)}$ 
by $[f]$.
\erem


In particular, we have the
following lemma which relates the formal
loop groupoids: $\G_{\Aone_k}$ 
and $\G_{\O_{B\G_m}(2)}$.

\blem
\label{lem: relating loop groupoids}
Let $\G_{\Aone} := 0 \times_{\Aone} 0$ and
$\G_{\O_{B\G_m}(2)} := (B\G_m) \times_{\O_{B\G_m}(2)} (B\G_m)$.
Consider $\G_{\Aone}$ as a scheme
with a $\G_m$-action obtained by taking
the product $0 \times_{\Aone} 0$
in the category of $\G_m$-schemes,
with $\Aone$ having the weight $d$
scaling action. Then there is an
equivalence,
\[\G_{\Aone}/\G_m \xrightarrow{\simeq} \G_{\O_{B\G_m}(2)},\]
that exhibits $\G_{\O_{B\G_m}(2)}$
as the stack quotient of $\G_{\Aone}$
by $\G_m$.
\elem

\bproof
Note that the formation of group invariants is a sifted colimit,
so, almost by definition, quotient by $\G_m$ commutes with products.
\eproof


\subsection{Comparison with classical matrix factorizations}
\subsubsection{Classical category of (ungraded) matrix factorizations}
The connection between $\MF(M,f)$
and the ``classical" category of matrix factorizations\footnotemark
found in the literature is explained by
the following proposition of Preygel.

\footnotetext{There are many variations
on the definition of the category of
matrix factorizations depending on
what kind of finiteness conditions its objects
are required to satisfy or what kind of structure one would like
the category to have.
We will not need any particular such
definition, but here is the broad idea
behind all of them (stated
very imprecisely): given a function,
$f: X \to \Aone$, the category of matrix
factorizations of $f$ consists of pairs, $(\F, d)$,
where $\F$ is a $\Ze/2$-graded sheaf
on $X$ and $d$ is an odd endomorphism of $\F$
such that $d^2$ is multiplication by $f$.}

\bprop[{\cite[Proposition 3.4.1]{PreygelT}}]
\label{prop: MF is idempotent completion of Sing}
The natural functor, 
\[\PreMF(M,f) \to \MF(M,f),\]
factors through the quotient functor, $\Coh(M_0) \to
\Sing(M_0) := \Coh(M_0)/\Perf(M_0)$. 
Furthermore, the induced functor 
of $2$-periodic (i.e.\ $\Perf \kb$-module) 
categories,
\beqn
\label{eqn: Sing to MF}
\Sing(M_0) \to \MF(M,f),
\eeqn 
is an idempotent completion.
\eprop

On the other hand, the category of
singularities of $M_0$ is equivalent
to the classical category of matrix factorizations
by a generalization of Orlov's
Theorem (\cite{Orlov}) found in
\cite{BRTV}. Before stating the
theorem, we recall a few
definitions and results from \textit{loc. cit.}
that are necessary in order to state
the theorem.

\begin{notation}
In what follows, let $A$ be a commutative Noetherian 
regular local ring, $S := \Spec A$, and $\Sch_{/S}$
be the category of schemes of finite type over $S$.
\end{notation}

\bdef
\label{def: LG model}
The category of Landau--Ginzburg \textit{models} over $S$,
denoted by $\LG_S$, is defined to be the subcategory
of $(\Sch_{/S})_{/\Aone_S}$ spanned by those pairs,
\[(p: X \to S, f: X \to \Aone_S),\]
where $p$ is a flat morphism.
\edefn


The category $\LG_S$ has a natural
symmetric monoidal structure, 
denoted by $\boxplus$\footnotemark, coming from
the monoidal structure on $(\Sch_{/S})_{/\Aone_S}$,
\[\boxplus: (\Sch_{/S})_{/\Aone_S} \times (\Sch_{/S})_{/\Aone_S} \to (\Sch_{/S})_{/\Aone_S},\]
given on pairs by
\[((X,f), (Y,g)) \mapsto (X \times_S Y, f \boxplus g).\]

\footnotetext{Details are found in \cite[Construction 2.4]{BRTV}.}

The authors of \cite{BRTV} construct both the
category of singularities and matrix factorizations
as lax symmetric monoidal $\infty$-functors,
\[\Sing, \MF^{\cl}: \LG_S^{\boxplus, \op} \to \dgcat^{\idm}_A.\]
Details of the construction written in a
cogent way and with an enormous
amount of rigor may be found in \S2.2 and \S2.3
of \textit{op. cit.} Assuming familiarity 
with the construction of these functors,
we state the promised generalization
of Orlov's Theorem.

\bthm[{\cite[Theorem 2.49]{BRTV}}]
\label{thm: Orlov's theorem}
There is a lax symmetric natural 
transformation of $\infty$-functors
\[\Orl^{-1, \otimes}: \Sing \to \MF^{\cl} : \LG_S^{\op}
	\to \dgcat^{\idm}_A\]
with the following properties:
\begin{enumerate}
\item
\label{item 1}
$\Orl^{-1, \otimes}$ identifies the symmetric 
monoidal structure of $\Sing(S,0)$ given in \cite[Proposition 2.45]{BRTV} 
with the one of $\MF^{\cl}(S, 0)$ as in \cite[Remark 2.12]{BRTV}.

\item 
\label{item 2}
$\Orl^{-1, \otimes}$ is an equivalence when 
restricted to the subcategory of Landau--Ginzburg pairs $(M,f)$
where $f$ is a non-zero divisor on $M$ (i.e.\ the induced morphism $f \cdot : \O_M \to \O_M$ is a
monomorphism), $M/S$ is separated, and $M$ has the resolution property (i.e.\ every
coherent $\O_M$-module is a quotient of a vector bundle, e.g.\ $M$ regular).

\item In particular, from (\labelcref{item 1}) and (\labelcref{item 2}), 
$\MF^{\cl}(X, f)$ and $\Sing(X,f)$ are then equivalent as
$A[\b^{\pm 1}]$-linear idempotent complete dg-categories.
\end{enumerate}
\ethm

\brem
\label{rem: Orlov's theorem for LG pairs}
As noted in \cite[Remark 2.5]{BRTV},
$\Orl^{-1,\otimes}$ restricts to an
equivalence on pairs $(X,f)$ when $X$
is regular. Furthermore, there is
an equivalence $\Sing(X,f) \simeq \Sing(X_0)$,
where the latter category is the
absolute category of singularities
of the derived zero locus of $f$.

Regularity is not stable under basechange, however,
so $\Orl^{-1,\otimes}$ in general is not
lax symmetric monoidal when restricted to
the subcategory spanned by regular Landau--Ginzburg
models. When $A = k$, however,
the notions of smoothness and regularity
coincide. In this case regular Landau--Ginzburg
models \emph{do} form a monoidal subcategory
of $\LG_S^{\op}$, so $\Orl^{-1,\otimes}$ \emph{does}
restrict to a lax symmetric monoidal equivalence
on regular Landau--Ginzburg models. Note also that, in this case,
\cref{def: LG model} coincides with 
the definition of a relative
Landau--Ginzburg pair over 
the pair $(S = \Spec k, \Aone_k)$ (c.f.
\cref{def: relative LG pair}).
\erem

Using \cref{prop: MF is idempotent completion of Sing},
\cref{thm: Orlov's theorem}, and \cref{rem: Orlov's 
theorem for LG pairs}, we deduce the following
corollary.

\bcor
\label{cor: MF is idempotent completion of MFcl}
Let $A = k$,
and let $\Orl: \MF^{\cl} \to \Sing$
denote the inverse equivalence
to $\Orl^{-1, \otimes}$, restricted
to relative Landau--Ginzburg pairs
over $(\Spec k, \Aone_k)$.
Then the composition,
\[\MF^{\cl}(M,f) \xrightarrow[\simeq]{\Orl} \Sing(M_0) 
	\xrightarrow{\labelcref{eqn: Sing to MF}} \MF(M,f),\]
realizes $\MF(M,f)$ as the idempotent
completion of $\MF^{\cl}(M,f)$.
\ecor

\subsubsection{Classical category of graded matrix factorizations}
In this section we state and prove the graded counterpart
to \cref{prop: MF is idempotent completion of Sing}.

\bprop
\label{prop: graded MF is idempotent completion of graded Sing}
Let $(M/B\G_m,f)$ be a relative Landau--Ginzburg
pair over $(B\G_m, \O_{B\G_m}(d))$. Then the
natural functor,
\[\PreMF(M/B\G_m,f) \to \PreMF(M/B\G_m,f) 
	\otimes_{\Coh(\G)} \Sing(\G) =: \MF(M/B\G_m,f)\]
factors through the quotient functor
$\Coh(M_{0_{\O_{B\G_m}(d)}}) \to \Sing(M_{0_{\O_{B\G_m}(d)}})$.
Furthermore, the induced functor,
\beqn
\label{eqn: graded Sing to graded MF}
\Sing(M_{0_{\O_{B\G_m}(d)}}) \to \MF(M/B\G_m,f)
\eeqn
is an idempotent completion.
\eprop

\bproof
The proof of \cref{prop: MF is idempotent completion of Sing}
found under \cite[Proposition 3.1.4.1]{PreygelThesis} may be easily
modified to work in the $\G_m$-equivariant setting (essentially by
working everywhere with graded modules).
\eproof

\brem
In particular, when $(M/B\G_m,f)$ is obtained
from a graded Landau--Ginzburg pair 
$(M,f)$ over $\Spec k$,
$M_{0_{\O_{B\G_m}(d)}} \simeq (M_0)/\G_m$,
the map \labelcref{eqn: graded Sing to graded MF}
may be written more legibly as
\beqn
\label{eqn: simplified idempotent completion map}
\Sing(M_0/\G_m) \to \MF(M/\G_m,f).
\eeqn
\erem

There is also graded analogue to the theorem of
Orlov relating the category of classical matrix 
factorizations and the category of singularities
which may be stated parallel to \cref{thm: Orlov's theorem},
as we now informally indicate.

In what follows let $A$ be a graded commutative Noetherian 
regular local ring; let $S$ be the stack $(\Spec A)/\G_m$; 
and let $\Sch_{S}$ be the category of schemes of finite type over $S$.

We may define the category of graded Landau--Ginzburg models
over $S$ of weight $d$ in parallel to the definition of plain Landau--Ginzburg
models.

\bdef
\label{def: graded LG model}
Let $\O_S(d)$ denote the total space of the line bundle
$\pt^*\O_{B\G_m}(d)$, where $\pt: S \to B\G_m$ is the
structure morphism $(\Spec A)/\G_m \to \ast/\G_m$.

The category of graded Landau--Ginzburg models of weight 
$d$ over $S$, denoted by $\LG_{S,d}$, is defined 
to be the full subcategory of $(\Sch_{/S})_{/\O_{S}(d)}$ 
spanned by those pairs, 
\[(p: X \to S, f: X \to \O_{S}(d)),\]
where $p$ is a flat morphism.
\edefn

The category of graded Landau--Ginzburg models of weight $d$
has a symmetric monoidal structure on it,
similarly denoted $\boxplus$. By imitating the constructions
found in \cite{BRTV}, one may construct the classical category of
graded matrix factorizations and the graded category of singularities 
as lax symmetric monoidal functors,
\[\Sing, \MF^{\cl}: \LG_{S,d}^{\boxplus, \op} \to \dgcat^{\idm}_{\Perf S}.\]

Orlov's theorem about graded matrix factorizations may
now be stated as follows.

\bthm[{cf. \cite[Theorem 3.10]{OrlovGraded}}]
\label{thm: graded Orlov's theorem}
There is a lax symmetric natural 
transformation of $\infty$-functors
\[\Orl^{-1, \otimes}_{\gr, d}: \Sing \to \MF^{\cl} : \LG_{S,d}^{\op}
	\to \dgcat^{\idm}_{\Perf S}\]
with the following properties:
\begin{enumerate}
\item
\label{item 1}
$\Orl^{-1, \otimes}_{\gr, d}$ identifies the symmetric 
monoidal structure of $\Sing(S,0)$ 
with that of $\MF^{\cl}(S, 0)$.

\item 
\label{item 2}
$\Orl^{-1, \otimes}_{\gr, d}$ is an equivalence when 
restricted to the subcategory of Landau--Ginzburg pairs $(M,f)$
where $f$ is a non-zero divisor on $M$ (i.e.\ the induced morphism $f \cdot : \O_M \to \O_M$ is a
monomorphism), $M/S$ is separated, and $M$ has the resolution property (i.e.\ every
coherent $\O_M$-module is a quotient of a vector bundle).

\item In particular, from (\labelcref{item 1}) and (\labelcref{item 2}), 
$\MF^{\cl}(X, f)$ and $\Sing(X,f)$ are then equivalent as
$A[\b^{\pm 1}]_{\gr}$-linear idempotent complete graded small categories.
\end{enumerate}
\ethm

Using \cref{thm: graded Orlov's theorem} and 
\cref{prop: graded MF is idempotent completion of graded Sing},
one deduces the following graded analogue of \cref{cor: MF is
idempotent completion of MFcl}.

\blem
\label{lem: graded MF is idempotent completion of graded MFcl}
Let $A = k$ with the trivial grading,
and let $\Orl_{\gr, d}: \MF^{\cl} \to \Sing$
denote the inverse equivalence
to $\Orl^{-1, \otimes}_{\gr, d}$, restricted
to relative Landau--Ginzburg pairs
over $(B\G_m, \O_{B\G_m}(d))$.
Then the composition,
\[\MF^{\cl}(M/\G_m,f) \xrightarrow[\simeq]{\Orl} \Sing(M_0/\G_m) 
	\xrightarrow{\labelcref{eqn: graded Sing to graded MF}} \MF(M/\G_m,f),\]
realizes $\MF(M/\G_m,f)$ as the idempotent
completion of $\MF^{\cl}(M/\G_m,f)$.
\elem

\subsection{Alternative definitions}
\label{ssec: alternative definitions}

\subsubsection{$2$-periodic matrix factorizations}
Let $(M,f)$ be a Landau-Ginzburg pair.
An alternative means of defining $\MF(M,f)$
is suggested by the following lemma.

\blem[{\cite[Lemma 5.4.1.2]{PreygelThesis}}]
\label{lem: map gives formal group action}
The datum of a map of schemes
$f: M \to \Aone$ is equivalent to the datum
of an action of the formal
group $B\wih{\G}_a$ on $\Coh(M)$
as a $k$-linear category.
\elem

One might hope for a definition of the category
of matrix factorizations in terms of the formal group
action on $\Coh(M)$ corresponding to the superpotential
$f$. Such a definition indeed exists: $\PreMF(M,f)$ is obtained
by taking $B\wih{\G}_a$-invariants of $\Coh(M)$.

\blem[{\cite[Lemma 5.4.1.3]{PreygelThesis}}]
\label{lem: alt definition of PreMF}
There is an equivalence of $k[\b]$-linear
categories,
\[\PreMF(M,f) \simeq \Coh(M)^{B\wih{\G}_a}.\]
\elem

\subsubsection{Graded matrix factorizations}
On the other hand, $B\wih{\G}_a$ has a natural
$\G_m$-action coming from the weight $d$ scaling action
on the additive group scheme $\G_a$, so one might
expect variants of the above lemmas in the graded setting.
Such variants do exist, at least in one direction in the case
of \cref{lem: map gives formal group action}.

Let $(M/\G_m, f)$ be a graded Landau-Ginzburg pair
of weight $d$; in particular, $f: M \to \Aone$ is equivariant
with respect to the weight $d$ scaling action on $\Aone$.

\blem
\label{lem: equivariant map gives graded formal group action}
The choice of a map $f: M/\G_m \to \Aone/\G_m \simeq \O_{B\G_m}(d)$
over $B\G_m$ determines an 
action of the graded formal group $B\wih{\G}_a/\G_m$
on $\Coh(M/\G_m) \simeq \Coh(M)^{\G_m}$ as a
graded (i.e. $\Rep(\G_m)$-linear) dg-category,
where $\G_m$ acts on $B\wih{\G}_a$ via the weight
$d$ scaling action on $\G_a$.
\elem

\bproof
Pullback along the map $f$ gives a symmetric monoidal
functor 
\[f^*: \QCoh(\Aone/\G_m) \to \QCoh(M/\G_m)\]
which makes $\QCoh(M/\G_m)$ into a $\QCoh(\Aone/\G_m)$-algebra.
This induces a map of graded $\bbE_2$-algebras,
\[\Map_{\QCoh(\Aone/\G_m)}(\mathbbm{1}_{\QCoh(\Aone/\G_m)}, 
	\mathbbm{1}_{\QCoh(\Aone/\G_m)}) \to \HH^{\bullet}_{\Rep(\G_m)}(\QCoh(M/\G_m)),\]
where $\mathbbm{1}_{\QCoh(\Aone/\G_m)}$ is the monoidal unit of $(\QCoh(\Aone/\G_m), \otimes)$.
Clearly the monoidal unit is the structure sheaf of $\Aone/\G_m$, which we
also express as $k[x]_{\gr} := \Sym_{\O_{B\G_m}}k(d)$, where $k(d)$ denotes
$k$ placed in $\G_m$-weight $d$.
By \cref{prop: BVhat/G_m-actions}, the space of graded $\bbE_2$-algebra maps
$\Sym_{\O(B\G_m)} k(d) \to \HH^{\bullet}_{\Rep \G_m}(\Coh(M/\G_m))$
is equivalent to the space of $B\wih{\G}_a/\G_m$-actions on $\Coh(M/\G_m)$,
for the $\G_m$-action on $B\wih{\G}_a$ induced by the weight $d$
scaling action on $\G_a$. 
Thus, the map $f$ induces a $B\wih{\G}_a/\G_m$-action on $\Coh(M/\G_m)$,
as desired.
\eproof

Note that, in the context of \cref{lem: equivariant map gives graded formal group action},
the invariants category $\Coh(M/\G_m)^{B\wih{\G}_a/\G_m}$
has a natural $(\O_{B\G_m})^{B\wih{\G_a}/\G_m} \simeq \O_{B^2\wih{\G}_a} 
\simeq k[\b]_{\gr}$-linear structure. 

\blem
By the previous lemma, a choice of graded Landau--Ginzburg pair
$(M, f, d)$ induces a $B\wih{\G}_a/\G_m$-action on $\Coh(M/\G_m)$.
There is an equivalence of graded $k[\b]$-linear\footnotemark
\footnotetext{i.e. $\Perf_{\gr} k[\b]$-linear graded categories, where
$\b$ has $\G_m$-weight $d$. With these choices
$\Perf_{\gr} k[\b]$ and $\Coh(\G_{\O_{B\G_m}(d)})$
are equivalent as graded categories.}
categories,
\[\PreMF(M/\G_m, \O_{B\G_m}(d), f) \simeq \Coh(M/\G_m)^{B\wih{\G}_a/\G_m}.\]
\elem

\bproof
By an analogue of \cite[Lemma 5.2.0.6]{PreygelThesis},
if $G$ is a graded formal group acting 
on a small graded dg-category $\scrC$,
then the natural functor $i: (\Ind\scrC)^G \to \Ind(\scrC)$ is
monadic, so there is a natural equivalence 
\[(\Ind\scrC)^G \simeq (i \circ i^L)\mod(\Ind(\scrC)),\]
where $i^L$ is the left adjoint to $i$ and preserves
compact objects. The natural functor $\scrC^G \to (\Ind\scrC)^G$
is fully faithful with essential image those objects in $(\Ind\scrC)^G$
whose image under $i$ belong to $\scrC$.
Now since $(\Ind\scrC)^G$ is compactly generated
and the monad $i \circ i^L$ preserves compact objects ($i$
also has a right adjoint and therefore preserves colimits),
$i \circ i^L$ restricts to a monad on the small
category of invariants, $\scrC^G$. The equivalence
\[\scrC^G \simeq (i \circ i^L)\mod(\scrC)\]
now follows from the corresponding equivalence
for ind-categories.

The monad $i \circ i^L$ for the $B\wih{\G}_a/\G_m$-action
on $\Coh(M/\G_m)$ coming from the choice
of graded Landau--Ginzburg pair $(M,f,d)$ is 
naturally identified with the endofunctor given by tensoring
with the algebra $\O_{M_0/\G_m}
\simeq \O_{M/\G_m} \otimes_{\Aone/\G_m} \O_{B\G_m}$
in $\Coh(M/\G_m)$ by \cite[Lemma 3.3.1]{PreygelT}.
Thus, as graded dg-categories, we have
\begin{align*}
\Coh(M/\G_m)^{B\wih{\G}_a/\G_m} 	&\simeq \O_{M_0/\G_m}\mod(\Coh(M/\G_m)) \\
							&\simeq \Coh(M_0/\G_m).
\end{align*}
On the other hand, by \cref{prop: BVhat/G_m-actions}, the $B\wih{\G}_a/\G_m$-invariants
of $\Coh(M/\G_m)$ are naturally a module over $\Perf(B\G_m)^{B\wih{\G}_a} \simeq (\Coh(\Omega_0 {\Aone}^{\vee}/\G_m), \circ)$.
It can seen by inspection 
that this $(\Coh(\Omega_0 {\Aone}^{\vee}/\G_m), \circ)$-linear structure
on $\Coh(M_0/\G_m)$ is the one obtained from the geometric action
of $\Omega_0 {\Aone}^{\vee}/\G_m$
on $M_0/\G_m$. Thus, the $\Perf_{\gr} k[\b]$-linear structure
on $\Coh(M_0/\G_m)$ induced here coincides with that described in
\cref{ssec: relative matrix factorizations}.
\eproof


\section{Descent}
\label{sec: descent}

\subsection{Schematic sites}

\subsubsection{}
Let $S$ be a derived stack over $k$, a fixed
field of characteristic $0$.

\bdef
Let $f: \calX \to \calY$ be a morphism of prestacks
over $S$ (which we also call an $S$-morphism of
prestacks). We say that $f$ is
schematic if for any affine $S$-scheme $X \to S$ and
$S$-morphism $X \to \calY$, the Cartesian product
$X \times_{\calY} \calX$
is a representable by an $S$-scheme.
\edefn

\brem
In the case when $S = \Spec k$, this definition
agrees with the definition of schematic morphism
of prestacks found in \cite[\S 3.6.1]{GR17}.
\erem

\bdef
Let $f: \calX \to \calY$ be a schematic $S$-morphism
of prestacks. We say that $f$ is smooth (resp. \'etale,
an open embedding) 
if for any affine $S$-scheme $X \to S$ and
$S$-morphism $X \to \calY$, the induced map
$X \times_{\calY} \calX \to X$
is a smooth (resp. \'etale, an open embedding)
map of $S$-schemes.
\edefn

\bdef
\label{def: site}
Let $\calX$ be a derived stack over a base
stack $S$. Then $(\calX/S)_{\sm}$ (resp. $(\calX/S)_{\et}$,
$(\calX/S)_{\Zar}$) will denote the
small site of $S$-morphisms $f: \calU \to \calX$ 
such that $f$ is schematic, bounded, and smooth (resp. \'etale, Zariski); 
covers are defined as usual (i.e.\ surjectivity on 
geometric points). Note that any morphism 
between objects of $(\calX/S)_{\sm}$ (resp.
$(\calX/S)_{\et}$, $(\calX/S)_{\Zar}$) is 
of bounded Tor dimension (resp. \'etale, an open embedding). 
By Nisnevich descent, we mean descent for the Grothendieck
topology generated by Nisnevich distinguished squares.
\edefn

In the present work, we are primarily interested in the
Zariski topology and in the case
when either $S = \Spec k$ or $S = BG$ for an affine algebraic group $G$.

\brem
When $S = \Spec k$, we denote $(\calX/S)_{\et}$ and $(\calX/S)_{\Zar}$
simply by $\calX_{\et}$ and $\calX_{\Zar}$, in which case
\cref{def: site} reduces to \cite[Definition A.2.1]{PreygelT}.
\erem

\begin{exmp}
Let $X$ be a scheme with an action of an
affine algebraic group $G$. Then the objects of
$(X/G)_{\et}$ are $G$-equivariant \'etale maps of schemes $U \to X$,
and the objects of $(X/G)_{\Zar}$ are open embeddings of $G$-invariant
subschemes of $X$. Covers in the former are \'etale covers of schemes
where the maps are all $G$-equivariant, and covers in the latter are
simply covers by $G$-invariant opens.
\end{exmp}

\blem
Let $\calX = X/G$ for the action of an affine algebraic
group $G$ on a scheme $X/k$. Then
$(X/G)_{\Zar}$ and $([X/G]/BG)_{\Zar}$ are equivalent
sites.
\elem

\bproof
It suffices to show that for every schematic open embedding
of prestacks over $BG$ of the form $\calY \to U/G$ for some $G$-scheme
$U$, $\calY$ has the form $Z/G$ for another $G$-scheme $Z/k$.

Note that the property of being an open embedding
is preserved under pullback,
as is the property of being a $G$-torsor,
so we have the following Cartesian diagram,
\[\begin{tikzcd}
	{U \times_{U/G} \calY} & U \\
	\calY & {U/G}
	\arrow["{G\textrm{-torsor}}"', from=1-1, to=2-1]
	\arrow["{\textrm{open}}", hook, from=2-1, to=2-2]
	\arrow["{\textrm{open}}", hook, from=1-1, to=1-2]
	\arrow[from=1-2, to=2-2]
	\arrow["\lrcorner"{anchor=center, pos=0.125}, draw=none, from=1-1, to=2-2]
\end{tikzcd}\]
From the diagram, is follows that $\calY \simeq (U \times_{U/G} \calY)/G$.
On the other hand, $U \times_{U/G} \calY$ is a $k$-scheme because it
has an open embedding into the $k$-scheme $U$. Thus $\calY$ has the
form of stack quotient of a $G$-scheme over $k$.
\eproof

\subsubsection{}
\label{sssec: varpi}
If $X$ is a $G$-scheme, 
there is an obvious comparison map of
sites, $\varpi: X_{\et} \to (X/G)_{\et}$
(resp. $\varpi: X_{\Zar} \to (X/G)_{\Zar}$).
In particular, there is a pushforward functor for
sheaves with coefficients in an arbitrary $\infty$-category
$\scrC$,
\begin{align*}
\varpi_*: \Shv(X_{\et}; \scrC) &\to \Shv((X/G)_{\et}; \scrC) \\
\varpi_*: \Shv(X_{\Zar}; \scrC) &\to \Shv((X/G)_{\Zar}; \scrC)
\end{align*} 
In later sections, we will use the latter functor
in order to compare the sheaves of Hochschild invariants
of graded and ungraded matrix factorizations for a graded
Landau--Ginzburg pair.


\subsection{Ind-coherent sheaves}
\label{ssec: ind-coherent sheaves}

The following theorem is a variant over a
more general base of the theorem stated 
and proven in appendix of \cite{PreygelT}.
The proof found therein is equally valid, \textit{mutatis mutandis},
over a more general base stack $S$.

\bthm[{c.f. \cite[Theorem A.2.5]{PreygelT}}]
\label{thm: IndCoh satisfies descent}
Suppose that $\calX/S$ is a Noetherian, perfect derived
$S$-stack that has affine diagonal. Let $\IndCoh(-/S)$ denote
the presheaf of $\QCoh(S)$-module categories,
\[U \mapsto \IndCoh(U), \hspace{1cm} [f:U' \hook U] \mapsto f^!\] 
Then,
\begin{enumerate}[label=(\roman*), ref=\roman*]
\item $\IndCoh(-/S)$ satisfies Nisnevich descent, and finite \'etale descent over $S$.

\item $\IndCoh(-/S)$ satisfies representable \'etale descent over $S$.

\item if furthermore $\calX$ is Deligne-Mumford over $S$, $\IndCoh(-/S)$ satisfies smooth descent over $S$.

\item if furthermore $\calX$ is Deligne-Mumford almost of finite presentation over $S$,
$\IndCoh(-/S)$ satisfies smooth descent over $S$, and $\IndCoh(\calX/S)$ coincides with  $\QC^!(\calX/S)$.
\end{enumerate}
\ethm

\brem
In particular, if $\calX$ is a \labelcref{stack condition}
stack, almost of finite presentation over $S$, $\IndCoh(-/S)$
is a $\QCoh(S)\mod$-valued sheaf on $\calX_{\sm}$, $\calX_{\et}$,
and $\calX_{\Zar}$.
\erem

\begin{notation}
When the base $S$ is clear from context,
we prefer to write $\IndCoh(\calX)$ instead
of $\IndCoh(\calX/S)$.
\end{notation}

\begin{notation}
To make the distinction between
the category $\IndCoh(\calX)$ and the
corresponding sheaf on $\calX$ (in whichever
topology), we frequently denote the
latter by $\un{\IndCoh}(\calX)$.
\end{notation}


\subsection{Matrix factorizations}

\Cref{thm: IndCoh satisfies descent} allows
one easily to show that $\PreMF(M/S,\scrL,f)$,
$\PreMF^{\infty}(M/S,\scrL,f)$ and $\MF^{\infty}(M/S,\scrL,f)$,
satisfy relative \'etale descent on
$M$, as we do in the following proposition.

\bprop[{c.f. \cite[Proposition A.3.1]{PreygelT}}]
\label{prop: relative MF satisfies descent}
Suppose $(M/S,\scrL, f)$ is a relative Landau--Ginzburg
pair over $(S,\scrL)$. Then,
\begin{enumerate}
\item \label{item: PreMF is sheaf}
The assignments,
\[U \mapsto \PreMF(U/S,\scrL, f|_U) \hspace{1cm} \textrm{and} \hspace{1cm} U \mapsto \PreMF^{\infty}(U/S,\scrL, f|_U),\]
determine sheaves of $\Coh(\G)$-linear and $\IndCoh(\G)$-linear $\infty$-categories, 
respectively, on $(M/S)_{\et}$ (hence $(M/S)_{\Zar}$).

\item \label{item: MF is sheaf}
The assignment,
\[U \mapsto \MF^{\infty}(U/S,\scrL,f|_U),\]
determines a sheaf of $\Sing^{\infty}(\G)$-linear 
$\infty$-categories on $(M/S)_{\et}$ (hence $(M/S)_{\Zar}$).
\end{enumerate}
\eprop

\bproof
Since an \'etale morphism $s: U \to M$
is in particular smooth, the composition
$U \to M \to S$ is a smooth $S$-scheme, 
so $(U, f \circ s)$ is a Landau--Ginzburg pair over $(S,\scrL)$.

We first prove (\labelcref{item: PreMF is sheaf}).
Note that any \'etale cover of $M_{\et}$ restricts to 
an \'etale cover of $M_{0_{\scrL}}$, and that a $k[\b]$-linear presheaf
is a sheaf if and only if it is a sheaf forgetting the extra linear 
structure. It therefore suffices to show that $\Coh$ and $\IndCoh$ 
are sheaves on $(M_{0_{\scrL}})_{\et}$; the former follows from 
the analogous theorem for $\QCoh$ and the local definition 
of $\Coh$ (since $X_{0_{\scrL}}$ is coherent), and the latter 
follows from the analogous theorem for $\IndCoh$ (\cite[Theorem A.2.5]{PreygelT}).

To prove part (\labelcref{item: MF is sheaf}), it suffices to note 
that $- \otimes_{\IndCoh(\G)} \Sing^{\infty}(\G): \IndCoh(\G)\mod
\to \Sing^{\infty}(\G)\mod$ commutes with limits,
since $\Sing^{\infty}(\G)$ is compactly generated by
the small $\Coh(\G)$-module category $\Sing(\G)$ and therefore 
a dualizable object of $\IndCoh(\G)\mod$ by \cite[Lemma 4.2.1]{PreygelT}.
\eproof

\begin{notation}
To make the distinction between
the categories $\PreMF^{\infty}(M/S,\scrL,f)$ and
$\MF^{\infty}(M/S,\scrL,f)$ and the
corresponding sheaves, we denote the
latter by $\unPreMF^{\infty}(M/S,\scrL,f)$
and $\unMF^{\infty}(M/S,\scrL,f)$, respectively.
\end{notation}

\brem
Let $Z \subset M_{0_{\scrL}}$ be a closed
substack. Unwinding the definitions, we have that
\begin{align*}
\Gamma_Z(\unPreMF^{\infty}(M/S,\scrL,f)) 			&\simeq \PreMF^{\infty}_Z(M/S,\scrL,f) \\
\Gamma_Z(\unMF^{\infty}(M/S,\scrL,f)) 			&\simeq \MF^{\infty}_Z(M/S,\scrL,f),
\end{align*}
where $\PreMF^{\infty}_Z(M/S,\scrL,f)$ and $\MF^{\infty}_Z(M/S,\scrL,f)$
were defined in \cref{def: PreMF with supports}
and in \cref{def: MF with supports}, respectively.
\erem

\subsection{Hochschild invariants}
In general, the pointwise Hochschild homology of a
sheaf of dualizable categories is not itself a sheaf,
such as seen in the following example.

\begin{exmp}
Let $X$ be a perfect stack over $k$. Then $\HH_{\bullet}^k(\QCoh(X))
\simeq \O(\calL X)$, but $\calL X$ is not local on $X$
with respect to the smooth (or even \'etale) topology. That is,
one cannot expect to recover loops in $X$ from loops in a smooth
cover of $X$.
For example, suppose $X = BG$. Then the map $\ast \to BG$
is a smooth cover, but it is clear that $\O(\calL(BG))
\simeq \O(G/G)$ cannot be recovered from $\O(\calL(\ast)) 
\simeq \O(\ast) \simeq k$.
\end{exmp}

Whenever this happens, we simply sheafify the presheaf
obtained by the pointwise application of Hochschild homology.

\begin{conv}
\label{conv: pointwise presheaf}
Let $\F$ be a sheaf valued in a fixed
symmetric monoidal category $\scrC$ such
that the sections of $\F$ over any open are 
dualizable objects of $\scrC$.

If $F: \scrC \to \scrD$ is any functor which does not
preserve limits, we write $F(\F)$ to denote the presheaf obtained by
the pointwise application of $F$ to the sections of $\F$---that is,
the presheaf determined by the assignment
\[U \mapsto F(\F(U))\]
on the level of objects.

We will denote the sheafification of $F(\F)$ by $\un{F}(\F)$.
Alternatively, if the notation for $\F$ itself includes an underline,
such as $\unPreMF^{\infty}$ or $\unMF^{\infty}$, then we
will denote sheafification of $F(\F)$ by simultaneously underlining $F$
and omitting the underline from the notation for $\F$.

For example, $\HH_{\bullet}^{\kb}(\unMF^{\infty}(M,f))$
will denote the $\kb\mod$-valued presheaf on either
$M_{\et}$ or $M_{\Zar}$ obtained by the
pointwise application of the $\kb$-linear Hochschild homology
functor to the sections of $\unMF^{\infty}(M,f)$,
and $\unHH_{\bullet}^{\kb}(\MF^{\infty}(M,f))$ will
denote its sheafification.
\end{conv}

\subsubsection{Descent for Hochschild homology of $\PreMF$}
Though we will not need this fact, one can show
that the presheaf $\HH_{\bullet}^{k[\b]}(\unPreMF(M,f))$
is a sheaf on the \'etale (hence Zariski) site of $M$.

\bprop
\label{prop: HH of PreMF is a sheaf}
Let $(M,f)$ be a Landau--Ginzburg pair. Then 
$\HH_{\bullet}^{k[\b]}\unPreMF^{\infty}(M,f)$ 
is a $k[\b]\mod$-valued sheaf on $M_{\et}$.
\eprop

\bproof
Note that by \cref{prop: natural BG_a maps are equivalences}
below, the natural comparison map,
\[\HH_{\bullet}^{k[\b]}(\PreMF^{\infty}(M,f)) \to \HH_{\bullet}^k(\IndCoh(M))^{B\wih{\G}_a}\]
is an equivalence. 
Since the forgetful functor $\oblv: B\wih{\G}_a\mod \to \Vect_k$
is monadic, it creates limits. It follows that a presheaf $\F$ valued
in $B\wih{\G}_a\mod$ is a sheaf if and only if the $\Vect_k$-valued
presheaf $\oblv(\F)$ is a sheaf. Since $\HH_{\bullet}^k(\IndCoh(M))$
is a sheaf on $M_{\et}$ by \cite[Theorem A.1.2.5]{PreygelThesis}
(or \cref{prop: relative MF satisfies descent} above), 
we therefore have that $\HH_{\bullet}^{k[\b]}(\unPreMF^{\infty}(M,f))$ 
is a sheaf on $M_{\et}$.
\eproof

Similarly,  $\unPreMF^{\infty}(M/\G_m,f)$ also gives a sheaf.

\bprop
\label{prop: HH of graded PreMF is a sheaf}
Let $(M/\G_m,f)$ be a relative Landau--Ginzburg pair
over $(B\G_m, \O_{B\G_m}(d))$. Then
$\HH_{\bullet}^{k[\b]_{\gr}}\unPreMF^{\infty}(M/\G_m,f)$ 
is a $k[\b]_{\gr}\mod$-valued sheaf on $(M/\G_m)_{\et}$.
\eprop

\bproof
The proof of \cref{prop: HH of PreMF is a sheaf}
adapts easily to the graded setting.
\eproof

%



\section{Ind-coherent sheaves and graded matrix factorizations}
\label{sec: ind-coherent sheaves and graded matrix factorizations}

It is well-known that the category of coherent sheaves
on a global complete intersection is equivalent to the
category of graded matrix factorizations on a certain 
associated graded Landau--Ginzburg pair. The connection
between coherent sheaves and graded matrix factorizations
was first shown by Isik in \cite{Isik}, and, 
since then, variations on this theme have 
appeared in works such as \cite{Shipman, Hirano, OR}.
It is claimed without proof in \cite[\S H.1.7]{AG15}
that Isik's equivalence of categories exchanges the singular
support of coherent sheaves with the set-theoretic
support of graded matrix factorizations by
intertwining the module structures on each category.
The author of \cite{Toda} produces an equivalence\footnotemark
\footnotetext{It is not immediately clear what the relation is between
Toda's equivalence and the others mentioned.}
between coherent sheaves and graded matrix factorizations
in \cite[Theorem 2.3.3]{Toda}, and then proves
in \cite[Proposition 2.3.9]{Toda} that
this equivalence exchanges singular support with
support.

In this section, we produce an equivalence, in our present
framework, between $\IndCoh(\Z)$ for $\Z$ given as the derived
zero fiber of a map $f: X \to V$, and the graded
matrix factorizations of a graded Landau--Ginzburg
pair denoted $(X \times V^{\vee}, \wit{f})$. 
Here the advantage of our framework
for matrix factorizations manifests itself: the
equivalence we produce is easily shown to sheafify
over the \'etale topology of $\Z$, and to exchange
the singular support of coherent sheaves 
for set-theoretic support of matrix factorizations.

\subsubsection{Conventions}
Throughout this section, we work
over a fixed field, $k$, of characteristic zero.

\begin{con}
\label{con: dimensional reduction}
Let $f: X \to V$ be a
map from a smooth scheme, $X$, to a
finite dimensional vector space, $V$. Let
$\wit{f}: X \times V^{\vee} \to \Aone$
denote the regular function given by
$(x,\lambda) \mapsto \langle \lambda, f(x) \rangle$.
Consider $X \times V^{\vee}$ as
a $\G_m$-scheme via the weight $2$
scaling action on $V^{\vee}$, and
$\Aone$ as a $\G_m$-scheme
also with the weight $2$ scaling action.
Clearly, $\wit{f}$ is equivariant with
respect to these actions, so
$(X \times V^{\vee},\wit{f})$
comprises a graded Landau--Ginzburg pair
of weight $2$. As described in
\S\labelcref{sssec: graded Landau--Ginzburg pairs},
this graded Landau--Ginzburg pair
is equivalent data to the relative Landau--Ginzburg pair
over $(B\G_m, \O_{B\G_m}(2))$ given by
\beqn
\label{eqn: main relative LG pair}
((X \times V^{\vee})/\G_m, \wit{f}).
\eeqn
\end{con}

\begin{notation}
In the sequel, we denote the derived
zero fiber of $f$ by $\Z$; the classical
zero fiber of $f$ by $Z$; and the
derived zero fiber of $\wit{f}$ by $\wit{\Z}$.
\end{notation}

\begin{convention}
$\QCoh(-)$ and $A\mod$ for commutative $A$ 
are always given the standard $\otimes$-symmetric
monoidal structures. For $\calG$ a group scheme, 
$\IndCoh(\calG)$ is given the convolution monoidal
structure.
\end{convention}

\subsection{Sheafifying Isik's theorem}
The purpose of this section is to prove the
following result.

\bthm
\label{thm: fancy Isik's theorem}
Let $p: \Z \times V^{\vee}/\G_m \to \Z$ denote the
standard projection, and let $i: \Z \times V^{\vee}/\G_m
\to \wit{\Z}/\G_m$ denote the natural inclusion. Then
the map,
\[i_*^{\IndCoh}p^!: \IndCoh(\Z) \xrightarrow{\simeq} \Sing^{\infty}(\wit{\Z}/\G_m),\]
is an equivalence of $\QCoh(X \times V^{\vee}/\G_m)$-modules.
\ethm

\bproof
By \cref{prop: showing equivariance}, the functor $i_*^{\IndCoh}p^!:
\IndCoh(\Z) \to \Sing^{\infty}(\wit{\Z}/\G_m)$ is $\QCoh(X \times V^{\vee}/\G_m)$-linear.
Since the forgetful functor from $\QCoh(X \times V^{\vee}/\G_m)$-linear
categories to plain dg-categories is conservative, it suffices to show
that $i_*^{\IndCoh}p^!$ is an equivalence on underlying dg-categories.
Since $p$ is smooth, $i_*^{\IndCoh}p^!$ and $i_*^{\IndCoh}p^{\IndCoh, *}$
differ by a factor of a line bundle, so each is an equivalence if and only if
the other is. It therefore follows from \cref{prop: showing equivalence}
that $i_*^{\IndCoh}p^!$ is an equivalence, and we are done.
\eproof

The $\QCoh(X \times V^{\vee}/\G_m)$-module structures
referenced in the theorem are defined in detail below: the module
structure on $\IndCoh(\Z)$ is defined in the Arinkin--Gaitsgory 
paper of the singular support of ind-coherent sheaves (\cite{AG15})
for a global complete intersection, and the module
structure on $\IndCoh(\wit{Z}/\G_m)$ corresponds to
the inclusion $\wit{\Z}/\G_m \hook X \times V^{\vee}/\G_m$.

\subsubsection{$1$-affineness}
In the foundational paper \cite{Gaitsgory1affineness},
Gaitsgory constructs an adjunction,
\[\begin{tikzcd}
	{\lLoc_{\calY}: \QCoh(\calY)\mod} && {\ShvCat(\calY): \Gamma_{\calY}}
	\arrow[shift right=1, from=1-1, to=1-3]
	\arrow[shift right=1, from=1-3, to=1-1],
\end{tikzcd}\]
for an arbitrary prestack $\calY$. Objects of $\ShvCat(\calY)$
are quasicoherent sheaves of dg-categories on $\calY$.\footnotemark 
\footnotetext{Briefly, the functor $\calY \mapsto \ShvCat(\calY)$ is given
by Kan extending the functor on affine derived schemes that sends
an affine $S$ to $\QCoh(S)\mod$.
See \cite[\S 1.1]{Gaitsgory1affineness} for more details.}
Prestacks for which the above adjunction is actually
an adjoint equivalence are called \emph{$1$-affine}.
Using the $1$-affineness of $X \times V^{\vee}/\G_m$,
obtain the following corollary of \cref{thm: fancy Isik's theorem}.

\bcor
\label{cor: sheafy Isik's theorem}
There is an equivalence,
\[\un{\IndCoh}(\Z) \xrightarrow{\simeq} \unMF^{\infty}(X \times V^{\vee}/\G_m, \wit{f}),\]
of objects in $\ShvCat(X \times V^{\vee}/\G_m)$.
\ecor

\bproof
Using the equivalence $\MF^{\infty}(X \times V^{\vee}/\G_m) \simeq
\Sing^{\infty}(\wit{\Z}/\G_m)$, we view $i_*^{\IndCoh}p^!$ as an equivalence
of $\QCoh(X \times V^{\vee}/\G_m)$-modules $\IndCoh(\Z) \xrightarrow{\simeq}
\MF^{\infty}(X \times V^{\vee}/\G_m, \wit{f})$. We conclude by noting that
$X \times V^{\vee}/\G_m$ is $1$-affine by \cite[Theorem 2.2.4]{Gaitsgory1affineness}.
\eproof

The sections of the sheaf $\un{\IndCoh}(\Z)$ with support
in a Zariski closed subset $Y \subset X \times V^{\vee}/\G_m$
is precisely the subcategory $\IndCoh_Y(\Z)$ of ind-coherent
sheaf with singular support contained in $Y$. On the other hand,
a section of $\unMF^{\infty}(X \times V^{\vee}/\G_m, \wit{f})$
has a representative $\F \in \IndCoh(\wit{\Z}/\G_m)$, and the support
of this section is exactly the support of $\F$ as an ind-coherent sheaf.
Thus, we see that the equivalence of sheaves 
$i_*^{\IndCoh}p^!$ exchanges the Arinkin--Gaitsgory 
singular support for the usual support of sheaves.

\subsubsection{}
In the remainder of this section, we collect the results
used in the proof of \cref{thm: fancy Isik's theorem}.

\subsection{$\QCoh(X \times V^{\vee}/\G_m)$-linearity}

\subsubsection{}
As explained earlier in \cref{construction: singular support QCoh action}, 
$\IndCoh(\Z)$ has a $\QCoh(V^{\vee}/\G_m)$-module
structure corresponding under Koszul duality to the action of based loops in $V$.
From now on, any reference to a $\QCoh(V^{\vee}/\G_m)$-module
structure on $\IndCoh(\Z)$ will refer to that one.
\begin{enumerate}
\item 
\label{eyetem 1}
Consider $\IndCoh(\Z)$ as a $\QCoh(V^{\vee}/\G_m) \otimes \QCoh(V^{\vee}/\G_m)$-module
structure by restricting its $\QCoh(V^{\vee}/\G_m)$-module
structure along the multiplication map $\QCoh(V^{\vee}/\G_m)
\otimes \QCoh(V^{\vee}/\G_m) \xrightarrow{\otimes}
\QCoh(V^{\vee}/\G_m)$ (i.e. pullback along the
diagonal embedding of $V^{\vee}/\G_m$).
\item 
\label{eyetem 2}
Consider $\IndCoh(\Z \times V^{\vee}/\G_m) \simeq \IndCoh(\Z) \otimes \IndCoh(V^{\vee}/\G_m)$
as a $\QCoh(V^{\vee}/\G_m) \otimes \QCoh(V^{\vee}/\G_m)$-module via the
obvious $\QCoh(V^{\vee}/\G_m)$-actions on each of the factors.
\end{enumerate}

We now describe the 
$\QCoh(V^{\vee}/\G_m) \otimes \QCoh(V^{\vee}/\G_m)$-module 
structure with which we equip $\IndCoh(\wit{\Z}/\G_m)$
for the remainder of the paper.

We first establish some notation. 

\begin{notation}
\label{notation: buncha functions}
Let $F: V \times V^{\vee} \to 
\Aone \times V^{\vee}$ denote the map
$(v, \lambda) \mapsto (\lambda(v), \lambda)$, 
and let $\Omega(F): \Omega_0 V \times V^{\vee} 
\to \Omega_0 \Aone \times V^{\vee}$ 
denote the map induced by $F$ 
on the based loop spaces. Note that
$\Omega(F)$ is map of group derived
schemes over $V^{\vee}$. Let $\pi: \Omega_0 V 
\times V^{\vee} \to \Omega_0 V$ 
be the standard projection onto the first factor. 
\end{notation}

\begin{enumerate}[resume]
\item
\label{eyetem 3}
Consider $\wit{f}$ as a map $X \times V^{\vee}/\G_m
\to \O_{B\G_m}(2)$, as in \cref{con: dimensional reduction}. 
Since $\wit{\Z}/\G_m$ is the zero locus of this map, 
$\IndCoh(\wit{\Z}/\G_m)$ obtains an action 
of $\IndCoh(\Omega_0 \Aone/\G_m)$ via based loops acting on $\wit{\Z}/\G_m$.
This action extends to
an action of $\IndCoh(\Omega_0 \Aone/\G_m) \otimes \QCoh(V^{\vee}/\G_m)$,
where the action of the second factor corresponds
to the projection $\wit{\Z}/\G_m \to V^{\vee}/\G_m$.
Since $\Omega(F)$ is a morphism of
group derived stacks over $V^{\vee}/\G_m$,
the pushforward $\Omega(F)_*^{\IndCoh}$ 
is uniquely symmetric
monoidal for the convolution monoidal
structures on the source and target (see e.g. \cite[\S F.4.8]{AG15}).
Restriction along $\Omega(F)_*^{\IndCoh}$ gives
$\IndCoh(\wit{\Z}/\G_m)$ a $\IndCoh(\Omega_0 V/\G_m) \otimes \QCoh(V^{\vee}/\G_m)$-module
structure. Under Koszul duality, however, $\IndCoh(\Omega_0 V/\G_m) \simeq \QCoh(V^{\vee}/\G_m)$,
and by this equivalence $\IndCoh(\wit{\Z}/\G_m)$ is a module over $\QCoh(V^{\vee}/\G_m) \otimes
\QCoh(V^{\vee}/\G_m)$.
\end{enumerate}

\subsubsection{}
\label{sssec: two lemmas}
The following two lemmas imply that $i_*^{\IndCoh}p^!$
is a map of $\QCoh(V^{\vee}/\G_m) \otimes \QCoh(V^{\vee}/\G_m)$-modules.

\blem
\label{lem: p^! intertwines}
$p^!: \IndCoh(\Z) \to \IndCoh(\Z \times V^{\vee}/\G_m)$ is
$\QCoh(V^{\vee}/\G_m) \otimes \QCoh(V^{\vee}/\G_m)$-linear
with respect to module structures 
(\labelcref{eyetem 1})
and (\labelcref{eyetem 2}).
\elem

\bproof
We will prove a more general fact that will imply the lemma.

Let $\calX$ and $\calU$ be derived stacks over $k$, where $\calU$ is
smooth, and suppose that $\IndCoh(\calX)$ has a $\QCoh(\calU)$-module
structure. Let $\pr_1: \calX \times \calU \to \calX$ denote the standard
projection onto the first factor. We will show that 
\[\pr_1^!: \IndCoh(\calX) \to \IndCoh(\calX \times \calU) \simeq \IndCoh(\calX) \otimes_k \QCoh(\calU)\]
has the canonical structure of a map of $\QCoh(\calU)^{\otimes 2}$-modules, where
\begin{itemize}
\item[$-$] the module structure on the source is given by restricting the $\QCoh(\calU)$-module
structure along the multiplication map $\mult: \QCoh(\calU) \otimes_k \QCoh(\calU) \to \QCoh(\calU)$,
\item[$-$] and the one of the target is given by evident actions of $\QCoh(\calU)$ on each
of the two factors.
\end{itemize}

The module structures on both the source and target of $\pr_1^!$ are
given each by a simplicial diagram of functors,
\beqn
\label{eqn: simp 1}
\begin{tikzcd}
	\cdots & {\left(\QCoh(\calU)^{\otimes 2}\right)^{\otimes 2} \otimes \IndCoh(\calX)} & {\QCoh(\calU)^{\otimes 2} \otimes \IndCoh(\calX)} & {\IndCoh(\calX)}
	\arrow[shift right, from=1-3, to=1-4]
	\arrow[shift left, from=1-3, to=1-4]
	\arrow[from=1-4, to=1-3]
	\arrow[shift right=2, from=1-2, to=1-3]
	\arrow[shift left=2, from=1-2, to=1-3]
	\arrow[from=1-2, to=1-3]
	\arrow[shift right, from=1-3, to=1-2]
	\arrow[shift left, from=1-3, to=1-2]
	\arrow[shift right, from=1-1, to=1-2]
	\arrow[shift right=3, from=1-1, to=1-2]
	\arrow[shift left=3, from=1-1, to=1-2]
	\arrow[shift left, from=1-1, to=1-2]
	\arrow[shift right=2, from=1-2, to=1-1]
	\arrow[shift left=2, from=1-2, to=1-1]
	\arrow[from=1-2, to=1-1]
\end{tikzcd}
\eeqn
\beqn
\label{eqn: simp 2}
\begin{tikzcd}
	\cdots & {\QCoh(\calU)^{\otimes 2} \otimes \left(\IndCoh(\calX) \otimes \QCoh(\calU)\right)} & {\IndCoh(\calX) \otimes \QCoh(\calU)}
	\arrow[shift right, from=1-2, to=1-3]
	\arrow[shift left, from=1-2, to=1-3]
	\arrow[from=1-3, to=1-2]
	\arrow[shift right=2, from=1-1, to=1-2]
	\arrow[shift left=2, from=1-1, to=1-2]
	\arrow[from=1-1, to=1-2]
	\arrow[shift right, from=1-2, to=1-1]
	\arrow[shift left, from=1-2, to=1-1]
\end{tikzcd}
\eeqn

In principle, giving $\pr_1^!$ the structure of a map of modules
amounts to specifying a map of simplicial objects between \labelcref{eqn: simp 1}
and \labelcref{eqn: simp 2} such that $\pr_1^!$ is the map on $0$-simplices.
In practice, such date is hard write in any elegant way, so we will content
ourselves will writing down the component of such a map on the level of $1$-simplices,
and see that similar maps may be defined on the level of higher simplices.

The map on $1$-simplices is given as follows.
\[\begin{tikzcd}
	{\QCoh(\calU)^{\otimes 2} \otimes \IndCoh(\calX)} & {\IndCoh(\calX)} \\
	{\QCoh(\calU)^{\otimes 2} \otimes (\IndCoh(\calX) \otimes \QCoh(\calU))} & {\IndCoh(\calX) \otimes \QCoh(\calU)}
	\arrow["\proj"', shift right, from=1-1, to=1-2]
	\arrow["\act", shift left, from=2-1, to=2-2]
	\arrow["{(\id_{\QCoh(\calU)} \otimes \pr_1^!) \circ (\mult \otimes \id_{\IndCoh(\calX)})}"', from=1-1, to=2-1]
	\arrow["{\pr_1^!}", from=1-2, to=2-2]
	\arrow["\act", shift left, from=1-1, to=1-2]
	\arrow["\proj"', shift right, from=2-1, to=2-2]
	\arrow[from=1-2, to=1-1]
	\arrow[from=2-2, to=2-1].
\end{tikzcd}\]
It is evident that the diagram is canonically commutative
once we note that $\pr_1^! \simeq \id_{\IndCoh(\calX)} \otimes \unit_{\QCoh(\calU)}$
and unravel the definitions on the $\proj$ and $\act$ functors.
\eproof

\blem
Let $i: \Z \times V^{\vee}/\G_m \to \wit{\Z}/\G_m$
denote the obvious inclusion.
Then $i_*^{\IndCoh}$ is $\QCoh(V^{\vee}/\G_m) 
\otimes \QCoh(V^{\vee}/\G_m)$-linear with
respect to the module structures (\labelcref{eyetem 2})
and (\labelcref{eyetem 3}).
\elem

\bproof
We begin with a toy model in the ungraded setting.

Take $T : W_1 \to W_2 \in \Vect_k^{\heart}$ a linear 
map of finite dimensional vector spaces. Suppose we
have $q : Q \to W_1$, a map of schemes over $k$. Then $\Omega_0 W_1$ acts on $q ^{-1}(0)$ 
while $\Omega_0 W_2$ acts on $(T \circ q )^{-1}(0)$. The natural
map $q^{-1}(0) \to (T \circ q)^{-1}(0)$ is then $\Omega_0 W_1$-equivariant, 
where $\Omega_0 W_1$ acts on the target via the
map $\Omega(T) : \Omega_0 W_1 \to \Omega_0 W_2$.


We can repeat the story above with a base stack $S$ in place of $\Spec k$. 
Here $W_i$ are replaced by finite rank vector bundles on $S$, $T$ 
is replaced by a map of vector bundles, $Q$ is an $S$-scheme, and so on. 
Note that $\Omega_0 W_1$ and $\Omega_0 W_2$ are replaced by suitable group
derived schemes over $S$. Now take $W_1 = V \times S$, $W_2 = \Aone \times S$,
$Q = X \times S$, and $q = f \times \id_S$. Then the map $T$ is determined
by a map $s: S \to V^{\vee}$.
Let $\wit{Z}_s := (T \circ q)^{-1}(0)$. 
Then, similar to above, the natural map $i: \Z \times S \to \wit{Z}_s$
is $\Omega_0 V \times S$-equivariant, where $\Omega_0 V \times S$
acts on the target via the map $\Omega(T): \Omega_0 V \times S 
\to \Omega_0 \Aone \times S$.
Note that this immediately implies that $i_*^{\IndCoh}$ is
a map of $\IndCoh(\Omega_0 V) \otimes \QCoh(S)$-module
categories, for the module structures induced by the
group derived $S$-scheme actions described above.

The story above has an obvious counterpart in the $\G_m$-equivariant setting.
Since $T: W_1 \to W_2$ is linear, it is in particular
$\G_m$-equivariant, so it descends to a map
$T: W_1/\G_m \to W_2/\G_m$, where we abuse notation
by also denoting the map $T$. We also obtain a map
$\Omega(T): \Omega_0 W_1/\G_m \to \Omega_0 W_2/\G_m$
for the same reason. Taking $W_1 = V \times S$, $W_2 = \Aone \times S$, etc.
to be the same as in the previous paragraph, we find that the map
$\Z \times S \to \wit{\Z}_s$ intertwines $\Omega_0 V/\G_m \times S$-actions
along the map of group $S$-schemes, $\Omega(T): \Omega_0 V/\G_m \times S \to 
\Omega_0 \Aone/\G_m \times S$. As before, this immediately
implies that $i_*^{\IndCoh}$ is a map of $\IndCoh(\Omega_0 V/\G_m) \otimes
\QCoh(S)$-modules. Note that the map $T$ in this setting
is determined by a map $s: S \to V^{\vee}/\G_m$.

Finally, taking $S = V^{\vee}/\G_m$, $s = 
\id_{V^{\vee}/\G_m}: V^{\vee}/\G_m \to V^{\vee}/\G_m$, 
and unwinding yields the claim (compare $T$ to the map
$F$ from \cref{notation: buncha functions}).
\eproof

\subsubsection{}
Observe that any $\QCoh(V^{\vee}/\G_m) \otimes \QCoh(V^{\vee}/\G_m)$-module
has a $\QCoh(V^{\vee})$-module structure given by restriction
along the pullback,
\[\pr_2^*: \QCoh(V^{\vee}/\G_m) \to \QCoh(V^{\vee}/\G_m \times V^{\vee}/\G_m),\]
where $\pr_2: V^{\vee}/\G_m \times V^{\vee}/\G_m \to V^{\vee}/\G_m$
denotes the standard projection onto the first factor.\footnotemark
\footnotetext{Alternatively, given any $k$-linear symmetric monoidal
categories $\scrC$ and $\scrD$, there are natural symmetric monoidal
$k$-linear functors $\scrC \to \scrC \otimes \scrD$ and $\scrD \to \scrC \otimes \scrD$.
The map $\pr_2^*$ is the latter of these functors for $\scrC
= \scrD = \QCoh(V^{\vee}/\G_m)$.}

\brem
\label{rem: QCoh action 1}
Note that the $\QCoh(V^{\vee}/\G_m)$-module
structure on $\IndCoh(\Z)$ obtained from (\labelcref{eyetem 1})
in this way is the standard $\QCoh(V^{\vee}/\G_m)$-module structure
on $\IndCoh(\Z)$ discussed in \cref{construction: singular support QCoh action}.
Indeed, this $\QCoh(V^{\vee}/\G_m)$-module structure is given by restricting the standard
one along the map $\pr_2^*\Delta^*: \QCoh(V^{\vee}/\G_m) \to \QCoh(V^{\vee}/\G_m)$,
but $\pr_2^*\Delta^* \simeq (\pr_2 \circ \Delta)^* = \id_{V^{\vee}/\G_m}^* = \id_{\QCoh(V^{\vee}/\G_m)}$.
\erem

\brem
\label{rem: QCoh action 2}
On the other hand, the $\QCoh(V^{\vee}/\G_m)$-module
structure on $\IndCoh(\wit{\Z}/\G_m)$ obtained from (\labelcref{eyetem 3})
in this way is evidently one corresponding to the projection $\wit{\Z}/\G_m
\to V^{\vee}/\G_m$.
\erem

These observations in conjunction with the lemmas
of the previous section (\S \labelcref{sssec: two lemmas})
immediately imply the following corollary.

\bcor
\label{cor: equi 1}
$i_*^{\IndCoh}p^!$ is $\QCoh(V^{\vee}/\G_m)$-linear
for the actions described in \cref{rem: QCoh action 1}
and \cref{rem: QCoh action 2}.
\ecor

\subsubsection{}
On the other hand, $\IndCoh(\Z)$ and $\IndCoh(\wit{\Z}/\G_m)$
also have natural $\QCoh(X)$-actions.

\begin{enumerate}[label=(\roman*), ref=\roman*]
\item 
\label{enum 1}
$\IndCoh(\Z)$ obtains a natural $\QCoh(X)$-module
structure via the pullback for $\QCoh$ along the inclusion
$\Z \hook X$ and the canonical action of $\QCoh(\Z)$ on $\IndCoh(\Z)$.

\item 
\label{enum 2}
$\IndCoh(\wit{\Z}/\G_m)$ obtains a $\QCoh(X)$-module
structure similarly via the pullback for $\QCoh$ along the projection map
$\wit{\Z}/\G_m \to X$ and the canonical action of $\QCoh(\wit{\Z}/\G_m)$
on $\IndCoh(\wit{\Z}/\G_m)$.
\end{enumerate}

\blem
\label{lem: equi 2}
$i_*^{\IndCoh}p^!$ is $\QCoh(X)$-linear
for the $\QCoh(X)$-actions (\labelcref{enum 1}) and
(\labelcref{enum 2}) above.
\elem

\bproof
In addition to the $\QCoh(X)$-actions
considered above, we also consider the
$\QCoh(X)$-action on $\IndCoh(\Z \times V^{\vee}/\G_m)$
corresponding to the pullback along the
composition $\Z \times V^{\vee}/\G_m \xrightarrow{p} \Z
\xrightarrow{i_{\Z}} X$.
It is obvious from the following commutative
diagram that $p^!$ is $\QCoh(X)$-equivariant
with respect to this action and the one on $\IndCoh(\Z)$:  
\[\begin{tikzcd}
	{X \times V^{\vee}/\G_m} & {\Z \times V^{\vee}/\G_m} \\
	X & \Z
	\arrow["{\pr_1}", from=1-1, to=2-1]
	\arrow["p", from=1-2, to=2-2]
	\arrow["{i_{\Z}}"', hook', from=2-2, to=2-1]
	\arrow["{i_{\Z} \times \id}"', hook', from=1-2, to=1-1].
\end{tikzcd}\]
Similarly, it is obvious from the following
commutative diagram that $i_*^{\IndCoh}$
is $\QCoh(X)$-equivariant:
\[\begin{tikzcd}
	& X \\
	{\Z \times V^{\vee}/\G_m} & {\wit{\Z}/\G_m}
	\arrow["{\pr_1 \circ i_{\wit{\Z}/\G_m}}"', from=2-2, to=1-2]
	\arrow["i", from=2-1, to=2-2]
	\arrow["{i_{\Z} \circ \pr_1}", from=2-1, to=1-2].
\end{tikzcd}\]
Since both $p^!$ and $i_*^{\IndCoh}$ have
$\QCoh(X)$-equivariant structures, it follows
that their composition $i_*^{\IndCoh}p^!$ does,
as well.
\eproof

The $\QCoh(X)$-module and $\QCoh(V^{\vee}/\G_m)$-module
structures on $\IndCoh(\wit{\Z}/\G_m)$ together induce 
the $\QCoh(X \times V^{\vee}/\G_m)$-module
structure on $\IndCoh(\wit{\Z}/\G_m)$ corresponding
to pullback along the inclusion $\wit{\Z}/\G_m \hook X \times V^{\vee}/\G_m$,
as shown in the following lemma.

\blem
\label{lem: equi 3}
Let $p: \calZ \to \calX$ and $q: \calZ \to \calY$
be maps, where $\calX$ and $\calY$ are perfect stacks
and $\calZ$ is a prestack locally almost of finite type.
Then the $\QCoh(\calX) \otimes \QCoh(\calY)$-module
structure on $\IndCoh(\calZ)$ corresponding to the tensor
product of module structures corresponding to $p$ and $q$ 
is the same as the one corresponding to the map
$p \times q: \calZ \to \calX \times \calY$ under the equivalence
$\QCoh(\calX) \otimes \QCoh(\calY) \simeq \QCoh(\calX \times \calY)$.
\elem

\bproof
Recall that the equivalence $\QCoh(\calX) \otimes
\QCoh(\calY) \simeq \QCoh(\calX \times \calY)$ is
induced by the exterior product $- \boxtimes - :=
\pr_1^*(-) \otimes \pr_2^*(-)$.
We have the following commutative diagram
by universality,
\[\begin{tikzcd}
	& \calZ \\
	& {\calX \times \calY} \\
	\calX && \calY
	\arrow["{\pr_1}", from=2-2, to=3-1]
	\arrow["{\pr_2}"', from=2-2, to=3-3]
	\arrow[from=1-2, to=2-2]
	\arrow["q", from=1-2, to=3-3]
	\arrow["p"', from=1-2, to=3-1],
\end{tikzcd}\]
from which it follows that the composition $(p \times q)^* \circ (\pr_1^* \otimes \pr_2^*):
\QCoh(\calX) \otimes \QCoh(\calY) \to \QCoh(\calZ)$ is
isomorphic to the functor $p^* \otimes q^*:
\QCoh(\calX) \otimes \QCoh(\calY) \to \QCoh(\calZ)$.
This suffices to prove the claim.
\eproof

\subsubsection{}

Putting it all together, we obtain the first main proposition of
this section,
establishing the $\QCoh(X \times V^{\vee}/\G_m)$-linearity
of the functor in \cref{thm: fancy Isik's theorem}.

\begin{enumerate}[label=(\alph*), ref=\alph*]
\item
\label{eyetem a}
Consider $\IndCoh(\Z)$ with the $\QCoh(X \times V^{\vee}/\G_m)$-module
structure induced by $\QCoh(X)$-module structure (\labelcref{enum 1})
and the standard $\QCoh(V^{\vee}/\G_m)$-module structure
of \cref{construction: singular support QCoh action}.
\item
\label{eyetem b}
The quotient functor $Q: \IndCoh(\wit{\Z}/\G_m) \to \Sing^{\infty}(\wit{\Z}/\G_m)$
is symmetric monoidal, so $\Sing^{\infty}(\wit{\Z}/\G_m)$ obtains a 
$\QCoh(X \times V^{\vee}/\G_m)$-module by composing the pullback
\[i_{\wit{\Z}/\G_m}^*: \QCoh(X \times V^{\vee}/\G_m) \to \QCoh(\wit{\Z}/\G_m)\] with
the composition $\QCoh(\wit{\Z}/\G_m) \xrightarrow{\Upsilon} \IndCoh(\wit{\Z}/\G_m)
\to \Sing^{\infty}(\wit{\Z}/\G_m)$.
\end{enumerate}

\bprop
\label{prop: showing equivariance}
By abuse of notation, let $i_*^{\IndCoh}p^!$ denote
the composition $Q i_*^{\IndCoh}p^!$. Then
$i_*^{\IndCoh}p^!$ is $\QCoh(X \times V^{\vee}/\G_m)$-linear 
for the (\labelcref{eyetem a}) module structure on the source and the
(\labelcref{eyetem b}) module structure on the target.
\eprop

\bproof
It follows from \cref{cor: equi 1}, \cref{lem: equi 2}, 
and \cref{lem: equi 3} that the
functor $i_*^{\IndCoh}p^!: \IndCoh(\Z) \to \IndCoh(\wit{\Z}/\G_m)$ is
$\QCoh(X \times V^{\vee}/\G_m)$-linear. On the other hand,
$Q$ is $\QCoh(X \times V^{\vee}/\G_m)$-linear almost by
definition, so the proposition follows.
\eproof

\subsection{Equivalence}
It remains to show that the map $i_*^{\IndCoh}p^!$ is an equivalence,
for which it suffices to prove the following proposition by the argument
given in the proof of \cref{thm: fancy Isik's theorem}

\bprop
\label{prop: showing equivalence}
The map $i_*^{\IndCoh}p^{\IndCoh, *}: \IndCoh(\Z) \to \Sing^{\infty}(\wit{\Z}/\G_m)$ 
is an equivalence of $k$-linear categories.
\eprop

\bproof
Note that since $i$ is proper locally almost of finite
type and $p$ is eventually coconnective, $i_*^{\IndCoh}$
and $p^*_{\IndCoh}$ each preserve compact objects.
Moreover, since each functor is continuous, their restrictions to $\Coh$
determines functors,
\begin{align*}
i_*: \Coh(\Z \times V^{\vee}/\G_m) \to \Coh(\wit{\Z}/\G_m) \\
p^*: \Coh(\Z) \to \Coh(\Z \times V^{\vee}/\G_m),
\end{align*}
whose ind-extensions are $i_*^{\IndCoh}$ and $p^*_{\IndCoh}$,
respectively. Thus, it suffices to show that $i_*p^*: \Coh(\Z)
\to \Sing(\wit{\Z}/\G_m)$ is an equivalence.

But now recall that we have an equivalence,
\[\Orl^{-1, \otimes}_{\gr}: \Sing(\wit{\Z}/\G_m) \xrightarrow{\simeq} \MF^{\cl}(\wit{Z}/\G_m)\]
by \cref{lem: graded MF is idempotent completion of graded MFcl}.
Thus, it suffices to show that 
\beqn
\label{eqn: hlurb}
\Orl^{-1,\otimes}_{\gr} \circ (i_*p^*): \Coh(\Z) \to \MF^{\cl}(\wit{Z}/\G_m)
\eeqn
is an equivalence. 

Unwinding the construction of this functor equivalence and comparing
to the constructions in \cite{Shipman}, we see that \labelcref{eqn: hlurb}
is precisely the map appearing in \cite[Theorem 3.5]{Shipman}, which
is shown in \textit{loc. cit.} to be an equivalence.
\eproof



\section{Comparing Hochschild invariants}
\label{sec: comparing Hochschild invariants}
Given a graded Landau--Ginzburg pair $(M,f,d)$, one may choose
to ignore the equivariant structure of $f$ and work
instead with just the underlying ungraded Landau--Ginzburg pair
$(M,f)$.
The data of a graded Landau--Ginzburg pair and its
underlying ungraded pair separately determine
categories of graded and ungraded matrix 
factorizations. When $d=2$ and the $\G_m$-action
on $M$ is even, these categories are $k$-linear
and $2$-periodic, respectively.
In general, these categories are quite different;
it is typically not possible to recover $\MF(M,f)$
from $\MF(M/\G_m,f)$ by $2$-periodization.
In this section, we will show that 
these different categories \emph{do}, however, have
the same periodic cyclic homology after $2$-periodization.

The starting point for our endeavor is the observation
that Hochschild invariants behave well with respect to
group actions on categories, and that both $\MF(M,f)$
and $\MF(M/\G_m,f)$ admit descriptions as the fixed
point of a certain formal group action on the category
of coherent sheaves. Indeed, recall from 
\S\labelcref{ssec: alternative definitions} that
$\PreMF(M,f)$ is invariants for
certain $B\wih{\G}_a$-action on $\Coh(M)$,
and $\PreMF(M/\G_m,f)$ is invariants
for a certain $B\wih{\G}_a/\G_m$-action on $\Coh(M/\G_m)$.
The data of these formal group actions is obtained
from the data of $f$ as a map and as a $\G_m$-equivariant
map, respectively, and one may be obtained from the
other by passing under the forgetful functor from
the $\G_m$-equivariant to the non-equivariant setting.
More precisely, we have the following lemma.

\blem
\label{lem: same module structures}
Suppose that $(M,f,d)$ is a graded Landau--Ginzburg pair of
weight $d$ over $\Spec k$. By \cref{lem: map gives formal
group action}, the map $f$ induces a $B\wih{\G}_a$-module
structure on $\Coh(M)$, and by \cref{lem: equivariant map gives
graded formal group action}, the $\G_m$-equivariant structure
on $f$ induces a $B\wih{\G}_a/\G_m$-module
structure on $\Coh(M/\G_m)$. We claim that 
the $B\wih{\G}_a$-module $\Coh(M)$ naturally identifies with the
image of the $B\wih{\G}_a/\G_m$-module 
$\Coh(M/\G_m)$ under the enhanced forgetful functor
$\oblv: B\wih{\G}_a/\G_m\bmod \to B\wih{\G}_a\bmod$. 
\elem

\bproof
This is obvious from the constructions of the
two module structures.
\eproof

Forgetful functors in this context are always symmetric monoidal,
and Hochschild homology, being defined in terms of the symmetric monoidal
structure of a category, commutes with such maps, as recorded
in the following lemma.

\blem
\label{lem: HH of symmetric monoidal}
Suppose that $F: \scrC \to \scrD$ is a symmetric monoidal
functor of symmetric monoidal, $(\infty,2)$-categories. Let
$c \in \scrC$ be a dualizable object. Then $F(c)$ is dualizable
in $\scrD$, and there is a natural
isomorphism,
\beqn
\label{eqn: HH of symmetric monoidal}
F(\HH_{\bullet}^{\scrC}(c)) \simeq \HH_{\bullet}^{\scrD}(F(c)),
\eeqn
of functors from $\scrC$ to $S^1\mod(\End_{\scrD}(1_{\scrD}))$.
\elem

\bproof
Follows from the functoriality of categorical traces.
\eproof

We will use \cref{lem: same module structures}
as well as various functorialities of Hochschild invariants
such as \cref{lem: HH of symmetric monoidal} in order
to relative the periodic cyclic homology of $\MF(M,f)$
to that of $\MF(M/\G_m,f)$.

\brem
\label{rem: HH(C) = HH(Ind(C))}
A comparison of periodic cyclic homology
for large categories matrix factorizations
can be deduced from the case of small 
categories using the following standard fact.

Suppose that $\scrC$ is a small idempotent-complete
dg-category which is dualizable as an object
of $\dgcat_k^{\idm}$. Then $\HH_{\bullet}^k(\scrC) 
\simeq \HH_{\bullet}^k(\Ind(\scrC))$.
\erem

\subsection{$\G_m$-actions on small categories}
Let $G$ be an affine algebraic group or 
formal group. Taking the $G$-invariants of a small $G$-category
does not commute with taking ind-categories.
Because we will be taking invariants of formal
group actions on categories, we prefer to work
with small categories throughout this section.
The theory of small graded categories is well-known,
but not well-documented in the literature, so we
establish some basic notation and facts below.

\bdef
Let $\G_m\bmod$ denote the category of 
comodules in $\dgcat_k^{\idm}$ over $\Perf(\G_m)
\in \coAlg(\dgcat_k^{\idm})$ whose underlying dg-category
lies in $\dgcat_k^{\idm}$,
and let $\Perf(B\G_m)\mod$ denote the category of
modules in $\dgcat_k^{\idm}$ over the symmetric monoidal category $(\Perf(B\G_m), \otimes)
\in \Alg(\dgcat_k^{\idm})$.
\edefn

There are functors of $\G_m$-variants and de-equivariantization
comprise a pair of adjoint functors,
\[\begin{tikzcd}
	{(-)^{\G_m}: \G_m\bmod} && {\Perf(B\G_m)\mod: \oblv}
	\arrow[shift right=1, from=1-1, to=1-3]
	\arrow[shift right=1, from=1-3, to=1-1]
\end{tikzcd}\]
where we have chosen to denote the
de-equivariantization functor, $- \otimes_{\Perf(B\G_m)} \Perf$, 
by $\oblv$.

\begin{warn}
The reader might be more familiar with the set-up
in which we consider
\emph{large} $\QCoh(\G_m)$-comodules
and $\Rep(\G_m)\mod$-modules---i.e. comodules and modules, respectively,
whose underlying dg-categories lie in $\dgcat_k^{\infty}$.
In this setting, the $\G_m$-invariants--de-equivariantization
adjunction is an adjoint equivalence 
because $B\G_m$ is $1$-affine in
the sense of \cite{Gaitsgory1affineness}.
The adjunction considered here is not an adjoint
equivalence since it is not necessarily the case that
$\scrC^{\G_m}$ is compactly generated if
$\scrC$ is compactly generated. In particular,
if $\scrC$ is small,
$\Ind(\scrC^{\G_m})$ is not necessarily
equivalent to $\Ind(\scrC)^{\G_m}$.
\end{warn}

By abuse of notation, we will also use $\oblv$ to 
denote the composition of $\oblv$ with the 
forgetful functor $\G_m\bmod \to \dgcat_k^{\idm}$.

\brem
We will refer to objects of $\G_m\bmod$ as small ``$\G_m$-categories,"
and objects of $\Perf(B\G_m)\mod$ as small ``graded categories."
\erem


\subsection{Computations and the comparison isomorphism}

\subsubsection{}
In this section, we will make substantial use of the material on
derived formal group actions on categories found in
\S\labelcref{ssec: Formal group actions on categories}.

\subsubsection{}
\label{sssec: grading on BG_a^}
A choice of grading 
on the polynomial algebra $\O(\G_a)$ equips the
derived formal group $B\wih{\G}_a$ with
an algebraic $\G_m$-action.

For the remainder of this section, we use the
$\G_m$-action on $B\wih{\G}_a$ corresponding
to double the usual scaling action on $\G_a \simeq \Aone$,
i.e. $t \cdot x = t^2x$.

\subsubsection{}
For the remainder of this section, we fix a graded
Landau--Ginzburg pair $(M,f,2)$.

\subsubsection{De-equivariantization}

\bdef
Let $B\wih{\G}_a/\G_m\bmod$
denote the category of $B\wih{\G}_a/\G_m$-modules\footnotemark
in $\Perf(B\G_m)\mod$.
Similarly, let $B\wih{\G}_a\bmod$,
denote the category of $B\wih{\G}_a$-modules in $\dgcat_k^{\idm}$.
\footnotetext{See \cref{def: modules over a formal group} for the general
formalism of formal groups acting on categories.}
\edefn

\begin{warning}
The reader should be careful to distinguish between
\[B\wih{\G}_a\bmod \; \textnormal{(resp.} \; B\wih{\G}_a/\G_m\bmod\textnormal{)},\] 
by which we mean $B\wih{\G}_a\mod(\dgcat_k^{\idm})$
(resp. $B\wih{\G}_a/\G_m\mod(\Rep(\G_m)\mod)$), and
\[B\wih{\G}_a\mod \; \textnormal{(resp.} \; B\wih{\G}_a/\G_m\mod\textnormal{)},\] 
by which mean $B\wih{\G}_a\mod(\Vect_k)$ (resp. $B\wih{\G}_a/\G_m\mod(\Rep(\G_m))$).
\end{warning}

Evidently, the forgetful functor $\oblv: \Perf(B\G_m)\mod \to \dgcat_k^{\idm}$
induces a functor,
\[\oblv: B\wih{\G}_a/\G_m\bmod \to B\wih{\G}_a\bmod.\]
One might hope that taking invariants commutes
with the passage from the graded to the ungraded
setting in the sense that the following diagram
commutes,
\[
\begin{tikzcd}
	{B\wih{\G}_a/\G_m\bmod} & {} & {\left(k[\b]\mod_{\gr}\right)\bmod} \\
	{B\wih{\G}_a\bmod} & {} & {\left(k[\b]\mod\right)\bmod}
	\arrow["\oblv", from=1-3, to=2-3]
	\arrow["\oblv", from=1-1, to=2-1]
	\arrow["{(-)^{B\wih{\G}_a}}", from=2-1, to=2-3]
	\arrow["{(-)^{B\wih{\G}_a/\G_m}}", from=1-1, to=1-3].
\end{tikzcd}
\]
In fact, this is generally \emph{not} the case, and in particular,
for a graded scheme $X$ equipped with an equivariant
function $f: X \to \Aone$,
$\oblv(\Coh(X/\G_m)^{B\wih{\G}_a/\G_m})$
is not equivalent to $\Coh(X)^{B\wih{\G}_a}$.
On the other hand, there is a natural
transformation of functors,
\beqn
\label{eqn: oblv and invariants map}
\oblv((-)^{B\wih{\G}_a/\G_m}) \to (\oblv(-))^{B\wih{\G}_a},
\eeqn
coming from the universal property of group invariants.

\subsubsection{Hochschild homology}
The usual functors of categorical dimension (see
\S\labelcref{sec: Hochschild invariants}) in
$\Rep(\G_m)\mod$ and $\dgcat_k^{\idm}$
induce functors, 
\begin{align*}
&\HH_{\bullet}^{\Rep(\G_m)}: B\wih{\G}_a/\G_m\bmod \to S^1\mod(B\wih{\G}_a/\G_m\mod) \\
&\HH_{\bullet}^k: B\wih{\G}_a\bmod \to S^1\mod(B\wih{\G}_a\mod),
\end{align*}
respectively,
by \cite[Proposition 5.3.3.11]{PreygelThesis}
and \cref{prop: graded version of Toly Prop 5.3.3.11},
respectively. As above, one might hope
that the natural transformations obtained from the universal
property of group invariants,
\begin{align}
\label{map 1} \HH_{\bullet}^{k[\b]_{\gr}}\left((-)^{B\wih{\G}_a/\G_m}\right) &\to \left(\HH_{\bullet}^{\Rep(\G_m)}(-)\right)^{B\wih{\G}_a/\G_m} \\
\label{map 2} \HH_{\bullet}^{k[\b]}\left((-)^{B\wih{\G}_a}\right) &\to \left(\HH_{\bullet}^k((-))\right)^{B\wih{\G}_a}
\end{align}
are natural isomorphisms of
functors from $B\wih{\G}_a/\G_m\bmod$ to $S^1\mod(k[\b]\mod_{\gr})$,
and from $B\wih{\G}_a\bmod$ to $S^1\mod(k[\b]\mod)$,
respectively.

\bprop
\label{prop: natural BG_a maps are equivalences}
Let $(M,f,d)$ be a graded Landau--Ginzburg
pair. Then components of the natural transformations \labelcref{map 1}
and \labelcref{map 2} at
the objects $\Coh(M/\G_m)$ and $\Coh(M)$, respectively,
are equivalences.
\eprop

\bproof
The claim for \labelcref{map 2} is the content of
\cite[Theorem 6.1.2.5(i)]{PreygelThesis}. In order to
prove the claim for \labelcref{map 1}, we imitate
the proof found therein, in the graded setting. 

In order to do so, we use
\cite[Theorem 4.1.3.6]{PreygelThesis} in order
to furnish a $k[\b]_{\gr}$-linear equivalence,
\[\Fun^L_{k[\b]_{\gr}}(\PreMF^{\infty}(M/\G_m, f), \PreMF^{\infty}(M/\G_m, f)) 
	\xrightarrow{\simeq} \PreMF^{\infty}((M \times M)/\G_m, -f \boxplus f),\]
under which the identity functor
$\id_{\PreMF^{\infty}(M/\G_m,f)}$ is sent to
$\over{\Delta}_*\omega_{M/\G_m}$ and the trace $\tr(-)$
of an endofunctor of $\PreMF^{\infty}(M/\G_m,f)$ corresponds to
\[\Hom_{k[\b]}(\over{\Delta}_*\O_{M/\G_m}, -),\]
where $\over{\Delta}: M/\G_m \to ((M/\G_m)^2)_0$ corresponds to
the inclusion of $M$ into the zero locus of $-f \boxplus f$
via the diagonal.

Now, assuming a graded version of Preygel's
\cite[Proposition 3.1.2.1]{PreygelThesis}, 
we have that, for any $\F, \calG \in \PreMF^{\infty}((M/\G_m)^2,-f \boxtimes f)$,
there is a natural graded $B\wih{\G}_a/\G_m$-action (equivalently,
a graded $S^1$-action) on
the graded $k$-linear mapping spectrum between $i_*\F$ and
$i_*\calG$ as objects of $\IndCoh((M/\G_m)^2)$,
$\Hom_{B\G_m}(i_*\F, i_*\calG)$, where 
$i: (-f \boxplus f)^{-1}(0) \to (M/\G_m)^2$ is the
canonical inclusion; and, moreover, that
\[\Hom_{k[\b]_{\gr}}(\F, \calG) \simeq \Hom_{B\G_m}({i_{(M/\G_m)^2_0}}_*\F, {i_{(M/\G_m)^2_0}}_*\calG)^{B\wih{\G}_a/\G_m}.\]
In particular, we have 
\[\Hom_{k[\b]_{\gr}}(\over{\Delta}_*\O_{M/\G_m}, \over{\Delta}_*\omega_{M/\G_m}) 
	\simeq \left(\Hom_{B\G_m}(\Delta_*\O_{M/\G_m}, \Delta_*\omega_{M/\G_m})\right)^{B\wih{\G}_a/\G_m},\]
where $\Delta: M/\G_m \to (M/\G_m)^2$ is
the diagonal embedding. On the other hand, by
\cref{prop: integral kernel theorem}, the $\Rep(\G_m)$-linear
Hochschild homology of $\IndCoh(M/\G_m)$ corresponds
precisely to $\Hom_{B\G_m}(\Delta_*\O_{M/\G_m}, \Delta_*\omega_{M/\G_m})$ 
It follows that
\[\HH_{\bullet}^{k[\b]_{\gr}}(\PreMF^{\infty}(M/\G_m,f)) 
	\simeq \left(\HH_{\bullet}^{\Rep(\G_m)}(\IndCoh(M/\G_m))\right)^{B\wih{\G}_a/\G_m}.\]
The result now follows from the equivalence $\PreMF(M/\G_m,f) \simeq
\Coh(M/\G_m)^{B\wih{\G}_a/\G_m}$.
\eproof


It follows from \cref{lem: HH of symmetric monoidal}
and \cref{prop: natural BG_a maps are equivalences} 
that applying the functor $\HH_{\bullet}^k$ to the natural
map \labelcref{eqn: oblv and invariants map}
for $\scrC = \Coh(M/\G_m)$
obtains the following morphism of
$k[\b]$-modules,
\beqn
\label{eqn: natural map of HH inv}
\oblv\left(\left(\HH_{\bullet}^{\Rep(\G_m)}(\Coh(M/\G_m))\right)^{B\wih{\G}_a/\G_m}\right) 
	\to \left(\oblv\left(\HH_{\bullet}^{\Rep(\G_m)}(\Coh(M/\G_m))\right)\right)^{B\wih{\G}_a}.
\eeqn
It is easy to see from the construction of
the above map that it is the same map obtained
from the universal property of $B\wih{\G}_a$-invariants.
That is, the following diagram commutes up
to natural transformation,
\beqn
\label{eqn: decat invariant diagram}
\begin{tikzcd}
	{S^1\mod(B\wih{\G}_a/\G_m\mod)} & {} & {S^1\mod(k[\b]\mod_{\gr})} \\
	{S^1\mod(B\wih{\G}_a\mod)} & {} & {S^1\mod(k[\b]\mod)}
	\arrow["{(-)^{B\wih{\G}_a}}", from=2-1, to=2-3]
	\arrow["{(-)^{B\wih{\G}_a/\G_m}}", from=1-1, to=1-3]
	\arrow["\oblv"', from=1-1, to=2-1]
	\arrow["\oblv"', from=1-3, to=2-3],
\end{tikzcd}
\eeqn
coming from the universal property of
$B\wih{\G}_a$-invariants, and \labelcref{eqn: natural map
of HH inv} is given by
component of this natural transformation
at the object
$\HH_{\bullet}^{\Rep(\G_m)}(\Coh(M/\G_m)) \in S^1\mod(B\wih{\G}_a/\G_m\mod)$.

\blem
\label{lem: when 2-morphism is invertible}
The natural transformation,
\beqn
\label{eqn: oblv and invariants natural transformation}
\oblv\left((-)^{B\wih{\G}_a/\G_m}\right) \to \left(\oblv(-)\right)^{B\wih{\G}_a},
\eeqn
indicated in diagram \labelcref{eqn: decat invariant diagram},
restricts to a natural isomorphism on the
full subcategory on objects whose underlying
$B\wih{\G}_a/\G_m$-module belongs belongs
to $B\wih{\G}_a/\G_m\mod^+$.
\elem

\bproof
This is a general fact about $B\wih{\G}_a$-actions.
Suppose that $F: \scrC \to \scrD$ preserves colimits;
$\scrC$ and $\scrD$ have right complete t-structures
that are compatible with filtered colimits; and $F$ is t-exact.
Then $F$ restricted to $\scrC^+$ intertwines $B\wih{\G}_a$-invariants.
\eproof

\bcor
\label{cor: the final piece}
\labelcref{eqn: natural map of HH inv}
is an equivalence of objects in 
$S^1\mod(k[\b]\mod)$.
\ecor

\bproof
It suffices by \cref{lem: when 2-morphism is invertible}
to show that $\HH_{\bullet}^{\Rep(\G_m)}(\Coh(M/\G_m))$
is bounded above in the t-structure on 
$B\wih{\G}_a/\G_m\mod(\Rep(\G_m))$, or, equivalently,
in the t-structure on $\Rep(\G_m)$.

Since $B\G_m$ is a very good stack, we may write
\begin{align*}
\HH_{\bullet}^{\Rep(\G_m)}(\Coh(M/\G_m)) 	&= \HH_{\bullet}^{\Rep(\G_m)}(\IndCoh(M/\G_m)) \\ 
									&\simeq \Hom_{B\G_m}(\Delta_*^{\IndCoh}\O_{M/\G_m}, \Delta_*^{\IndCoh}\omega_{M/\G_m}) \\
									&\simeq \Hom_{B\G_m}(\O_{M/\G_m}, \Delta^!\Delta_*^{\IndCoh}\omega_{M/\G_m}) \\
									&\simeq {\pt_{M/\G_m}}_*\Delta^!\Delta_*^{\IndCoh}\omega_{M/\G_m}
\end{align*}  
where $\Delta: M/\G_m \to M/\G_m \times_{B\G_m} M/\G_m$ is the diagonal map
and $\pt_{B\G_m}: M/\G_m \to B\G_m$ is the structure map to $B\G_m$.
Since $\Delta$ is a proper map locally almost of finite type,
$\Delta_*^{\IndCoh}\omega_{M/\G_m}$
belongs to $\Coh(M^2)$, hence is t-bounded. Moreover,
$\Delta$ is also of bounded Tor dimension, so $\Delta^!$
preserves $\Coh$, hence $\Delta^!\Delta_*^{\IndCoh}\omega_{M/\G_m}$
is t-bounded. 

Let $p: M \to M/\G_m$ denote the canonical
quotient map. Using base change, it suffices to
show that the pushforward of $p^*\Delta^!\Delta_*^{\IndCoh}\omega_{M/\G_m}$
along the terminal map $\pt :M \to \ast$ is t-bounded.
But this follows from the assertion, proven in
\cite[Lemma 75.6.1]{Stacks} that $\pt_*$ has finite
cohomological amplitude since $M$ is finite type (so
in particular quasi-compact). 

In particular, it follows that $\HH_{\bullet}^{\Rep(\G_m)}(\Coh(M/\G_m)) \in
\IndCoh(B\G_m)^+ \simeq \Rep(\G_m)^+$, as desired.
\eproof

Altogether, we obtain the following proposition.

\bprop
\label{prop: graded and ungraded PreMF have the same HH}
Let $(M,f,2)$ be a graded Landau--Ginzburg pair
of weight $2$. Then there is an equivalence,
\beqn
\label{eqn: HH of gr/usual MF}
\oblv\left(\unHH_{\bullet}^{\kb_{\gr}}(\MF^{\infty}(M/\G_m,f))\right) 
	\simeq \unHH_{\bullet}^{\kb}(\MF^{\infty}(M,f)),
\eeqn
of $S^1\mod(\kb\mod)$-valued sheaves on $(X/\G_m)_{\et}$.
\eprop

\begin{notation}
\label{notation: enhanced Coh}
We first introduce some notation for use in our proof
of \cref{prop: graded and ungraded PreMF have the same HH}
which does not appear elsewhere in this work.

As shown in a previous section, the choice of a
regular function $f$ on a classical scheme
$X$ of finite type is equivalent to the 
datum of a $B\wih{\G}_a$-action
on $\Coh(X)$. Similarly, if $X$ is a $\G_m$-scheme, the
choice of a $\G_m$-equivariant function $f$ of weight $2$
is equivalent to the datum of of a $B\wih{\G}_a/\G_m$-action
on $\Coh(X/\G_m)$, where $B\wih{\G}_a$ is taken
to have the $\G_m$-action specified in \S\labelcref{sssec: grading on BG_a^}.
In fact, the assignment of a pair $(X,f) \mapsto \Coh(X) \in B\wih{\G}_a\bmod$ 
(respectively, $(X/\G_m, f) \mapsto \Coh(X/\G_m) \in B\wih{\G_a}/\G_m\bmod$)
is functorial in the pair $(X,f)$ (respectively, $(X/\G_m, f)$). More
precisely, for a fixed (graded) Landau--Ginzburg pair
$(X,f)$, there are presheaves,
\begin{align*}
\Coh^{f,\enh}(-): (X_{\et})^{\op} &\to B\wih{\G}_a\bmod \\
\Coh^{f,\enh}(-/\G_m): ((X/\G_m)_{\et})^{\op} &\to B\wih{\G}_a/\G_m\bmod,
\end{align*}
which assign to a $\G_m$-equivariant \'etale map $U \to X \in (X/\G_m)_{\et}$
the $B\wih{\G}_a$-module $\Coh(U)$ and the $B\wih{\G}_a/\G_m$-module
$\Coh(U/\G_m)$, respectively.

Moreover, these functors are compatible
in the sense that there is a commutative diagram,
\[\begin{tikzcd}
	{((X/\G_m)_{\et})^{\op}} & {((X/\G_m)_{\et})^{\op}} \\
	{B\wih{\G}_a\bmod} & {B\wih{\G}_a/\G_m\bmod}
	\arrow["{(-)^{\G_m}}", from=2-1, to=2-2]
	\arrow["{{\Coh^{f,\enh}(-)}|_{((X/\G_m)_{\et})^{\op}}}"', from=1-1, to=2-1]
	\arrow["{\Coh^{f,\enh}(-/\G_m)}", from=1-2, to=2-2]
	\arrow[shift left=1, no head, from=1-1, to=1-2]
	\arrow[no head, from=1-1, to=1-2].
\end{tikzcd}\]
\end{notation}

\bproof[Proof of \cref{prop: graded and ungraded PreMF have the same HH}]
It suffices to show the corresponding equivalence of presheaves,
since $\oblv$ commutes with sheafification.

In order to ease the exposition and make our
notation more consistent with that introduced
in \cref{notation: enhanced Coh} above, we denote the (pre-)sheaves
$\unPreMF(M,f)$ and $\unMF(M,f)$,
by $\PreMF(-,f)$ and $\MF(-,f)$,
where the scheme $M$ is clear from context.

From the various results
obtained above, we have the following chain of
natural isomorphisms of presheaves
$((M/\G_m)_{\et})^{\op} \to S^1\mod(k[\b]\mod)$:
\begin{align}
\label{1} \oblv\left(\HH_{\bullet}^{k[\b]_{\gr}}(\PreMF(-/\G_m, f))\right) 	&:= \oblv\left(\HH_{\bullet}^{k[\b]_{\gr}}\left(\Coh(-/\G_m)^{B\wih{\G}_a/\G_m}\right)\right) \tag{1} \\
\label{2}													&\xrightarrow{\simeq} \oblv\left(\left(\HH_{\bullet}^{\Rep(\G_m)}\left(\Coh(-/\G_m)\right)\right)^{B\wih{\G}_a/\G_m}\right) \tag{2} \\
\label{3}													&\xrightarrow{\simeq} \left(\oblv\left(\HH_{\bullet}^{\Rep(\G_m)}\left(\Coh(-/\G_m)\right)\right)\right)^{B\wih{\G}_a} \tag{3} \\
\label{4}													&\simeq \left(\HH_{\bullet}^{k[\b]}\left(\oblv(\Coh(-/\G_m))\right)\right)^{B\wih{\G}_a} \tag{4} \\
\label{5}													&\simeq \left(\HH_{\bullet}^{k[\b]}\left(\Coh(-)\right)\right)^{B\wih{\G}_a} \tag{5} \\
\label{6}													&\xleftarrow{\simeq} \HH_{\bullet}^{k[\b]}\left(\Coh(-)^{B\wih{\G}_a}\right) \tag{6} \\
\label{7}													&= \HH_{\bullet}^{k[\b]}(\PreMF(-,f)). \tag{7}
\end{align}
We now justify each natural
isomorphism in this chain.
\begin{itemize}
\item[\labelcref{1}] is a definition.
\item[\labelcref{2}] is obtained by the precomposition of the natural transformation
\labelcref{map 1} with the functor $\Coh^{f,\enh}(-/\G_m): ((M/\G_m)_{\et})^{\op} 
\to B\wih{\G}_a/\G_m\mod$, and this composition is a natural isomorphism by
\cref{prop: natural BG_a maps are equivalences}.
\item[\labelcref{3}] is obtained
by composing the presheaf $\HH_{\bullet}^{\Rep(\G_m)}(\Coh(-/\G_m))$
with the natural transformation
\labelcref{eqn: oblv and invariants natural transformation},
and this composition is an isomorphism by \cref{cor: the final piece}.
\item[\labelcref{4}] is obtained from the natural isomorphism
\labelcref{eqn: HH of symmetric monoidal}.
\item[\labelcref{5}] is obtained using the equivalence $\oblv(\Coh(M/\G_m)) \simeq \Coh(M)$,
which is natural in $M$.
\item[\labelcref{6}] is obtained in the same way as \labelcref{2}, but
instead composing the functor ${\Coh^{f,\enh}(-)}|_{((M/\G_m)_{\et})^{\op}}$ 
with \labelcref{map 2}. This composition is a natural isomorphism,
also by \cref{prop: natural BG_a maps are equivalences}.
\item[\labelcref{7}] is the identification found in \cref{lem: same module structures}.
\end{itemize}

Because the functor $- \otimes_{k[\b]} \kb$ is symmetric
monoidal, we have the following natural isomorphisms,
using \cref{lem: HH of symmetric monoidal} and \cref{rem: HH(C) = HH(Ind(C))},
\begin{align*}
\HH_{\bullet}^{k[\b]}(\PreMF(-,f)) \otimes_{k[\b]} \kb 	&\simeq \HH_{\bullet}^{\kb}(\PreMF(-,f) \otimes_{k[\b]} \kb\mod) \\
											&=: \HH_{\bullet}^{\kb}(\MF(-,f)).
\end{align*}
On the other hand, by the chain of natural isomorphisms
\labelcref{1}--\labelcref{7} written above, we obtain
\begin{align*}
\simeq \oblv\left(\HH_{\bullet}^{\kb_{\gr}}(\MF(-/\G_m,f))\right)		&\simeq \oblv\left(\HH_{\bullet}^{\kb_{\gr}}(\PreMF(-/\G_m,f) \otimes_{k[\b]_{\gr}} \kb\mod_{\gr})\right) \\
													&\simeq \oblv\left(\HH_{\bullet}^{\kb_{\gr}}(\PreMF(-/\G_m,f)) \otimes_{k[\b]_{\gr}} \kb_{\gr}\right) \\
													&\simeq \oblv\left(\HH_{\bullet}^{\kb_{\gr}}(\PreMF(-/\G_m,f))\right) \otimes_{k[\b]} \kb \\
													&\simeq \HH_{\bullet}^{k[\b]}(\PreMF(-,f)) \otimes_{k[\b]} \kb.
\end{align*}
Thus we have exhibited a natural isomorphism between the
left-hand and right-hand sides of \labelcref{eqn: HH of gr/usual MF},
as desired. This completes the proof.
\eproof

Taking Tate invariants for the $S^1$-action 
on Hochschild homology obtains, by definition,
the periodic cyclic homology. The following corollary compares
the periodic cyclic homologies of graded and $2$-periodic matrix factorizations.

\bcor
\label{cor: HP of gr/usual MF}
Let $k[u] \simeq C^*(BS^1;k)$ denote another copy of
the algebra of cochains on $BS^1$.
There is an equivalence of sheaves of $k[\b^{\pm 1}, u^{\pm 1}]$-modules
on $(M/\G_m)_{\et}$,
\[\oblv\left(\unHP_{\bullet}^{\kb_{\gr}}(\MF^{\infty}(M/\G_m,f))\right) \simeq \unHP_{\bullet}^{\kb}(\MF^{\infty}(M,f)).\]
\ecor

\bproof
As before, it suffices to show the corresponding equivalence of presheaves
since $\oblv$ commutes with sheafification.

Note that each of the standard t-structures on $\kb\mod_{\gr}$
and $\kb\mod$ is right-complete and compatible with
filtered colimits. 
Since the functor $\oblv: \kb\mod_{\gr} \to \kb\mod$
preserves colimits (as a left adjoint), and is t-exact with
respect to these t-structures, its restriction to $\kb\mod_{\gr}^+$
intertwines $S^1$-invariants on objects with $S^1$-action.
That is, the following diagram commutes,
\[\begin{tikzcd}
	{S^1\mod({\kb\mod_{\gr}}^+)} & {k[u]\mod(\kb\mod_{\gr})_{\gr}} \\
	{S^1\mod(\kb\mod)} & {k[u]\mod(\kb\mod)}
	\arrow["\oblv"', from=1-1, to=2-1]
	\arrow["{(-)^{S^1}}", from=2-1, to=2-2]
	\arrow["{(-)^{S^1}}", from=1-1, to=1-2]
	\arrow["\oblv", from=1-2, to=2-2].
\end{tikzcd}\]
Thus we obtain the equivalence of presheaves,
\[\oblv\left(\HH_{\bullet}^{\kb_{\gr}}(\unMF^{\infty}(M/\G_m,f))^{S^1}\right) \simeq \HH_{\bullet}^{\kb}(\unMF^{\infty}(M,f))^{S^1}.\]
The $S^1$-invariants have a canonical $k[u]$-module
structure, and Tate invariants for $S^1$-actions
are given, by definition, by inverting the generator $u$.
Because $\oblv$ is symmetric monoidal,
it follows that,
\begin{align*}
\oblv\left(\HH_{\bullet}^{\kb_{\gr}}(\unMF^{\infty}(M/\G_m,f))^{S^1} \otimes_{k[u]_{\gr}} \ku_{\gr} \right) 	&\simeq \HH_{\bullet}^{\kb}(\unMF^{\infty}(M,f))^{S^1} \otimes_{k[u]} \ku \\
\oblv\left(\HP_{\bullet}^{\kb_{\gr}}(\unMF^{\infty}(M/\G_m,f))\right) 								&\simeq \HP_{\bullet}^{\kb}(\unMF^{\infty}(M,f)).
\end{align*}
\eproof

\subsection{Even $\G_m$-actions}
If the $\G_m$-action on $M$ is even, 
the story in the previous section admits 
a convenient simplification.
Namely, if the graded Landau--Ginzburg pair $(M,f)$
has an even action, then its graded matrix factorizations
category is naturally a plain dg-category
under shearing (see \cref{sec: shearing}
for the necessary background on shearing), 
and $\HP_{\bullet}^{\kb}(\MF^{\infty}(M,f))$
is simply the $2$-periodization of the complex
$\HP_{\bullet}^k(\MF^{\infty}(M/\G_m,f))$,
as we explain in this section.

\subsubsection{Even actions on schemes}
We begin by recalling the definition of an even
$\G_m$-action on a scheme.

\bdef
We will say that a scheme $X$ has an \emph{even action}
of $\G_m$ if the induced action by the subgroup of square
roots of unity $\mu_2 \subset \G_m$ is trivial.
\edefn

While much of the discussion until this point 
was valid for very general kinds of schemes,
we are ultimately concerned with the case of
a Landau--Ginzburg pair $(M,f)$ in which $M$
is a smooth quasi-projective variety.
It turns out that any even $\G_m$-action on
such a scheme actually comes from
squaring another (possibly not even) $\G_m$-action
on that scheme, as shown in the lemma following
the definition below.

\bdef
Let $\act: \G_m \times X \to X$ be a $\G_m$-action
on the classical scheme $X$. For every natural number
$d \geq 1$, let $\act_d: \G_m \times X \to X$ be the
action on $X$ given by precomposing $\act$ with the
$d$th power map, $(-)^d: \G_m \to \G_m$.
When $d=2$, we refer to this new action as
\emph{the square} of the $\act$.  
\edefn

\blem
Let $X$ be a normal algebraic $k$-variety.
Then any even $\G_m$-action on $X$ is obtained
as the square of some other $\G_m$-action on $X$.
\elem

\bproof
Suppose first that $X$ is any affine scheme
(not necessarily normal), so $X = \Spec R$
for an ordinary commutative $k$-algebra $R$. Then an even
$\G_m$-action on $X$ is easily seen to be equivalent
to a grading $R \simeq \oplus_{n \in \bbZ} R_n$ such that
$R_n = 0$ if $n$ is odd. Define a new grading $R = \oplus_{n \in \bbZ} A_n$ by
setting $A_n := R_{2n}$. It is immediate from the definitions that
original $\G_m$-action on $X$ is the square of
the $\G_m$-action on $X$ defined by the new grading
$\oplus_{n \in \bbZ} A_n$.

Now assume $X$ is as in the statement of the lemma,
and let $\act$ denote the $\G_m$-action on $X$.
By Sumihiro's theorem (\cite{Sumihiro}), $X$ admits
a cover by $\G_m$-invariant open affine subsets.
Fix such a cover $\{U_i\}$.
By the argument of the previous paragraph, the
restriction of $\act$ to each open
affine $U_i$ in the cover is the square of another
action $\act_i$ on $U_i$, defined in the way above. 
It therefore suffices to
show that for any $i, j \in I$,
\beqn
\label{eqn: in-house}
\act_i|_{\G_m \times U_{ij}} = \act_j|_{\G_m \times U_{ij}}.
\eeqn
Observe that, since both $U_i$ and $U_j$ are
$\G_m$-invariant (with respect to $\act$), 
their fiber product $U_{ij} := U_i \times_X U_j$
is as well.
Let 
\begin{itemize}
\item[$-$] $U_i = \Spec(R(i))$
\item[$-$] $U_j = \Spec(R(j))$
\item[$-$] $U_{ij} = \Spec(R(ij))$,
\end{itemize}
and let
\begin{itemize}
\item[$-$] $R(i) = \oplus_{n \in \bbZ} R(i)_n$ denote the grading
coming from $\act|_{\G_m \times U_i}$, and $\oplus_{n \in \bbZ} \calR(i)_n$
denote the grading coming from $\act_i$
\item[$-$] $R(j) = \oplus_{n \in \bbZ} R(j)_n$ denote the grading
coming from $\act|_{\G_m \times U_j}$, and $\oplus_{n \in \bbZ} \calR(j)_n$
the grading coming from $\act_j$.
\end{itemize}
The ring
$R(ij)$ has two distinct gradings, 
$\oplus_{n \in \bbZ} A_n$ and $\oplus_{n \in \bbZ} B_n$, 
induced by $\act_i$ and $\act_j$, respectively.
Checking the equality \labelcref{eqn: in-house} amounts
to checking that $A_n = B_n$,
for all $n \in \bbZ$. From the definition
of each grading, we obtain,
\begin{align*}
A_n 				&:= {\calR(i)_n}|_{U_{ij}} \\ 
				&:= {R(i)_{2n}}|_{U_{ij}} \\
				&= {R(j)_{2n}}|_{U_{ij}} \\
				&=: {\calR(j)_n}|_{U_{ij}} \\
				&=: B_n,
\end{align*}
where ${R(i)_{2n}}|_{U_{ij}} = {R(j)_{2n}}|_{U_{ij}} $
because both are, by the definition, the ${2n}^{\text{th}}$ graded 
piece of the decomposition
corresponding to the $\act$ $\G_m$-action on $U_{ij}$.
This proves the claim.
\eproof

Two different $\G_m$-actions $\act$ and $\act'$ 
on $X$ will produce two different $\G_m$-category structures
on $\IndCoh(X)$. When $\act' = \act_d$, however,
there is a nice description of how these two
structures relate to each other.

\blem
\label{lem: invariants and maps}
Let $f: G \to H$ be a map of algebraic group
schemes. Let $\scrC \in H\bmod$, and
let $\wit{\scrC}$ denote the object in
$G\bmod$ induced by the map $f: G \to H$.
Then
\[\scrC^{H} \otimes_{\Rep(H)} \Rep(G) \simeq \wit{\scrC}^{G}.\]
\elem

\bproof
This follows immediately from the naturality
of the invariants functor.
\eproof

\noindent Let $X$ denote the above scheme equipped with the
$\act$ action, and let $\wit{X}$ denote the same
scheme equipped with the $\act_d$ action.
Then $\IndCoh(X)$ and $\IndCoh(\wit{X})$ are both naturally
$\G_m$-categories. On the other hand,
the power map $(-)^d$ induces another 
$\G_m$-category structure on $\IndCoh(X)$
which we denote by $\wit{\IndCoh(X)}$.
Clearly, $\wit{\IndCoh(X)}$ and $\IndCoh(\wit{X})$ 
are equivalent $\G_m$-categories.
It therefore follows from \cref{lem: invariants and maps}
that
\beqn
\label{eqn: Coh is in image}
\IndCoh(\wit{X}/\G_m) \simeq \IndCoh(X/\G_m) \otimes_{\Rep(\G_m)} \Rep(\G_m),
\eeqn
where the tensor product is taken along
the functor, ${\wp}_d := (B(-)^d)^*:\Rep(\G_m) \simeq 
\QCoh(B\G_m) \to \QCoh(B\G_m) \simeq \Rep(\G_m)$.

\subsubsection{Shearing and even actions}
We have the following commutative diagram,
\[\begin{tikzcd}
	{\Rep(\G_m)} & {\Rep(\G_m)} \\
	{\Rep(\G_m)} & {\Rep_{\epsilon}^{\super}(\G_m)}
	\arrow["\mathfrak{p}_2", from=1-2, to=2-2]
	\arrow["{{\wp}_2}"', from=1-1, to=2-1]
	\arrow["{(-)^{2\shear}}", from=1-1, to=1-2]
	\arrow["{(-)^{\shear}}", from=2-1, to=2-2].
\end{tikzcd}\]
where $\mathfrak{p}_2$ is the symmetric monoidal functor 
defined as follows. The composition
of ${\wp}_2$ with the inclusion $\Rep(\G_m)
\subset \Rep^{\super}(\G_m)$ factors
through the inclusion $\Rep_{\epsilon}^{\super}(\G_m) \hook
\Rep^{\super}(\G_m)$. The map $\mathfrak{p}_2$ is defined to
be the first map in this factorization.
This commutative diagram induces a commutative
diagram of categories of modules over its terms:
\[\begin{tikzcd}
	{\Rep(\G_m)\mod} & {\Rep(\G_m)\mod} \\
	{\Rep(\G_m)\mod} & {\Rep_{\epsilon}^{\super}(\G_m)\mod}
	\arrow["{\mathfrak{p}_2 := - \otimes_{\Rep(\G_m)} \Rep_{\epsilon}^{\super}(\G_m)}", from=1-2, to=2-2]
	\arrow["{{\wp}_2 := - \otimes_{\Rep(\G_m)} \Rep(\G_m)}"', from=1-1, to=2-1]
	\arrow["{(-)^{2\shear}}", from=1-1, to=1-2]
	\arrow["{(-)^{\shear}}", from=2-1, to=2-2],
\end{tikzcd}\]
where we abuse notation by reusing
the notation ${\wp}_2$ and $\mathfrak{p}_2$.
Moreover, because all the functors in this
diagram are symmetric monoidal,
if $\sfA$ is an algebra object in $\Rep(\G_m)$,
we obtain, in the usual way, a commutative diagram,
\beqn
\label{eqn: top rectangle}
\begin{tikzcd}
	{\sfA\mod(\Rep(\G_m)\mod)} & {\sfA^{2\shear}\mod(\Rep(\G_m)\mod)} \\
	{\sfA(2)\mod(\Rep(\G_m)\mod)} & {\sfA(2)^{\shear}\mod(\Rep_{\epsilon}^{\super}(\G_m)\mod)}
	\arrow["{\mathfrak{p}_2}", from=1-2, to=2-2]
	\arrow["{{\wp}_2}"', from=1-1, to=2-1]
	\arrow["{(-)^{2\shear}}", from=1-1, to=1-2]
	\arrow["{(-)^{\shear}}", from=2-1, to=2-2],
\end{tikzcd}
\eeqn
where $\sfA(2) := {\wp}_2(\sfA) = \sfA \otimes_{\Rep(\G_m)} \Rep(\G_m)$.

\begin{notation}
${}$
\begin{itemize}
\item Let $k[\a] \simeq C^*(BS^1;k)$ denote the
graded nonconnective derived $k$-algebra where the
generator $\a$ has $\G_m$-weight $1$.

\item Let $k[\b] \simeq C^*(BS^1;k)$ denote the
graded nonconnective derived $k$-algebra
where the generator $\b$ has $\G_m$-weight $2$.
\end{itemize}

Now take $\sfA$ to be $k[\a^{\pm 1}]\mod_{\gr}$. Then
\begin{itemize}
\item  $\sfA(2) \simeq \kb\mod_{\gr}$; 

\item $\sfA^{2\shear} \simeq k[x^{\pm}]\mod_{\gr}$, 
where $k[x^{\pm 1}]$ is the ordinary
graded commutative ring where $x$ has $\G_m$-weight $1$;

\item and $\sfA(2)^{\shear} \simeq k[z^{\pm 1}]\mod(\Rep_{\epsilon}^{\super}(\G_m))
=: k[z^{\pm 1}]\mod_{\gr}$, 
where $k[z^{\pm 1}]$ is the ordinary graded commutative
$k$-algebra, where $z$ has $\G_m$-weight $2$,
viewed as a trivial super $k$-algebra where
everything has even parity.
\end{itemize}
\end{notation}

Observe that $k[x^{\pm 1}]\mod_{gr}$ is equivalent,
as a $\Rep(\G_m)$-module, to $\Vect_k$ with the trivial module
structure. In turn, we have $(k[x^{\pm 1}]\mod_{\gr})\mod(\Rep(\G_m)\mod)
\simeq \dgcat_k^{\infty}$ as symmetric monoidal categories.
Thus, we obtain the following commutative diagram,
\beqn
\label{eqn: biggest diagram}
\begin{tikzcd}
	{(k[\a^{\pm 1}]\mod_{\gr})\mod} & {\dgcat_k^{\infty}} \\
	{(\kb\mod_{\gr})\mod} & {(k[z^{\pm 1}]\mod_{\gr})\mod(\Rep_{\epsilon}^{\super}(\G_m)\mod)} \\
	{(\kb\mod)\mod}
	\arrow["{{\wp}_2}"', from=1-1, to=2-1]
	\arrow["\mathfrak{p}_2", from=1-2, to=2-2]
	\arrow["{(-)^{2\shear}}","\simeq"', from=1-1, to=1-2]
	\arrow["{(-)^{\shear}}","\simeq"', from=2-1, to=2-2]
	\arrow["\oblv"', from=2-1, to=3-1]
	\arrow["\oblv \circ (-)^{\unshear}", from=2-2, to=3-1],
\end{tikzcd}
\eeqn
where the top square of this diagram is
simply the diagram \labelcref{eqn: top rectangle}
for $\sfA = k[\a^{\pm 1}]\mod_{\gr}$. 
We note in the following lemma that
$(k[\a^{\pm 1}]\mod_{\gr})\bmod$
and $\dgcat_k^{\infty}$ embed as full, symmetric
monoidal subcategories of $(\kb\mod_{\gr})\bmod$
and $(k[z^{\pm 1}]\mod_{\gr})\mod(\Rep_{\epsilon}^{\super}(\G_m)\mod)$,
respectively, via the vertical arrows in the top square of
this diagram.

\blem
The functors ${\wp}_2$ and $\mathfrak{p}_2$ are fully
faithful symmetric monoidal embeddings.
\elem

\bproof
To prove the claim, it suffices to show only that one
of either ${\wp}_2$ or $\mathfrak{p}_2$ is fully faithful, since
they form a commutative square with the equivalences
$(-)^{\shear}$ and $(-)^{2\shear}$. In order
to see that ${\wp}_2$ is fully faithful, we observe that
the square power map on $\G_m$ is a central isogeny of
group schemes, so restriction along it induces a fully
faithful functor,
\[\res_{(-)^2}: \Rep(\G_m) \to \Rep(\G_m).\]
Under the equivalence $\QCoh(B\G_m) \simeq \Rep(\G_m)$,
the functor of restriction along the square power map $(-)^2:
\G_m \to \G_m$ corresponds to the pullback of $\QCoh$ along
the map of classifying stacks $B(-)^2: B\G_m \to B\G_m$,
which is therefore fully faithful. 
Unwinding the definition of ${\wp}_2$, we
see that this suffices to show that it is fully faithful.

That $\mathfrak{p}_2$ and ${\wp}_2$ are symmetric monoidal
can be seen from the fact that pullback of $\QCoh$ is symmetric monoidal.
\eproof

In particular, the Hochschild homology
of an object in $\dgcat_k^{\infty}$ is sent under the
embedding $\mathfrak{p}_2$ to its Hochschild homology as an object in
$(\kz\mod_{\gr})\mod(\Rep_{\epsilon}^{\super}(\G_m)\mod)$.


\brem
Unwinding the definitions, it is easy to
see that equivalence \labelcref{eqn: Coh is in image} 
implies that $\MF^{\infty}(M/\G_m,f) := \IndCoh(M_0/\G_m)[\b^{-1}]$ 
is the image of $\MF^{\infty}(\wit{M}/\G_m, f) :=
\IndCoh(\wit{M}_0/\G_m)[\a^{-1}]$
under the functor ${\wp}_2$.
On the other hand, we may identify objects in
$(\ka\mod_{\gr})\mod$ 
and $(\kb\mod_{\gr})\mod$
with their images under shearing, so that
computing the Hochschild homology of $\MF^{\infty}(M/\G_m,f)$
amounts to computing the Hochschild
homology of the \emph{shear} of $\MF^{\infty}(\wit{M}/\G_m,f)$ as
as a plain dg-category (i.e. object of $\dgcat_k^{\infty}$).
\erem

\begin{notation}
We introduce the following notation.
\begin{itemize}
\item Let $\calMF^{\infty}(M/\G_m,f)$ denote $\MF^{\infty}(M/\G_m,f)^{\shear} \in \kz_{\gr}\mod(\Rep_{\epsilon}^{\super}(\G_m))$.
\item Let $\calMF^{\infty}(\wit{M}/\G_m,f)$ denote $\MF^{\infty}(\wit{M}/\G_m,f)^{2\shear} \in \dgcat_k^{\infty}$.
\end{itemize}
\end{notation}

\bprop
\label{prop: simplified HP computation}
Let $(M,f)$ be an evenly graded Landau--Ginzburg pair,
and let $\wit{M}$ denote the same scheme equipped with the
$\G_m$-action whose square gives $M$. Then there is an equivalence,
\[\HP_{\bullet}^k(\calMF^{\infty}(\wit{M}/\G_m, f)) \otimes_{k} \kb \simeq \HP_{\bullet}^{\kb}(\MF^{\infty}(M,f)),\]
as $k[\b^{\pm 1}, u^{\pm 1}]$-modules.
\eprop

\bproof
We begin by showing the corresponding equivalence
for Hochschild homology. \Cref{prop: graded and ungraded PreMF have the same HH}
implies the following equivalence, 
\[\oblv\left(\HH_{\bullet}^{\kz_{\gr}}(\calMF^{\infty}(M/\G_m,f))^{\unshear}\right) \simeq \HH_{\bullet}^{\kb}(\MF^{\infty}(M,f)).\]
On the other hand, it follows from the commutativity of
\labelcref{eqn: biggest diagram} that
\begin{align*}
\oblv\left(\HH_{\bullet}^{\kz_{\gr}}(\calMF^{\infty}(M/\G_m,f))^{\unshear}\right) &\simeq \oblv\left({\wp}_2\left(\HH_{\bullet}^{\ka_{\gr}}(\MF^{\infty}(\wit{M}/\G_m,f))\right)\right) \\
														&\simeq \oblv\left(\mathfrak{p}_2\left(\HH_{\bullet}^k(\calMF^{\infty}(\wit{M}/\G_m,f))\right)^{\unshear}\right) \\
														&\simeq \HH_{\bullet}^k(\calMF^{\infty}(\wit{M}/\G_m,f)) \otimes_k \kb,
\end{align*}
where the last isomorphism is the claim
that the composition $\oblv \circ (-)^{\unshear} \circ \mathfrak{p}_2$
is naturally isomorphic to the $2$-periodization of dg-categories,
$- \otimes_k \kb: \dgcat_k^{\infty} \to (\kb\mod)\mod$, which can
be seen by unwinding the definitions. Finally,
observe that all the equivalences
above preserve the $S^1$-actions on Hochschild
homology, so the proposition follows
by taking Tate invariants.
\eproof

Identifying $\calMF^{\infty}(\wit{M}/\G_m,f)$
with $\MF^{\infty}(M/\G_m,f)$ by abuse of notation, we can state the assertion of
\cref{prop: simplified HP computation} as computing the
equivalence,
\[\HP_{\bullet}^k(\MF^{\infty}(M/\G_m,f)) \otimes_k \kb \simeq \HP_{\bullet}^{\kb}(\MF^{\infty}(M,f)),\]
for an evenly graded Landau--Ginzburg pair. We observe that
the action of $\G_m$ on a $\G_m$-equivariant \'etale open $U \to M$
is not necessarily even. Consequently, the results of this section
on evenly graded Landau--Ginzburg pairs do not pass to presheaves
on $(M/\G_m)_{\et}$. On the other hand, they clearly do sheafify
over the conic Zariski site of $M$ (since the $\G_m$-action on $M$ restricts
to an even $\G_m$-action on every Zariski open of $M$), so as a variant of
\cref{cor: HP of gr/usual MF} we obtain the following corollary.

\bcor
\label{cor: even LG pair HP comparison}
Let $(M,f)$ be an evenly graded Landau--Ginzburg pair.
Then there is an equivalence of sheaves of $k[\b^{\pm 1}, u^{\pm 1}]$-modules
on $(M/\G_m)_{\Zar}$,
\[\unHP_{\bullet}^k(\MF^{\infty}(M/\G_m,f)) \otimes_k \kb \simeq \varpi_*\unHP_{\bullet}^{\kb}(\MF^{\infty}(M,f)).\]
\ecor


\section{Periodic cyclic homology of matrix factorizations}
\label{sec: periodic cyclic homology of matrix factorizations}
In the previous section we related the periodic
cyclic homology of graded matrix factorizations to that of
the usual $2$-periodic matrix factorizations. The purpose of
this computation is to leverage the well-known relationship
between matrix factorizations and vanishing cycles first formulated
by A. Efimov in \cite{Efimov}. In this section, we will prove a
sheafified variant of Efimov's theorem, formulated in terms 
of our present setting.

Our computation of
$\HP_{\bullet}^{\kb}(\unMF^{\infty}(M,w))$,
following Efimov,
will proceed in roughly two parts:
\begin{enumerate}
\item 
\label{item: algebraic}
An algebraic computation identifying $\HP_{\bullet}^{\kb}(\unMF^{\infty}(M,w))$
with an algebraic formal twisted de Rham complex on $M$.
\item 
\label{item: analytic}
An algebraic-to-analytic computation that identifies
this algebraic formal twisted de Rham complex with the vanishing cycles
sheaf with respect to the superpotential $w$.
\end{enumerate}

\begin{conventions} We adopt the following conventions
for the remainder of the section. We take the field
$k$ of the previous sections to be the field of complex numbers,
$\C$. We take $e := \C$ to be the
coefficients of any constructible sheaves that may appear.
We also fix a Landau--Ginzburg pair $(M,w)$ for the
remainder of this section.
\end{conventions}

\subsubsection{}


The computation identifying the periodic cyclic homology
of the matrix factorizations sheaf $\unMF^{\infty}(M,w)$ 
as the formal algebraic twisted de Rham
complex of $w$ can be found in \cite{PreygelT} in the
form of the following proposition.

\bprop[{c.f. \cite[Theorem 8.2.6 (iv)]{PreygelT}}]
\label{prop: no curved dg-algebras}
The HKR isomorphism induces the following 
equivalence of $k[\b^{\pm 1}, u^{\pm 1}]\mod$-valued
presheaves on $M_{\et}$,
\beqn
\label{eqn: Toly de Rham complex}
\HP_{\bullet}^{\kb}(\unMF^{\infty}(M,w)) \simeq 
	\left[\unOmega_M^{\bullet}[\b^{\pm 1}, u^{\pm 1}], 
		\b \cdot (-dw \wedge -) + u \cdot d\right],
\eeqn
where $\b$ and $u$ are each generators
for distinct copies of $C^*(BS^1;k)$.
\eprop

\bproof
This follows essentially immediately from
\cite[Theorem 6.1.2.5 (iv)]{PreygelThesis}
%
\eproof

\subsubsection{}


We will compare the right-hand side of \labelcref{eqn: Toly de Rham complex}
to vanishing cycles using a special case of the author's recent result
in \cite{KennyKSS} generalizing the theorem of Sabbah and Saito
(originally a conjecture of Kontsevich) that hypercohomology
of the formal twisted de Rham complex is given by vanishing cohomology.

In the special case formulated below, $\tau: M_{\an} \to M_{\Zar}$
denotes the canonical continuous map of topologies on $M$.

\bthm[{\cite[Theorem 8.21]{KennyKSS}}]
\label{thm: main theorem from KennyKSS}
Let $M$ be a smooth complex algebraic variety,
and $w: M \to \Aone$ a regular function on $M$.
Assume that the critical locus of $w$ is contained
in the zero fiber of $w$.
Then there is an equivalence of $\wE_{\C,0}^{\poly}\mod$-valued
Zariski sheaves on $M$,
\beqn
\label{eqn: Sabbah de Rham complex}
\left(\left[\unOmega_M^{\bullet}\llp \xi \rrp, (-dw \wedge -) + \xi \cdot d\right], \nabla_{\DR}\right) 
	\simeq \tau_*\muRH_{\Shv}\left(\varphi_w(\C_M), T_f \right),
\eeqn
where $\xi$ is a formal variable in degree $0$, and
where $\nabla_{\DR} = \partial_\xi + w/\xi^2$.
\ethm

\brem
See \cite{KennyKSS} for the
definitions of the various symbols appearing
in \cref{thm: main theorem from KennyKSS}
that have not appeared in the present paper until now.
Their precise meanings are 
irrelevant for our purposes, 
as the reader shall see below.
\erem

In fact, we will not need the full strength of even this
special case, since we consider the periodic cyclic homology
of $\MF(M,w)$ without its standard connection (i.e. without
its full $\wE_{\C,0}^{\poly}$-module structure). The following
corollary follows from \cref{thm: main theorem from
KennyKSS} by restriction along the algebra map, $\C\llp \xi \rrp
\hook \wE_{\C,0}^{\poly}$ and unwinding the definitions.

\bcor
\label{cor: Sabbah's thm without connection}
Let $M$ and $w$ be as in \cref{thm: main theorem from KennyKSS}.
View $\varphi_w(\C_M)$ as a sheaf on $M$ via
the pushforward along the closed inclusion $\crit(w) \hook M$.
Then there is an equivalence,
\[\left[\unOmega_M^{\bullet}\llp \xi \rrp, (-dw \wedge -) + \xi \cdot d\right]
	\simeq \tau_*\varphi_w(\C_M)\llp \xi \rrp,\]
of $\C\llp \xi \rrp\mod$-valued
sheaves on $M_{\Zar}$.
\ecor

%

\brem
The left-hand side of the equivalence in \cref{cor: Sabbah's thm without connection}
denotes the sheaf obtained as follows. 

For each $i$, $\un{\Omega}_M^i$ denotes
the sheaf of degree $i$ algebraic differential forms, viewed as an object in
the abelian category of Zariski sheaves on $M$ valued in $\C$ vector spaces
(i.e. $\Vect_{\C}^{\heart}$), which we denote by $\calA(M;\C)$. 
One obtains a strict complex of sheaves by tensoring each $\un{\Omega}_M^i$ over $\C$
with the ring of formal Laurent series $k\llp \xi \rrp$ and defining 
the differential $(-dw \wedge -) + \xi \cdot d$, where $d$ is the usual
de Rham differential:
\beqn
\label{eqn: strict complex}
\cdots 0 \to \un{\Omega}_M^0\llp \xi \rrp \to \cdots \to \un{\Omega}_M^i\llp \xi \rrp \xrightarrow{(-dw \wedge -) + \xi \cdot d} \un{\Omega}_M^{i+1}\llp \xi \rrp \to \cdots \un{\Omega}_M^{\dim M}\llp \xi \rrp \to 0 \to \cdots
\eeqn
The result is an object the abelian category $\Ch^b(\calA(M; \C\llp \xi \rrp))$,
which we view as an object of the bounded derived $\infty$-category $\calD^b(M;\C\llp \xi \rrp)$
via its image under the localization functor
\beqn
\label{eqn: localization functor}
\Ch^b(\calA(M; \C\llp \xi \rrp)) \to \calD^b(M;\C\llp \xi \rrp).
\eeqn
The bounded derived $\infty$-category $\calD^b(M; \C\llp \xi \rrp)$,
on the other hand, is equivalent to the category of sheaves on $M$
valued in $\calD^b(\C) \simeq \Vect^b_{\C}$. Under this equivalence,
the complex of sheaves \labelcref{eqn: strict complex} is
sent under the localization functor \labelcref{eqn: localization functor} to the sheaf determined by
the assignment,
\begin{align*}
U 	&\mapsto R\Gamma(U; \left[\unOmega_M^{\bullet}\llp \xi \rrp, (-dw \wedge -) + \xi \cdot d\right]) \\
	&\simeq  \left[\unOmega_U^{\bullet}\llp \xi \rrp, (-dw \wedge -) + \xi \cdot d\right]
\end{align*}
where $R\Gamma(U; -)$ denotes the derived functor of
sections over the Zariski open $U$. 
\erem

\subsubsection{}
The two expressions given by the
right-hand side of \labelcref{eqn: Toly de Rham complex}
and the left-hand side of \labelcref{eqn: Sabbah de Rham complex},
respectively, look very similar, but they are not
quite the same. In particular, the object given
by the former expression is a
doubly $2$-periodic complex, whereas the
object given by the latter expression is
a $2$-periodic complex with an action
of the ordinary ring
of formal Laurent series.
The way we explain away this discrepancy
is to observe that we are working with
\emph{complete} modules, and that complete
modules over the rings $k[\b^{\pm 1}, u]$
and $k\llb \xi \rrb[\b^{\pm 1}]$ are identified
under an equivalence of categories.
We explain the details below, starting by
recalling the definition of a complete 
$k[\b^{\pm 1}, u]$-module.

\bdef
A $k[\b^{\pm 1}, u]$-module $M$ is $u$-complete
if the natural map,
\[M \to \lim_{n \geq 0} M/{u}^n,\footnotemark\]
\footnotetext{$M/{u}^n$ denotes the
cofiber in $k[\b^{\pm 1}, u]\mod$ of the 
map given by multiplication by ${u}^n$.}
is an equivalence. 
\edefn

Let $k[\b^{\pm 1}, u]\cmod$ denote the full $\infty$-subcategory
of $k[\b^{\pm 1}, u]\mod$ spanned by $u$-complete modules.
We record the following observation, pointed out to the
author by Pavel Safronov.

\blem
\label{lem: identification of complete modules}
Let $\b$, $u$ denote generators
for two distinct copies of $C^*(BS^1;k)$, and let
$\xi$ denote a formal variable in degree $0$. There is
an equivalence between
$k[\b^{\pm 1}, u]\cmod$ and $k\llb \xi \rrb[\b^{\pm 1}]\cmod$,
where the latter denotes the category of $\xi$-complete modules
over $k\llb \xi \rrb[\b^{\pm 1}]$.
\elem

\bproof
Let $M$ be a $u$-complete module,
which we view as an object in 
\[\lim_{n \geq 0} \left(\left(k[\b^{\pm 1}] \otimes_k k[u]/{u}^n\right) \mod\right).\]
Observe that, for each $n$, we have an isomorphism
of rings 
\[k[u]/{u}^n \otimes_k k[\b^{\pm 1}]
	\xrightarrow{\simeq} k[\xi]/\xi^n \otimes_k k[\b^{\pm 1}],\]
given by sending $u \mapsto \xi\b$ and $\b \mapsto \b$.
Under this identification of rings, we have an equivalence
of categories, $\lim_{n \geq 0} \left(\left(k[\b^{\pm 1}] 
\otimes_k k[u]/{u}^n\right) \mod\right) \simeq  
\lim_{n \geq 0} \left(\left(k[\xi]/\xi^n \otimes_k k[\b^{\pm 1}]\right) \mod\right)$.
The latter category is precisely $k\llb \xi \rrb [\b^{\pm 1}]\cmod$.
\eproof

The operation of inverting $u$ (meaning tensoring
with $k[\b^{\pm 1}, u^{\pm 1}]$) is a functor,
$k[\b^{\pm 1}, u]\cmod \to k[\b^{\pm 1}, u^{\pm 1}]\mod$,
whose essential image we denote by $k[\b^{\pm 1}, u^{\pm 1}]\cmod$.
Similarly, the operation of inverting $u$ is a functor
$k\llb \xi \rrb[\b^{\pm 1}]\cmod \to k\llp \xi \rrp [\b^{\pm 1}]\mod$,
whose essential image we denote by $k\llp \xi \rrp[\b^{\pm 1}]\cmod$.
It is clear that the equivalence exhibited in
\cref{lem: identification of complete modules}
descends to an equivalence,
\beqn
\label{eqn: equivalence after inverting}
k[\b^{\pm 1}, u^{\pm 1}]\cmod \xrightarrow{\simeq} k\llp \xi \rrp[\b^{\pm 1}]\cmod
\eeqn

\bclaim
The $\C[\b^{\pm 1}, u]$-module
$\left[\Omega_M^{\bullet}[\b^{\pm 1}, u], 
\b \cdot (-dw \wedge -) + u \cdot d\right]$
is $u$-complete.
\eclaim

\bproof
This is easy to see by inspection. 
\eproof

Under the equivalence \labelcref{eqn: equivalence after inverting},
the $\C[\b^{\pm 1}, u^{\pm 1}]$-module 
$\left[\Omega_M^{\bullet}[\b^{\pm 1}, u^{\pm 1}], 
\b \cdot (-dw \wedge -) + u \cdot d\right]$
corresponds to the $\C\llb \xi \rrb[\b^{\pm 1}]$-module 
$\left[\Omega_M^{\bullet}\llb \xi \rrb[\beta^{\pm 1}], 
\b \cdot \left((-dw \wedge -) + \xi \cdot d\right)\right]$.
Thus, we may view the periodic cyclic homology
of matrix factorizations as a $2$-periodic module
over the ordinary field $\C\llp \xi \rrp$.

\blem
\label{lem: Tate invariants for twisted dR complex}
Consider $\left[\unOmega_M^{\bullet}\llp \xi \rrp, (-dw \wedge -) + \xi \cdot d\right]$
as a $S^1\mod(\C\llp \xi \rrp\mod)$-valued Zariski presheaf on $M$,
whose sections are all trivial $S^1$-modules. 
Then we have the following equivalence, 
\[\left[\unOmega_M^{\bullet}\llp \xi \rrp, (-dw \wedge -) + \xi \cdot d\right]^{tS^1} 
	\simeq \left[\unOmega_M^{\bullet}\llp \xi \rrp[\beta^{\pm 1}], \b \cdot \left((-dw \wedge -) + \xi \cdot d\right)\right],\]
of $\C\llp \xi \rrp[\beta^{\pm 1}]\mod$-valued Zariski presheaves on $M$.
\elem

\bproof
It suffices to check that the equivalence
holds on sections over a Zariski open $U \subset M$.

The trivial $S^1$-action on $\left[\Omega_U^{\bullet}\llp \xi \rrp, 
(-df \wedge -) + \xi \cdot d\right]$ corresponds under the
localization map $\Ch_{\C\llp \xi \rrp} \to \C\llp \xi \rrp\mod$ 
to a trivial mixed complex structure on $\Omega_U^{\bullet}$,
where both the $b$ and $B$ differentials are given by
$(-df \wedge -) + \xi \cdot d$.
There is a standard formula for computing the $S^1$-invariants
in terms of this mixed complex structure which, in this case, 
yields the complex $\left[\Omega_U^{\bullet}\llp \xi \rrp[\b], 
\b \cdot \left((-df \wedge -) + \xi \cdot d\right)\right]$.
The Tate invariants are then obtained by inverting $\b$,
which obviously gives the desired result.
\eproof

\subsubsection{}

All together, we obtain the following result.

\bprop
\label{prop: fancy Efimov theorem}
Let $M$ be a smooth complex projective variety, 
and let $w: M \to \Aone$ be a regular function such
that the critical locus of $f$ is contained in its zero locus.
Then there is an equivalence,
\[\HP_{\bullet}^{\Cb}(\unMF^{\infty}(M,w)) \simeq
	\tau_*\varphi_w(\C_M)\llp \xi \rrp[\beta^{\pm 1}]\]
of $\C\llp \xi \rrp[\beta^{\pm 1}]\mod$-valued
Zariski presheaves on $M$, supported on $\crit(w)$.
Sheafifying, we obtain an equivalence of Zariski sheaves
on $M$, also supported on $\crit(w)$,
\[\unHP_{\bullet}^{\Cb}(\MF^{\infty}(M,w)) \simeq \tau_*\varphi_w(\C_M)\llp \xi \rrp[\b^{\pm 1}].\]
\eprop

\bproof
By \cref{cor: Sabbah's thm without connection}, the
Tate invariants for the trivial $S^1$-action on $\tau_*\varphi_w(\C_M)\llp \xi \rrp$
is equivalent to the Tate invariants for the trivial $S^1$-action on
$[\unOmega_M^{\bullet}\llp \xi \rrp, (-dw \wedge -) + \xi \cdot d]$.
The latter is equivalent to $[\unOmega_M^{\bullet}\llp \xi \rrp[\b^{\pm 1}], 
\b \cdot((-dw \wedge -) + \xi \cdot d)]$ by \cref{lem: Tate invariants for
twisted dR complex}, which in turns corresponds to 
$\left[\unOmega_M^{\bullet}[\b^{\pm 1}, u^{\pm 1}], 
\b \cdot (-dw \wedge -) + u \cdot d\right]$
under the equivalence \labelcref{eqn: equivalence after inverting}.
We conclude
by observing that this last object is 
equivalent to $\HP_{\bullet}^{\Cb}(\unMF^{\infty}(M,w))$
by \cref{prop: no curved dg-algebras}.
\eproof

\subsubsection{Variant with supports}

Although the pointwise application of the functor 
of periodic cyclic homology, like
Hochschild homology, does not preserve sheaves, 
we nonetheless are able to find a nice expression for
$\HP_{\bullet}^{\Cb}(\MF^{\infty}_Z(M,w))$.

\blem
Let $Z$ be a Zariski closed
subset of $M$, and let $U = M \setminus Z$
denote its complementary open.
Then there is an equivalence
of $\C\llp \xi \rrp [\b^{\pm 1}]$-modules
\[\HP_{\bullet}^{\Cb}(\MF^{\infty}_Z(M,w)) \simeq 
	\Gamma_Z\left(Z(w);\varphi_w(\C_M)\right) \otimes_\C \C\llp \xi \rrp [\b^{\pm 1}].\]
\elem

\bproof
We recall that the Hochschild homology $\HH_{\bullet}^{\Cb}(-)$
is a localizing invariant, meaning that it sends localization
sequences in $(\Cb\mod)\mod$ to fiber sequences in $S^1\mod(\Cb\mod)$.
On the other hand, the Tate fixed points, as the cofiber of a map
between exact functors (namely, $S^1$-coinvariants and $S^1$-invariants), is also
exact, and therefore preserves fiber sequences. It follows that
the periodic cyclic homology sends localization sequences to fiber
sequences, and, in particular, we have that
\beqn
\label{eqn: French revolution}
\HP_{\bullet}^{\Cb}(\MF^{\infty}_Z(M,w)) \simeq 
	\fib\left(\HP_{\bullet}^{\Cb}(\MF^{\infty}(M,w)) \to \HP_{\bullet}^{\Cb}(\MF^{\infty}(U,w|_U)\right).
\eeqn
The right-hand side of \labelcref{eqn: French revolution} can be rewritten
using \cref{prop: fancy Efimov theorem} as
\beqn
\label{eqn: Google fiber}
\fib\left(\Gamma\left(M;\varphi_w(\C_M)\llp \xi \rrp [\b^{\pm 1}]\right) 
	\to \Gamma\left(U;\varphi_w(\C_M)\llp \xi \rrp [\b^{\pm 1}]\right)\right).
\eeqn
Because $- \otimes_{\C} \C\llp \xi \rrp [\b^{\pm 1}]_M$ preserves finite
limits and colimits (being a functor of stable $\infty$-categories), it follows that
\begin{align*}
\labelcref{eqn: Google fiber} &\simeq \Gamma_Z\left(M; \varphi_w(\C_M)\llp \xi \rrp [\b^{\pm 1}]\right) \\
						&\simeq \Gamma_Z\left(M; \varphi_w(\C_M)\right) \otimes_{\C} \C\llp \xi \rrp [\b^{\pm 1}].
\end{align*}
\eproof

\brem
The above corollary is a variant of \cite[Theorem 5.4]{Efimov} 
in which one allows various support
conditions (and disregards connections). Efimov wrote in \cite[\S6]{Efimov} that
such a variation holds, but did not produce a proof therein. 
\erem


\section{The main theorem and applications}
\label{sec: the main theorem and applications}


We now state the main theorem of this paper,
whose proof consists mainly in collecting and
fitting together the results of the previous sections.

\bthm
\label{thm: main theorem}
Let $\Z$ be a derived global complete
intersection over $k = \C$, presented as
the derived zero fiber of a map $f: X \to V$,
where $X$ has dimension $m$ and $V$ dimension
$n$. Let $\tau:
(T^*[-1]\Z)_{\an} \to (T^*[-1]\Z/\G_m)_{\Zar}$
denote the canonical map from the analytic topology 
to the $\G_m$-invariant Zariski topology.
There is an equivalence,
\[\unHP_{\bullet}^k(\IndCoh(\Z))[\b^{\pm 1}][m+n] \simeq \tau_*\varphi_{\T^*[-1]\Z}\llp \xi \rrp[\b^{\pm 1}],\]
of $\C\llp \xi \rrp[\b^{\pm 1}]\mod$-valued sheaves\footnotemark
\footnotetext{See \cref{conv: pointwise presheaf}.}
on $(T^*[-1]\Z/\G_m)_{\Zar}$.
\ethm

\bproof
By the $1$-affineness of $X \times V^{\vee}/\G_m$,
the $\QCoh(X \times V^{\vee}/\G_m)$-module $\IndCoh(\Z)$
is equivalent data to that of a sheaf $\un{\IndCoh}(\Z) \in
\ShvCat(X \times V^{\vee}/\G_m)$. By \cref{cor: sheafy Isik's theorem},
$\un{\IndCoh}(\Z) \simeq \unMF^{\infty}(X \times V^{\vee}/\G_m)$ as
objects of $\ShvCat(X \times V^{\vee}/\G_m)$.
Objects of $\ShvCat(X \times V^{\vee}/\G_m)$
satisfy fppf descent by \cite[Corollary 1.5.4]{Gaitsgory1affineness}, 
so, in particular, they are Zariski sheaves,
i.e. objects of 
$\Shv((X \times V^{\vee}/\G_m)_{\Zar}; \dgcat_k^{\infty})$,
and therefore equivalent as such. 


Moreover, since each is supported on 
$T^*[-1]\Z/\G_m \subset X \times V^{\vee}$, we may view them 
as Zariski sheaves $T^*[-1]\Z/\G_m$ under the
adjoint equivalence of categories,
\[\begin{tikzcd}
	{i_*: \Shv\left((T^*[-1]\Z/\G_m)_{\Zar}; \dgcat_k^{\infty}\right)} && {\Shv\left((X \times V^{\vee}/\G_m)_{\Zar}; \dgcat_k^{\infty}\right)_{T^*[-1]\Z}: i^*}
	\arrow[shift right=1, from=1-1, to=1-3]
	\arrow["\simeq"', shift right=1, from=1-3, to=1-1].
\end{tikzcd}\]

Applying $\unHP_{\bullet}^k(-)$, we obtain the following equivalence
of $k[u^{\pm 1}]\mod$-valued sheaves on $(T^*[-1]\Z/\G_m)_{\Zar}$,
\[\unHP_{\bullet}^k(\IndCoh(\Z)) \xrightarrow{\simeq} 
	\unHP_{\bullet}^k\left(\MF^{\infty}(X \times V^{\vee}/\G_m; \wit{f})\right),\]
from which we deduce the equivalence of 
$k[\b^{\pm 1}, u^{\pm 1}]\mod$-valued presheaves,
\[\unHP_{\bullet}^k(\IndCoh(\Z))[\b^{\pm 1}] \simeq 
	\varpi_*\unHP_{\bullet}^{\kb}\left(\MF^{\infty}(X \times V^{\vee}, \wit{f})\right),\]
using \cref{cor: even LG pair HP comparison}, relying on the fact that
$(X \times V^{\vee}/\G_m, \wit{f})$ is an evenly
graded Landau--Ginzburg pair.

Since sections of $\unHP_{\bullet}^k(\MF^{\infty}(X \times V^{\vee}, \wit{f}))$
are complete $k[\b^{\pm 1}, u^{\pm 1}]$-modules
(by computations done in the previous section), the equivalence
above may be identified, by \cref{lem: identification of complete modules}, 
with a corresponding equivalence 
of $k\llp \xi \rrp[\b^{\pm 1}]$-modules. 
By \cref{prop: fancy Efimov theorem}, the sheaf
of $k\llp \xi \rrp[\b^{\pm 1}]$-modules
$\unHP_{\bullet}^{\kb}(\MF^{\infty}(X \times V^{\vee}, \wit{f}))$ is
equivalent to $\tau_*\varphi_{\wit{f}}(\C_{X \times V^{\vee}})\llp \xi \rrp[\beta^{\pm 1}]$.
We conclude by noting that the 
perverse sheaf $\varphi_{\T^*[-1]\Z}$ is
equivalent to $\varphi_{\wit{f}}(\C_{X \times V^{\vee}}[m+n])$.
\eproof

\subsubsection{}
From our main theorem we obtain the following
corollary which provides a nice description of the
periodic cyclic homology of subcategories of ind-coherent
sheaves with prescribed singular support in terms of sections
of the twisted vanishing cycles sheaf.

\bcor
\label{cor: HP(IndCoh_Lambda)}
Let $\Z$ be a derived global complete intersection
over $k = \C$. Then, for any Zariski-closed, 
conic subset $\Lambda \subseteq T^*[-1]\Z$, 
there exists an equivalence,
\beqn
\label{eqn: HP(IndCoh_Lambda)}
\HP_{\bullet}^k(\IndCoh_{\Lambda}(\Z)) \otimes_{\C} \C[\b^{\pm 1}] \simeq 
	\left(\Gamma_{\Lambda}\varphi_{\T^*[-1]\Z}\right) \otimes_{\C} \C\llp \xi \rrp[\b^{\pm 1}],
\eeqn
where $\Gamma_{\Lambda}$ denotes the functor of
global sections with support contained in $\Lambda$. 
\ecor

\bproof
This follows by applying the functor $\Gamma_{\Lambda}$
to the equivalence of sheaves in \cref{thm: main theorem}
and rewriting the resultant left- and right-hand sides.

Indeed,
as a consequence of the projection formula, tensor
product with the constant sheaves 
$\un{\kb}$ and $\un{k\llp \xi \rrp[\b^{\pm 1}]}$
commutes with the functor $\Gamma_{\Lambda}$,
so we have 
\[\Gamma_{\Lambda}(\varphi_{\T^*[-1]\Z}\llp \xi \rrp [\b^{\pm 1}]) 
	\simeq \left(\Gamma_{\Lambda}\varphi_{\T^*[-1]\Z}\right) \otimes_{\C} \C\llp \xi \rrp [\b^{\pm 1}],\]
which is the right-hand side of \labelcref{eqn: HP(IndCoh_Lambda)}.
By the same reasoning, we have
\beqn
\label{eqn: wacky 2}
\Gamma_{\Lambda}(\unHP_{\bullet}^k(\IndCoh(\Z))[\b^{\pm 1}]) 
	\simeq \Gamma_{\Lambda}(\unHP_{\bullet}^k(\IndCoh(\Z))) \otimes_{\C} \C[\b^{\pm 1}].
\eeqn

Now note that $\un{\IndCoh}(\Z)$ was defined in a such
a way that $\Gamma_{\Lambda}(\un{\IndCoh}(\Z)) \simeq
\IndCoh_{\Lambda}(\Z)$. Since $\HP_{\bullet}^k(-)$
takes localization sequences to fiber sequences, and
since sheafification is exact, it follows that,
\[
\Gamma_{\Lambda}(\unHP_{\bullet}^k(\IndCoh(\Z))) 
	\simeq \HP_{\bullet}^k(\IndCoh_{\Lambda}(\Z)),
\]
which, combined with \labelcref{eqn: wacky 2},
yields the left-hand side of \labelcref{eqn: HP(IndCoh_Lambda)}. 
\eproof

\subsection{Microlocal homology}
Using \cref{cor: HP(IndCoh_Lambda)},
we obtain a beautiful connection between
the categories of ind-coherent sheaves
with prescribed singular support and
Nadler's microlocal homology.

The microlocal homology spaces are
invariants defined for a global complete intersection
variety $Z$ presented as the zero fiber of a
map $f: X \to V$ from a smooth scheme
$X$ to a finite dimensional vector space $V$. 
These invariants are indexed by the family of
closed conic subsets of $V^{\vee}$ and interpolate
between the Borel--Moore homology
and singular cohomology of $Z$: the space associated
to the conic subsets $\{0\}, V^{\vee} \subset V^{\vee}$, respectively,
is $H^{\bullet}(Z;\C)$, $H_{\bullet}^{\BM}(Z;\C)$. 
We reproduce the precise definition below.

\bdef[David Nadler]
\label{def: microlocal homology}
Let $f: X \to V$ be a function from
a smooth scheme $X$ of dimension $n$
to a vector space $V$ of dimension $m$.
Let $\Lambda \subseteq V^{\vee}$
be a closed conic subset. Then the 
\emph{microlocal homology} of $f$
with support $\Lambda$ is
defined as,
\[H_{\bullet}^{\Lambda}(f) := \Gamma_{\Lambda}
	(V^{\vee}; \mu_{0/V}(f_*\C_{X}[m])),\]
where $\Gamma_{\Lambda}(V^{\vee};-)$ is
the global sections of the ``sections with support in
$\Lambda$" functor, using the 
terminology of \cite{KS90}---alternatively, 
the functor ${\pt}_*(\Lambda \hook V^{\vee})^!$. 
\edefn

Let $\Z(f)$ denote the derived zero fiber of $f$.
In previous work (\cite{Kenny22}), we defined a
functor of microlocalization\footnotemark 
\footnotetext{akin to the functor of microlocalization 
along smooth submanifolds found in e.g. \cite{KS90}.}
along quasi-smooth closed immersions of derived
schemes which refines Nadler's microlocal homology
when taken along the quasi-smooth closed immersion
$\Z(f) \hook X$.

More precisely, for any quasi-smooth 
closed immersion $\Y \hook \X$ of derived
schemes over $\C$, and for any sufficiently well-behaved
coefficient ring $e$, we defined a functor
\[\mu_{\Y/\X}: \Shv_c(\X; e) \to \Shv_c(\bmN^{\vee}_{\Y/\X}; e)\] 
from the category of constructible $\Mod_e$-valued sheaves
on $\X$ to the category of constructible $\Mod_e$-valued 
sheaves on the derived conormal
bundle of $\Y \hook \X$. In the case when $\Y \hook \X$ is the inclusion
$\Z(f) \hook X$, the derived conormal bundle $\bmN^{\vee}_{\Z/X}$
is canonically equivalent to the trivial bundle $\Z(f) \times V^{\vee}$. 
In this case, global sections of $\mu_{\Z/X}(\C_X)$
with support recovers the microlocal homology of $f$ when $f$ is
proper, as the following proposition shows.

\bprop[{\cite[Theorem 6.17]{Kenny22}}]
Let $f$ be a proper map, and 
let $\Lambda \subset V^{\vee}$ be a closed conic subset.
Let $\wit{\Lambda} \subset Z(f) \times V^{\vee}$ be any
closed conic subset that $\pr_2(\wit{\Lambda}) = \Lambda$. 
Then there exists an equivalence,
\[\Gamma_{\wit{\Lambda}}(\mu_{\Z(f)/X})
	\simeq H_{\bullet}^{\Lambda}(f)[-n].\]
\eprop

Using \cite[Theorem 6.1]{Kenny22}
and \cref{cor: HP(IndCoh_Lambda)},
we obtain the following immediate corollary, which
provides a family
of statements that interpolate between 
the Feigin--Tsygan theorem relating $\HP(\QCoh(X))$ 
and $H_{\mathrm{dR}}^{\bullet}(X)$, and Preygel's theorem  
relating $\HP(\IndCoh(X))$ and $H^{\mathrm{BM}}_{\bullet}(X)$.

\begin{cor}
\label{cor: interpolating}
Let $\Z$ be a proper derived global complete intersection
over $k = \C$. Then
\[\HP_{\bullet}^k(\IndCoh_{\wit{\Lambda}}(\Z))[\b^{\pm 1}] \simeq 
	(H_{\bullet}^{\Lambda}(f)\llp \xi \rrp [\b^{\pm 1}])[m].\]
\end{cor}

\begin{rem}
Taking $\Lambda = 0$ in the statement of
\cref{cor: interpolating}
yields the Feigin--Tsygan theorem, while taking
$\Lambda = V^{\vee}$ yields Preygel's theorem.
In this way, \cref{cor: interpolating} can be
seen as ``sheafifying" periodic cyclic homology
of $\IndCoh(\Z(f))$ over the conic $-1$-shifted
cotangent bundle of $\Z(f)$, thereby interpolating
between the two results.
\end{rem}

\subsection{Enhanced periodic cyclic homology speculation}
A stronger formulation of \cref{thm: main theorem} which
we believe to be true but do not show is as follows. 

Note that $\IndCoh(\Z)$ is a compactly generated 
$\QCoh(X \times V^{\vee}/\G_m)$-module,
and as such is a dualizable object in $\QCoh(X \times V^{\vee}/\G_m)\mod$.
Its trace, denoted $\HH_{\bullet}^{\QCoh(X \times V^{\vee}/\G_m)}(\IndCoh(\Z))$,
is an object of $S^1\mod(\QCoh(\calL(X \times V^{\vee}/\G_m))$. There is a natural map
$\pi: (\calL (X \times V^{\vee}))/\G_m \to \calL(X \times V^{\vee}/\G_m)$.
Pullback along $\pi$ induces a map,
\[S^1\mod\left(\QCoh(\calL(X \times V^{\vee}/\G_m))\right) \xrightarrow{\pi^*} S^1\mod\left(\QCoh(\calL(X \times V^{\vee})/\G_m)\right).\]
Taking $S^1$-Tate fixed points of an object of
$S^1\mod(\QCoh(\calL(X \times V^{\vee}))$
gives an object of $\QCoh(\calL(X \times V^{\vee}))^{tS^1}$,
and this latter category is equivalent by \cite[Corollary 1.15]{BZNLoops}
to the category of $2$-periodized $\scrD$-modules on $X \times V^{\vee}$,
$\scrD_{X \times V^{\vee}}\mod \otimes_k k[\b^{\pm 1}]$.
Meanwhile, taking $S^1$-Tate fixed points of an object
of $S^1\mod(\QCoh(\calL(X \times V^{\vee})/\G_m))$
will produce a weakly $\G_m$-equivariant (also known as
a $\G_m$-monodromic) $2$-periodized $\scrD$-module on $X \times V^{\vee}/\G_m$.
Assuming a suitable formulation of the Riemann--Hilbert
correspondence for $2$-periodized $\scrD$-modules, we
posit the following conjecture.

\begin{conj}
\label{conj: stronger main theorem}
Let $\DR_{\per}: \scrD_{X \times V^{\vee}}[\b^{\pm 1}]\mod_{\rh} \to
\Shv^b_c(X \times V^{\vee}; \C[\b^{\pm 1}])$ denote the $2$-periodized
de Rham functor. Then
\[\DR_{\per}\left(\left(\pi^*\HH_{\bullet}^{\QCoh(X \times V^{\vee}/\G_m)}(\IndCoh(\Z))\right)^{tS^1}\right)
	\simeq \varphi_{\wit{f}}(\C_{X \times V^{\vee}})[\b^{\pm 1}].\]
\end{conj}

\brem
Note that \cref{thm: main theorem} can be deduced
from \cref{conj: stronger main theorem} by tensoring
up with $\C[u^{\pm 1}]$, taking global sections of each side,
and sheafifying the resulting presheaf.
\erem


\appendix


\section{Primer on $\IndCoh$}
\label{sec: primer on IndCoh}
This appendix is a brief review of the
portions of the theory of ind-coherent sheaves
needed in the main body of this work. The content,
if not the presentation, of this section is taken
from \cite{GaitsgoryIndCoh} and \cite{AG15}.
We use much of the content of this section
freely throughout the work without specific
reference.


\subsection{Basic definition for schemes}
Let $S$ be a Noetherian derived scheme in
the sense of \cite{GaitsgoryIndCoh}.

\bdef
Let $\Coh(S) \subset \QCoh(S)$ denote the full subcategory 
spanned by objects of bounded cohomological amplitude 
(with respect to the natural t-structure on $\QCoh(S)$) and 
coherent (over $\O_{\pi_0(S)}$) cohomology sheaves. 
\edefn

The assumption that $X$ is Noetherian guarantees
that $\Coh(S)$ is a stable subcategory of $\QCoh(S)$.

\bdef
\label{def: IndCoh for Noetherian schemes}
Let $\IndCoh(S) := \Ind(\Coh(S))$ denote the $k$-linear
(stable cocomplete) $\infty$-category of ind-objects of $\Coh(S)$.
This is known as the category of ``ind-coherent sheaves."
\edefn

\brem
It is true, but not obvious, that $\Coh(S)$ is
the subcategory of \emph{all} compact objects
of $\IndCoh(S)$ (\cite[Corollary 1.2.6]{GaitsgoryIndCoh}).
\erem

\subsection{Basic properties}
\subsubsection{}
By construction, we have a canonical morphism
in $\dgcat_k^{\infty}$,
\[\Psi_S: \IndCoh(S) \to \QCoh(S).\]
It is t-exact by \cite[Lemma 1.2.2]{GaitsgoryIndCoh}.

\subsubsection{}
$\Coh(S)$ carries a natural t-structure, so $\IndCoh(S)$
acquires a canonical t-structure characterized by the properties
that
\begin{itemize}
\item[$-$] it is compatible with filtered colimits, and
\item[$-$] the tautological embedding $\Coh(S) \to \IndCoh(S)$
is t-exact.
\end{itemize}

\subsubsection{}
If $S$ is eventually coconnective, $\Psi_S$ admits a
fully faithful left adjoint,
\[\Xi_S: \QCoh(S) \to \IndCoh(S),\]
which realizes $\QCoh(S)$ as a colocalization
of $\IndCoh(S)$ with respect to the full subcategory
of objects which are infinitely connective (\cite[Proposition 1.5.3]{GaitsgoryIndCoh}).
In fact, $\Psi_S$ admits a left adjoint \emph{if and only if}
$S$ is eventually coconnective (\cite[Proposition 1.6.2]{GaitsgoryIndCoh}).

\subsubsection{}
$\QCoh(S)$ has a natural symmetric monoidal structure,
and $\IndCoh(S)$ is naturally a module over $\QCoh(S)$
as follows. There is an obvious action of $\Perf S$ on $\QCoh(S)$,
and this action in fact preserves the subcategory $\Coh(S) \subset \QCoh(S)$.
This gives an action of $\Perf S$ on $\Coh(S)$, which furnishes
a natural $\QCoh(S)$-action on $\IndCoh(S)$ after passing
to ind-completions.

With this $\QCoh(S)$-module structure on $\IndCoh(S)$ and the tautological
one on $\QCoh(S)$, $\Psi_S$ and $\Xi_S$ are naturally functors of $\QCoh(S)$-module
categories (\cite[Lemma 1.4.2, Corollary 1.5.5]{GaitsgoryIndCoh}).

\subsection{Basic functorialities}
Let $f: S \to S'$ be a morphism of Noetherian
derived schemes. Depending on the properties of
the map $f$, there are several functors one
can easily define between $\IndCoh(S)$ and $\IndCoh(S')$.

\begin{outline}
\1 For arbitrary $f$, there is a functor,
\[f_*^{\IndCoh}: \IndCoh(S) \to \IndCoh(S').\]

\1 If $f$ is almost of finite type, there is
the ``extraordinary" pullback functor,
\[f^!: \IndCoh(S') \to \IndCoh(S).\]
	\2 If $f$ is proper, $(f_*^{\IndCoh}, f^!)$
	form an adjunction pair. 

\1 If $f$ is of bounded Tor dimension (= eventually coconnective),
the ordinary pullback of quasicoherent sheaves
preserves coherent complexes, so by ind-extension
we obtain the ``ordinary" pullback functor,
\[f^{\IndCoh, *}: \IndCoh(S') \to \IndCoh(S),\]
which is always left adjoint to $f_*^{\IndCoh}$.
	\2 If $f$ is an open embedding, $f^! = f^{\IndCoh,*}$.
\end{outline}

\subsubsection{Upgrading to a functor}
\label{sssec: upgrading to a functor}
One has several choices of how to upgrade the assignment
$S \mapsto \IndCoh(S)$
to a functor on $\Sch_k^{\Noeth}$ using the functors
$f_*^{\IndCoh}$, $f^!$, and $f^{\IndCoh,*}$:
\[\begin{tikzcd}
	{\IndCoh^!: ((\Sch_k^{\Noeth})_{\aft})^{\op}} & {\dgcat_k^{\infty}} && {S \mapsto \IndCoh(S),\; f \mapsto f^!} \\
	{\IndCoh: \Sch_k^{\Noeth}} & {\dgcat_k^{\infty}} && {S \mapsto \IndCoh(S),\; f \mapsto f_*^{\IndCoh}} & {} \\
	{\IndCoh^*: ((\Sch_k^{\Noeth})_{\Tor})^{\op}} & {\dgcat_k^{\infty}} && {S \mapsto \IndCoh(S),\; f \mapsto f^{\IndCoh,*}}
	\arrow[from=1-1, to=1-2]
	\arrow[from=2-1, to=2-2]
	\arrow[from=3-1, to=3-2]
\end{tikzcd}\]

where $(\Sch_k^{\Noeth})_{\aft}$
and $(\Sch_k^{\Noeth})_{\Tor}$
denote the categories of Noetherian derived $k$-schemes with, respectively,
morphisms almost of finite type and of bounded Tor dimension.

\brem
Note that ${\IndCoh^!}|_{(\Sch_k^{\Noeth})_{\open}} \simeq {\IndCoh^*}|_{(\Sch_k^{\Noeth})_{\open}}$,
where $(\Sch_k^{\Noeth})_{\open}$ denotes the category of Noetherian
schemes with morphisms given by open embeddings.
\erem

\subsubsection{}
The assignment $S \mapsto \QCoh(S)$ has a similar enhancement
to a functor $\QCoh^*: \Sch_k^{\Noeth} \to \dgcat_k^{\infty}$
which sends $f: S \to S'$ to $f^*: \QCoh(S') \to \QCoh(S)$.
By abuse of notation, we let $\QCoh^*$ denote its restriction
to $(\Sch_k^{\Noeth})_{\Tor}$.
The functor $\Psi_S$ of the previous section then
upgrades to a natural transformation,
\[\Psi: \IndCoh^* \to \QCoh^*.\]
If we restrict to the full subcategory $({^{<\infty}\Sch_k^{\Noeth}})_{\Tor}
\subset (\Sch_k^{\Noeth})_{\Tor}$ spanned by eventually
coconnective Noetherian schemes,
$\Xi_S$ also upgrades to a natural transformation,
\[\Xi: {\QCoh^*}|_{({^{<\infty}\Sch_k^{\Noeth}})_{\Tor}}
	\to {\IndCoh^*}|_{({^{<\infty}\Sch_k^{\Noeth}})_{\Tor}}.\]

\subsubsection{$\IndCoh$ as a functor on correspondences}
With more technical machinery at one's disposal, 
it is possible to obtain the functors
$\IndCoh^!$ and $\IndCoh$ above by restriction
from a certain functor out of the category of correspondences
in $\Sch_{\Noeth}$ with arbitrary vertical morphisms
and horizontal morphisms\footnotemark
\footnotetext{Let $S \xleftarrow{} S'' \to S'$ be a correspondence
from $S$ to $S'$. In this example, the 
vertical morphism is the arrow $S'' \to S$,
and the horizontal morphism is the arrow $S'' \to S'$.}
almost of finite type.
More precisely, let $\Corr(\Sch_k^{\Noeth})_{\text{all:aft}}$
denote the aforementioned category of correspondences.
There are natural inclusions $(\Sch_k^{\Noeth})_{\aft} \hook
\Corr(\Sch_k^{\Noeth})_{\text{all:aft}}$ and $\Sch_k^{\Noeth}
\hook \Corr(\Sch_k^{\Noeth})_{\text{all:aft}}$,
heuristically given, respectively, by the assignments,
\[\begin{tikzcd}[row sep=0]
	{\left(S \xrightarrow{f} S'\right)} & {\left(S \xleftarrow{\id} S \xrightarrow{f} S'\right)\text{{} and}} \\
	{\left(S \xrightarrow{f} S'\right)} & {\left(S' \xleftarrow{f} S \xrightarrow{\id} S\right)}
	\arrow[maps to, from=1-1, to=1-2]
	\arrow[maps to, from=2-1, to=2-2].
\end{tikzcd}\]
By \cite[Theorem 5.2.2]{GaitsgoryIndCoh}, there exists a functor,
\[\IndCoh_{\Corr}: \Corr(\Sch_k^{\Noeth})_{\text{all:aft}} \to \dgcat_k^{\infty}\]
such that ${\IndCoh_{\Corr}}|_{(\Sch_k^{\Noeth})_{\aft}} \simeq \IndCoh^!$
and ${\IndCoh_{\Corr}}|_{\Sch_k^{\Noeth}} \simeq \IndCoh$, canonically.

\brem
A salient feature of this formulation of $\IndCoh$ is that
it encodes base change into the very definition
of $\IndCoh$. For details, see \cite[\S 5]{GaitsgoryIndCoh}.
\erem

\brem
\label{rem: IndCoh correspondences}
Alternatively, we may eliminate the restrictions on which
correspondences we consider at the cost of restricting instead to
Noetherian schemes with are almost of finite type. Arbitrary maps
between such schemes are automatically almost of finite type,
so $\Corr(\Sch_{\aft}^{\Noeth})_{\mathrm{all:all}}$ is a full
subcategory of $\Corr(\Sch_k^{\Noeth})_{\mathrm{all:aft}}$ in
which the roles of the horizontal and vertical arrow of a correspondence
are interchangeable. Because of this feature, it will be useful to consider the functor
\[\IndCoh_{\Corr}: \Corr(\Sch_{\aft}^{\Noeth})_{\mathrm{all:all}} \to \dgcat_k^{\infty}\]
in the section below on Serre duality.
\erem

\subsection{Symmetric monoidal structure}
\label{ssec: symmetric monoidal structure}
Let $S$ be a Noetherian derived scheme.
Define the \emph{dualizing sheaf}
$\omega_S := \pt^!k$, where $\pt:
S \to \ast$ is the terminal map to a point.

There is a natural symmetric monoidal structure
on $\IndCoh(S)$ for which $\omega_S$ is the unit. We denote
this symmetric monoidal structure by $\overset{!}{\otimes}$.

\subsubsection{}
The symmetric monoidal structure on $\IndCoh(S)$ is
compatible with its $\QCoh(S)$-module structure.

\subsubsection{}
Among other facts, we note that $f^!$ has a natural symmetric monoidal
structure with respect to $\overset{!}{\otimes}$.
In fact, the functor $\IndCoh^!$ mentioned earlier
upgrades to a functor taking values in $\CMon(\dgcat_k^{\infty})$,
which we denote by
\[(\IndCoh^!, \overset{!}{\otimes}): 
	(\Sch_k^{\Noeth})_{\aft} \to \CMon(\dgcat_k^{\infty}).\]

\subsubsection{}
There is a symmetric monoidal functor,
\[\Upsilon_S: (\QCoh(S), \otimes) \to (\IndCoh(S), \overset{!}{\otimes})\]
given by the assignment $\F \mapsto \F \otimes \omega_S$.
In fact, there is a natural transformation of functors,
\[(\QCoh^*, \otimes) \xrightarrow{\Upsilon} (\IndCoh^!, \overset{!}{\otimes}): 
	(\Sch_k^{\Noeth})_{\aft} \to \CMon(\dgcat_k^{\infty}).\]
For a given Noetherian derived scheme $S$, the functor
$\Upsilon_S$ has a natural $\QCoh(S)$-linear structure.

\brem
If $S$ is eventually coconnective, $\Upsilon_S$ is fully faithful.
\erem

\subsection{$\IndCoh$ on prestacks}
A general prestack (object of $\PreStk$) is simply
a functor in $\Fun(\Alg_k, \Grpd_{\infty})$. $\IndCoh$ is
most easily defined for those prestacks which are \emph{locally almost
of finite type}. In essence, these are prestacks determined
by left Kan extension from their restrictions to the full subcategory
\[^{<\infty}\Alg_k^{\aft} \subset \Alg_k\] 
spanned by eventually coconnective $k$-algebras 
which are almost of finite type. In fact, we may
take this as a definition: $\PreStk_{\laft} := \Fun(^{<\infty}\Alg_k^{\aft}, \Grpd_{\infty})$.

Note that the spectrum of an algebra $A \in {^{<\infty}\Alg_k^{\aft}}$ is,
in particular, an affine Noetherian derived scheme,
so $\IndCoh(\Spec A)$ is defined by \cref{def: IndCoh
for Noetherian schemes} above.
Thus, we define $\IndCoh(X)$ for $X \in \PreStk_{\laft}$ to
be the limit of categories $\IndCoh(S)$ under $!$-pullback
indexed by maps $S \to X$, where $S \in {}^{<\infty}\Sch^{\aft}$.

\bdef
\label{def: IndCoh for prestacks}
Let $X \in \PreStk_{\laft}$. Then we define
\[\IndCoh(X) := \lim_{S \in {^{<\infty}\Sch^{\aft}},\, S \to X} \IndCoh(S).\]
\edefn

\brem
\Cref{def: IndCoh for Noetherian schemes} and \cref{def:
IndCoh for prestacks} agree when the prestack $X$
is representable by a Noetherian derived scheme.
\erem

Note that for a general laft prestack $X$, $\IndCoh(X)$ is \emph{not}
defined as the ind-completion of the subcategory $\Coh(X) \subset \QCoh(X)$.
In many favorable cases---including all the prestacks which appear
in the the body of this work---however, $\IndCoh(X)$ is compactly
generated by $\Coh(X)$.

\brem
In \cite{PreygelT} and \cite{PreygelThesis}, 
``$\on{Ind}\mathrm{DCoh}(X)$"
is used to denote the category defined in \cref{def: IndCoh
for prestacks}, while ``$\QC^!(X)$" is used to denote
the ind-completion of $\Coh(X)$.
\erem

\subsection{Properties of $\IndCoh$ on prestacks}

\subsubsection{Functorialities}
$\IndCoh$ on laft prestacks has the same functorialities
as it does in the case of Noetherian derived schemes
except that pushforward along a map $f$ of prestacks
is only defined if $f$ is schematic quasi-compact
and quasi-separated (\cite[\S 10.6]{GaitsgoryIndCoh}).
We denote the corresponding functors as follows,
\[\begin{tikzcd}
	{\IndCoh^!_{\PreStk}: (\PreStk_{\laft})^{\op}} & {\dgcat_k^{\infty}} && {X \mapsto \IndCoh(X),\; f \mapsto f^!} \\
	{\IndCoh_{\PreStk}: (\PreStk_{\laft})_{\sch}} & {\dgcat_k^{\infty}} && {X \mapsto \IndCoh(X),\; f \mapsto f_*^{\IndCoh}} \\
	{\IndCoh^*_{\Stk}: ((\Stk_{\Artin})_{\Tor})^{\op}} & {\dgcat_k^{\infty}} && {X \mapsto \IndCoh(X),\; f \mapsto f^{\IndCoh,*}}
	\arrow[from=1-1, to=1-2]
	\arrow[from=2-1, to=2-2]
	\arrow[from=3-1, to=3-2]
\end{tikzcd}\]
where $(\Stk_{\Artin})_{\Tor}$ is the
category of laft Artin stacks\footnotemark
\footnotetext{See \cite{GaitsgoryStacks}
for the definition of Artin stack in this context.} 
with morphism of bounded Tor dimension and $(\PreStk_{\laft})_{\sch}$
is the category of laft prestacks with
schematic quasi-compact quasi-separated morphisms.

If we consider, rather than maps between arbitrary prestacks,
maps between a narrower class of prestacks, such as geometric
stacks in the sense of \cite{CautisWilliams}, we are able may pushforward
along maps which are not necessarily schematic. Let $\GStk_{\laft}$ denote
the category of laft geometric stacks\footnotemark
\footnotetext{See \cite[Definition 3.1]{CautisWilliams} for the
definition of geometric stack.}
 with morphisms of finite cohomological
dimension. Then, we have a functor,
\[\begin{tikzcd}
{\IndCoh_{\GStk}: \GStk_{\laft}} & {\dgcat_k^{\infty}} && {X \mapsto \IndCoh(X),\; f \mapsto f_*^{\IndCoh}}
\arrow[from=1-1, to=1-2],
\end{tikzcd}\]
which coincides with the restriction of $\IndCoh_{\PreStk}$ to
subcategory $\GStk_{\laft} \subset (\PreStk_{\laft})_{\sch}$.\footnotemark
\footnotetext{See the discussion after \cite[Definition 5.16]{CautisWilliams}
for details.}

Every stack appearing in the main body of the paper
is evidently geometric.

\subsubsection{Correspondences}
One may formulate define in $\IndCoh$ as a functor
out of a certain category of correspondences of geometric
stacks, just as we had for schemes above. Namely,
by \cite[Definition 5.16]{CautisWilliams}, there is a functor,
\[\IndCoh_{\Corr}: \Corr(\GStk_{\laft})_{\text{fcd:ftd}} \to \dgcat_k^{\infty}\]
whose restriction to $\GStk_{\laft} \subset \Corr(\GStk_{\laft})_{\text{fcd:ftd}}$
is equivalent to $\IndCoh_{\GStk}$, and whose restriction to
$\GStk_{\laft})_{\text{ftd}}^{\op} \subset \Corr(\GStk_{\laft})_{\text{fcd:ftd}}$
is equivalent to $(\IndCoh^!_{\PreStk})|_{\GStk_{\laft}}$. Here $\text{fcd}$
stands for ``finite cohomological dimension" and $\text{ftd}$ stands
for ``finite Tor dimension."

\subsubsection{Symmetric monoidal structure}
Like in the case of schemes, $\IndCoh(X)$ for a laft
prestack $X$ has a natural symmetric monoidal structure,
for which the dualizing sheaf $\omega_X := \pt^!k$ is the
unit. Also like in the case of schemes, there is
an action of $(\QCoh(X),\otimes)$ on $\IndCoh(X)$ which
is compatible with this symmetric monoidal structure.
There is a natural transformation,
\[(\QCoh^*_{\PreStk}, \otimes) \xrightarrow{\Upsilon} (\IndCoh^!_{\PreStk}, \overset{!}{\otimes}):
	\PreStk_{\laft} \to \CMon(\dgcat_k^{\infty}).\]

\subsubsection{}
If $X$ is an Artin stack, there is a functor,
\[\Psi_X: \IndCoh(X) \to \QCoh(X)\]
just as in the case of schemes, with all the expected
properties.
Likewise, if $X$ is an eventually coconnective
Artin stack, $\Psi_X$ has a fully faithful left adjoint,
\[\Xi_X: \QCoh(X) \to \IndCoh(X).\]
These functors upgrade to natural transformations
$\Psi: \IndCoh_{\Stk}^* \to \QCoh_{\Stk}^*$
and $\Xi: \QCoh_{\Stk}^* \to \IndCoh_{\Stk}^*$,
as well.

\subsection{Support}
Ind-coherent sheaves, like quasicoherent sheaves,
have a notion of support.

\bdef
\label{def: support of ind-coherent sheaf}
Let $X$ be a laft Artin stack. The support
of an object $\F \in \IndCoh(X)$ is defined
to be the complement of the smallest open
embedding $j: U \hook X$ such that $j^{\IndCoh,*}\F\simeq 0$.
\edefn

\subsection{Relative enhancement}
The following material is not found in \cite{GaitsgoryIndCoh} or,
in fact, any reference of which we are aware. Nonetheless, with the
exemption of a few claims, which we dutifully note, the assertions
and constructions contained in this section are seen to 
follow easily from straightforward modification of those
found in \textit{op. cit.}

\subsubsection{Relative schemes}
We now formulate an enhancement of the constructions of the
previous section in which our spaces are structured over
a more general base stack, rather than just $\Spec k$.

If $f: Y \to X$ is a map of derived schemes, the dg-category
$\IndCoh(Y)$ obtains a $\QCoh(X)$-module structure by pulling back
its canonical $\QCoh(Y)$-module structure along the
symmetric monoidal functor $f^*: \QCoh(X) \to \QCoh(Y)$.
Moreover, given a map $s: Y \to Y'$ of stacks over
$X$, the functors $s_*^{\IndCoh}$, $s^!$, and $s^{\IndCoh,*}$,
whenever they are defined, are $\QCoh(X)$-linear.

Now fix a base derived stack $X$, and consider the category
$\Sch_X$ of derived stacks which are schemes relative to $X$.
The usual properties of maps between schemes have
their counterparts for maps between stacks, so we may
also consider the full subcategories,
\[\Sch_X^{\Noeth}, \; (\Sch_X^{\Noeth})_{\aft}, \; (\Sch_X^{\Noeth})_{\Tor}, \; \text{etc.} \subset \Sch_X.\] 
Since $\IndCoh(S)$ for $S \in \Sch_X$ is naturally a
$\QCoh(X)$-module, the functors of the previous section
have the follow counterparts relative to $X$:
\[\begin{tikzcd}
	{\IndCoh^!(-/X): ((\Sch_X^{\Noeth})_{\aft})^{\op}} & {\QCoh(X)\mod} && {S \mapsto \IndCoh(S),\; f \mapsto f^!} \\
	{\IndCoh(-/X): \Sch_X^{\Noeth}} & {\QCoh(X)\mod} && {S \mapsto \IndCoh(S),\; f \mapsto f_*^{\IndCoh}} \\
	{\IndCoh^*(-/X): ((\Sch_X^{\Noeth})_{\Tor})^{\op}} & {\QCoh(X)\mod} && {S \mapsto \IndCoh(S),\; f \mapsto f^{\IndCoh,*}}
	\arrow[from=1-1, to=1-2]
	\arrow[from=2-1, to=2-2]
	\arrow[from=3-1, to=3-2].
\end{tikzcd}\]
If $X$ is itself a derived Noetherian $k$-scheme,
these functors are nicely compatible with those
defined above in the sense that we have a commutative diagram,
\[\begin{tikzcd}
	{\Sch_X^{\Noeth}} && {\QCoh(X)\mod} \\
	{\Sch_k^{\Noeth}} && {\dgcat_k^{\infty}}
	\arrow["{\IndCoh(-/X)}", from=1-1, to=1-3]
	\arrow["\IndCoh", from=2-1, to=2-3]
	\arrow["i"', from=1-1, to=2-1]
	\arrow["\oblv", from=1-3, to=2-3]
\end{tikzcd}\]
as well as commutative diagrams for $\IndCoh^!$ and
$\IndCoh^*$, respectively. Here $i: \Sch_X^{\Noeth} \to \Sch_k^{\Noeth}$
is the functor sending an $X$-scheme to its underlying $k$-scheme,
and $\oblv$ is the forgetful functor
taking a $\QCoh(X)$-module to its underlying
dg-category.
On the level of objects, the commutativity of this
diagram translates to an equivalence,
\[\oblv\left(\IndCoh(S/X)\right) \simeq \IndCoh(S).\]

\subsubsection{Relative prestacks}
Fix a base stack $X$, and let $f: Y \to X$ be 
a prestack over $X$. Such data is equivalent
to its functor of points, an object in $\Fun(\Sch_X^{\aff}, \Grpd_{\infty})$,
where $\Sch_X^{\aff}$ denotes the category of prestacks which
are affine schematic over $X$.

In defining $\IndCoh(-/X)$ for such prestacks, 
we restrict ourselves, as above, to prestacks $f: Y \to X$ which are
locally almost of finite type. In essence, these are prestacks whose
functors of points are determined by left Kan extension
from their restrictions to the full subcategory,
\[^{<\infty}\Sch_{\aff/X}^{\aft} \subset \Sch_{\aff/X}.\]
In fact, this may be taken as the definition,
\[\PreStk_{\laft/X} := \Fun(^{<\infty}\Sch_{\aff/X}^{\aft}, \Grpd_{\infty}).\]

The reader should find the following definition well-motivated
by the discussion of the previous section.

\bdef
Let $Y \in \PreStk_{/X, \laft}$. Then we define
\[\IndCoh(Y/X) := \lim_{S \in {^{<\infty}\Sch_X^{\aft}},\, S \to Y} \IndCoh(S/X),\]
where the limit on the right-hand side is taken inside $\QCoh(X)\mod$,
and the morphisms in the index category are maps of prestacks over $X$.
\edefn

The functors $\IndCoh^!(-/X)$, $\IndCoh(-/X)$, and $\IndCoh^*(-/X)$
defined for schemes over $X$ have the expected extensions to the appropriate
categories of stacks and prestacks over $X$:
\[\begin{tikzcd}
	{\IndCoh^!_{\PreStk}(-/X): (\PreStk_{\laft/X})^{\op}} & {\QCoh(X)\mod} && {S \mapsto \IndCoh(S/X),\; f \mapsto f^!} \\
	{\IndCoh_{\PreStk}(-/X): (\PreStk_{\laft/X})_{\sch}} & {\QCoh(X)\mod} && {S \mapsto \IndCoh(S/X),\; f \mapsto f_*^{\IndCoh}} \\
	{\IndCoh^*_{\Stk}(-/X): ((\Stk_{\Artin/X})_{\Tor})^{\op}} & {\QCoh(X)\mod} && {S \mapsto \IndCoh(S/X),\; f \mapsto f^{\IndCoh,*}} \\
	{\IndCoh_{\GStk}(-/X): \GStk_{\laft/X}} & {\QCoh(X)\mod} && {X \mapsto \IndCoh(X),\; f \mapsto f_*^{\IndCoh}}
	\arrow[from=1-1, to=1-2]
	\arrow[from=2-1, to=2-2]
	\arrow[from=3-1, to=3-2]
	\arrow[from=4-1, to=4-2].
\end{tikzcd}\]

\begin{warn}
\label{warn: IndCoh abuse of notation}
By abuse of notation, we will often denote
$\IndCoh(-/X)$ (resp. $\IndCoh^!(-/X)$, $\IndCoh^*(-/X)$) 
by $\IndCoh(-)$ (resp. $\IndCoh^!(-)$, $\IndCoh^*(-)$).
In this case, it should be clear from context whether we are considering
$\IndCoh(-)$ as a module over $\QCoh$ of a base
stack or as a plain dg-category.
\end{warn}

\subsubsection{Symmetric monoidal structure}
If $Y$ is structured over a base stack $X$
with structure morphism $f$,
one has the \emph{relative dualizing sheaf}, $\omega_{Y/X} := 
f^!\O_X$. To our knowledge, there is no natural symmetric monoidal structure
on $\IndCoh(Y)$ for which $\omega_{Y/X}$ is the unit.

On the other hand, the $\QCoh(X)$-module structure on $\IndCoh(Y)$
is actually a $\QCoh(X)$-\emph{algebra} structure induced
by the symmetric monoidal functor,
\[\QCoh(X) \xrightarrow{f^*} \QCoh(Y) \xrightarrow{\Upsilon} \IndCoh(Y).\]
Thus, the symmetric monoidal structure on $\IndCoh(Y)$
is $\QCoh(X)$-linear for any stack and any chosen map $Y \to X$.

\subsubsection{Correspondences}
We mentioned earlier the formulation of $\IndCoh$
as a functor out a certain category of correspondences
which bakes base change into the definition of $\IndCoh$.
We denoted this functor $\IndCoh_{\Corr}$.

Because the functors of pullback and pushforward along maps
of $X$-prestacks are naturally $\QCoh(X)$-linear, the construction
of $\IndCoh_{\Corr}$ is easily seen to give a functor,
\[\IndCoh_{\Corr}(-/X): \Corr(\Sch_X^{\Noeth})_{\text{all:aft}/X} \to \QCoh(X)\mod.\]

Similarly, after observing that all Artin stacks with affine diagonal are geometric stacks,
we see there also exists a functor,
\[\IndCoh_{\Corr, \Stk}(-/X): \Corr(\Stk_{\Artin^{\Delta}/X})_{\text{ftd:fcd}/X} \to \QCoh(X)\mod,\]
where $\Artin^{\Delta}/X$ denotes the category of relative
Artin stacks with affine relative diagonal.

\subsection{Serre duality}
The Serre duality anti-equivalence $\bbD(-)$ on $\IndCoh$
is discussed in \cite[\S 9]{GaitsgoryIndCoh}. We now discuss
such an equivalence in the relative setting by imitating the
presentation therein.

\subsubsection{}
Fix a derived $X$-stack $Y$.
Since $\QCoh(X)$ is rigid symmetric monoidal,
an object of $\QCoh(X)\mod$ is dualizable if and only if
its underlying dg-category is. The underlying dg-category
of $\IndCoh(Y/X)$ is $\IndCoh(Y)$, which is compactly
generated, essentially by construction, and therefore dualizable.
Thus, the essential image of $\IndCoh_{\Corr, \Stk}$ lies
in the full subcategory $\QCoh(X)\mod_{\dual} \subset \QCoh(X)\mod$
spanned by dualizable objects of $\QCoh(X)$.

\subsubsection{}
Let
\[(-)^{\vee}: \QCoh(X)\mod_{\dual} \to (\QCoh(X)\mod_{\dual})^{\op}\]
denote the canonical anti-involution given on object by sending
a category to its $\QCoh(X)$-linear dual.

On the other hand, $\Corr(\Stk_{\Artin^{\Delta}/X})_{\text{ftd:fcd}}$ has an
involution,
\[\Swap: \Corr(\Stk_{\Artin^{\Delta}/X})_{\text{ftd:fcd}} \to \Corr(\Stk_{\Artin^{\Delta}/X})_{\text{ftd:fcd}}^{\op}\]
given by interchanging the roles of the vertical and
horizontal arrows of a correspondence (see \cref{rem: IndCoh correspondences}
above).

The following theorem is the relative counterpart to 
\cite[Theorem 9.1.4]{GaitsgoryIndCoh}, and can be proven
in the same way, \textit{mutatis mutandis}.

\bthm
\label{thm: Serre duality}
The diagram,
\[\begin{tikzcd}
	{\Corr(\Stk_{\Artin^{\Delta}/X})_{\mathrm{ftd:fcd}}} &&& {\QCoh(X)\mod_{\dual}} \\
	{\Corr(\Stk_{\Artin^{\Delta}/X})_{\mathrm{ftd:fcd}}^{\op}} &&& {(\QCoh(X)\mod_{\dual})^{\op}}
	\arrow["{\IndCoh_{\Corr, \Stk}(-/X)}", from=1-1, to=1-4]
	\arrow["{\IndCoh_{\Corr, \Stk}(-/X)^{\op}}", from=2-1, to=2-4]
	\arrow["{(-)^{\vee}}", from=1-4, to=2-4]
	\arrow["\Swap"', from=1-1, to=2-1],
\end{tikzcd}\]
is canonically commutative. Moreover this commutative
structure is canonically compatible with the involutivity
structures on the equivalences $\Swap$ and $(-)^{\vee}$.
\ethm

\subsubsection{}
At the level of objects, the assertion of \cref{thm: Serre duality}
furnishes a canonical equivalence,
\[\bbD_{/X}: \IndCoh(Y/X)^{\vee} \xrightarrow{\simeq} \IndCoh(Y/X),\]
of module categories over $\QCoh(Y)$.
This functor also goes by the names of Grothendieck or
coherent duality in the literature (in particular in \cite{PreygelThesis}).

The main feature of $\bbD$ is that it is an anti-involution;
that is, there is a canonical natural equivalence,
\[\bbD_{/X}^2 \simeq \id_{\IndCoh(Y)}.\]
It follows from the construction that,
on the level of objects,
\[\bbD_{/X}(\F) = \sHom_Y(\F, \omega_Y),\]
where $\F \in \IndCoh(Y/X)$.
In particular, on the level of underlying stable $\infty$-categories,
the dual of $\IndCoh(Y/X)$ is its opposite category.

The Serre duality functor we have defined here
is compatible with that defined in \cite[\S 9]{GaitsgoryIndCoh}
in the sense that the following diagram commutes,
\beqn
\label{eqn: duality diagram}
\begin{tikzcd}
	{\IndCoh(Y/X)^{\vee}} & {\IndCoh(Y/X)} \\
	{\IndCoh(Y)^{\vee}} & {\IndCoh(Y)}
	\arrow["\oblv"', from=1-1, to=2-1]
	\arrow["\oblv", from=1-2, to=2-2]
	\arrow["\bbD", from=2-1, to=2-2]
	\arrow["{\bbD_{/X}}", from=1-1, to=1-2].
\end{tikzcd}
\eeqn

\begin{warn}
We often omit the subscript ``$/X$" on the relative Serre duality
functor, but it should be clear from context whether we are referring
to the usual $k$-linear duality functor $\bbD$ or the $\QCoh(X)$-linear duality
functor $\bbD_{/X}$.
\end{warn}


\section{Very good stacks}
\label{sec: very good stacks}
The notion of very good stack was introduced in \cite{PreygelT}
to refer to those stacks over which an integral kernel theorem
for $\IndCoh$ holds, and over which the relative box product is
an isomorphism. We recall the precise notion below, as well as
show that $B\G_m$ is a very good stack.

\subsection{Notation}
We introduce the following conditions on
a derived stack $X$:
\beqn
\label{stack condition 0} \tag{$\star$}
X \, \textrm{is Noetherian, has affine diagonal, and is perfect.}
\eeqn
\beqn
\label{stack condition} \tag{$\star_F$}
X \, \textrm{is Noetherian, has finite diagonal, is perfect, and is Deligne-Mumford.}
\eeqn

We say that a morphism of derived Noetherian
stacks $X \to Y$ is \labelcref{stack condition 0} (resp. \labelcref{stack condition}) 
if $X \times_Y \Spec A$ is a \labelcref{stack condition 0} (resp. \labelcref{stack condition}) 
stack for any $\Spec A \to Y$ almost of finite presentation.

A \labelcref{stack condition 0} (resp. \labelcref{stack condition}) 
derived $S$-stack is a stack over a Noetherian
derived stack $S$ such that the structure morphism is 
a \labelcref{stack condition 0} (resp. \labelcref{stack condition}) morphism of derived stacks.

\subsection{The definition}
In what follows, we abuse notation by using
$\IndCoh(-)$ to denote both the plain dg-category
and its $\QCoh(S)$-module enhancement $\IndCoh(-/S)$
(see \cref{warn: IndCoh abuse of notation}).

\bdef[{\cite[\S 4.1.3.1]{PreygelThesis}}]
\label{def: very good stack}
Let $S$ be a derived stack.
We say that $S$ is a \emph{very good stack}
is for any $X$ and $Y$ almost finitely-presented
\labelcref{stack condition} derived stacks over $S = \Spec k$, and
$Z_{X} \subset X$ and $Z_{Y} \subset Y$ closed subsets,
we have the following commutative diagram of equivalences
of $\QCoh(S)$-linear dg-categories:
\[\begin{tikzcd}
	{\Fun^L_{\QCoh(S)}(\IndCoh(X)_{Z_{X}},\IndCoh(Y)_{Z_{Y}})} & {\IndCoh(X \times_S Y)_{Z_{X} \times_S Z_{Y}}} \\
	{\IndCoh(X)_{Z_{X}} \otimes_{\QCoh(S)} \IndCoh(Y)_{Z_{Y}}}
	\arrow["\simeq", "\Phi"', from=1-2, to=1-1]
	\arrow["\Psi", "\simeq"', from=2-1, to=1-1]
	\arrow["\simeq", "\boxtimes"', from=2-1, to=1-2],
\end{tikzcd}\]
where
\begin{itemize}
\item $\boxtimes_S$ denotes the external tensor product over $S$, and restricts to
an equivalence on compact objects,
\[\boxtimes_S: \Coh_{Z_X}(X) \otimes_S \Coh_{Z_Y}(Y) 
	\xrightarrow{\simeq} \Coh_{Z_X \times_S Z_Y}(X \times_S Y);\]
\item $\Phi(\K) := {\pr_2}_*(\pr_1^!(-) \overset{!}{\otimes} \K)$ is 
the integral transform with kernel $\K$;
\item and $\Psi$ is uniquely defined by its values on compact objects
$\F$ and $\calG$ by $\Psi(\F \otimes \calG) := \Hom_S(\bbD(\F), -) \otimes_{\QCoh(S)} \calG$. 
\end{itemize}
\edefn

The prototypical example of a very good stack is $S = \Spec k$,
which satisfies the properties of a very good stack by
\cite[Theorem A.2.2.4]{PreygelThesis}.

\subsection{$B\G_m$ is a very good stack}
It may be shown that a certain class of stacks
that include $B\G_m$ are very good stacks
by imitating Preygel's proofs of \cite[Theorem A.2.2.4]{PreygelThesis}.
In order to do so, we introduce the following few propositions.

The following proposition is found in \cite{PreygelThesis}.

\bprop[{\cite[Proposition A.2.1.1]{PreygelThesis}}]
\label{prop: primordial shriek}
Let $S$ be a regular \labelcref{stack condition 0} stack,
and let $X$ and $Y$ be \labelcref{stack condition 0} $S$-stacks.
Then the exterior product determines a well-defined and fully faithful
functor,
\[\boxtimes_S: \Coh(X) \otimes_S \Coh(Y) \to \Coh(X \times_S Y).\]
\eprop

From the above proposition one may deduce the following
lemmas, which are stated in \cite{PreygelThesis} in less generality,
though the proofs contained therein apply equally well to the
more general statements written below.

\blem[{\cite[Lemma A.2.2.1]{PreygelThesis}}]
\label{lem: shriek1}
Suppose $X$ is a \labelcref{stack condition} stack
over a regular \labelcref{stack condition 0} stack $S$.
Then for any $\F, \calG \in \Coh(X)$, there are natural
equivalences,
\begin{enumerate}[label=(\roman*), ref=\roman*]
\item $\bbD(\F \boxtimes_S \calG) \simeq \bbD(\F) \boxtimes_S \bbD(\calG)$
\item $\omega_X \overset{!}{\otimes} \calG \simeq \calG$.
\end{enumerate}
\elem

\bproof
Let $\pr_i: X \times_S X \to X$ denote the standard $i$th projection.
Note that $\pr_2^!\calG \simeq \omega_X \boxtimes_S \calG$. In particular, $\omega_{X \times_S X} 
\simeq \pr_2^!\omega_X \simeq \omega_X \boxtimes_S \omega_X$. Part (i) now
follows from \cref{prop: primordial shriek}
and the formula $\bbD(-) = \sHom_X(-, \omega_X)$. Part (ii) follows by
noting that
\[\omega_X \overset{!}{\otimes} \calG = \Delta^!(\omega_X \boxtimes_S \calG) = \Delta^!\pr_2^!\calG = \calG.\]
\eproof

\blem[{\cite[Lemma A.2.2.2]{PreygelThesis}}]
\label{lem: shriek2}
Suppose $X$ is a \labelcref{stack condition} stack
over a \labelcref{stack condition 0} stack $S$.
Then for any $\F, \calG \in \Coh(X)$, there is a natural
equivalence in $\QCoh(X)$,
\[\F \overset{!}{\otimes} \calG \simeq \sHom_X(\bbD\F, \calG).\]
\elem

\bproof
As above, let $\pr_i: X \times_S X \to X$ be the standard $i$th projection,
and let $\Delta: X \to X \times_S X$ be the diagonal morphism relative to $S$.

Since $\Delta$ is finite it is proper, so the functor $\Delta_*^{\IndCoh}$ preserves compact
objects by \cite[Corollary 3.3.6]{GaitsgoryIndCoh}. This implies that $\Delta^!$ is continuous,
so the adjunction $(\Delta_*, \Delta^!)$ on $\IndCoh$ restricts
to an adjunction on $\Coh$. Since the functor $\Coh(X) \to \QCoh(X)$
is fully faithful, we therefore may write
\begin{align*}
\F \overset{!}{\otimes} \calG &\simeq \Delta^!\left(\F \boxtimes_S \calG\right) \\
&\simeq {\pr_1}_*\Delta_*\sHom_{\QCoh(X)}\left(\O_X, \Delta^!\left(\F \boxtimes_S \calG\right)\right) \\
&\simeq {\pr_1}_*\sHom_{\QCoh(X \times_S X)}\left(\Delta_*\O_X, \F \boxtimes_S \calG\right).
\end{align*}
Since $\Delta_*\O_X$ belongs to $\Coh(X \times_S X)$,
we may apply Serre duality to rewrite this as
\begin{align*}
&\simeq {\pr_1}_*\sHom_{\QCoh(X \times_S X)}\left(\Delta_*\O_X, \bbD\bbD\left(\F \boxtimes_S \calG\right)\right) \\
& \simeq {\pr_1}_*\sHom_{\QCoh(X \times_S X)}\left(\Delta_*\O_X, \sHom_{X \times_S X}\left(\bbD\left(\F \boxtimes_S \calG\right), \omega_{X \times_S X}\right)\right) \\
& \simeq {\pr_1}_*\sHom_{\QCoh(X \times_S X)}\left(\Delta_*\O_X \otimes \bbD\left(\F \boxtimes_S \calG\right), \omega_{X \times_S X} \right) \\
&\simeq {\pr_1}_*\sHom_{X \times_S X}\left(\bbD(\F \boxtimes_S \calG), \bbD(\Delta_*\O_X)\right).
\end{align*}
Applying \cref{lem: shriek1}(i), we obtain
\begin{align*}
&\simeq {\pr_1}_*\sHom_{X \times_S X}\left(\bbD\F \boxtimes_S \bbD(\calG), \bbD\left(\Delta_*\O_X\right)\right) \\
&\simeq {\pr_1}_*\sHom_{X \times_S X}\left(\bbD\F \boxtimes_S \O_X, \sHom_{X \times_S X}\left(\O_X \boxtimes_S \bbD\calG, \bbD(\Delta_*\O_X)\right)\right).
\end{align*}
Undoing the above operations on the inner $\sHom$ in the previous line,
we obtain
\begin{align*}
&\simeq {\pr_1}_*\sHom_{X \times_S X}\left(\bbD\F \boxtimes_S \O_X, \sHom_{X \times_S X}\left(\Delta_*\O_X, \bbD\left(\O_X \boxtimes_S \bbD\calG\right)\right)\right) \\
&\simeq {\pr_!}_*\sHom_{X \times_S X}\left(\bbD\F \boxtimes_S \O_X, \Delta_*\left(\omega_X \overset{!}{\otimes} \calG\right)\right).
\end{align*}
Applying the relative $(\Delta^*, \Delta_*)$ adjunction, we obtain
\begin{align*}
&\simeq {\pr_1}_*\Delta_*\sHom_X\left(\bbD\F \otimes \O_X, \omega_X \overset{!}{\otimes} \calG\right).
\end{align*}
Finally we conclude by using \cref{lem: shriek1}(ii) 
to write the above expression as
\begin{align*}
&\simeq \sHom_X(\bbD\F, \calG).
\end{align*}
\eproof

We now show that any smooth stack with smooth affine diagonal is very good.

\bprop
\label{prop: integral kernel theorem}
Let $S$ be a smooth stack with smooth affine diagonal.
Suppose that $f: X \to S$ and $g: Y \to S$ are almost finitely-presented
\labelcref{stack condition} derived stacks over $S$, and let
$Z_{X} \subset X$ and $Z_{Y} \subset Y$ be closed subsets.
Then the following diagram commutes,
\[\begin{tikzcd}
	{\Fun^L_{\QCoh(S)}(\IndCoh(X)_{Z_{X}},\IndCoh(Y)_{Z_{Y}})} & {\IndCoh(X \times_S Y)_{Z_{X} \times_S Z_{Y}}} \\
	{\IndCoh(X)_{Z_{X}} \otimes_{\QCoh(S)} \IndCoh(Y)_{Z_{Y}}}
	\arrow["\Phi"', from=1-2, to=1-1]
	\arrow["\Psi", from=2-1, to=1-1]
	\arrow["\boxtimes_S"', from=2-1, to=1-2].
\end{tikzcd}\]
where
\begin{itemize}
\item $\Phi(\K) := {\pr_2}_*(\pr_1^!(-) \overset{!}{\otimes} \K)$ is 
the integral transform with kernel $\K$,
\item and $\Psi$ is uniquely defined by its values on compact objects
$\F$ and $\calG$ by $\Psi(\F \otimes \calG) := \Hom_S(\bbD(\F), -) \otimes_S \calG$. 
\end{itemize}
Furthermore, all of these functors are equivalences of
$\QCoh(S)$-linear categories, and the equivalence induced
by the external tensor product restricts to an equivalence,
\[\boxtimes_S: \Coh_{Z_X}(X) \otimes_S \Coh_{Z_Y}(Y) 
	\xrightarrow{\simeq} \Coh_{Z_X \times_S Z_Y}(X \times_S Y).\]
\eprop

\bproof
Because the $\QCoh(S)$-linear Serre duality functor is compatible with
the $k$-linear one (see the commutative diagram 
\labelcref{eqn: duality diagram}), it respects supports
since the $k$-linear one does, so it restricts to the $\Perf(S)$-linear equivalence
\[\bbD(-) : \Coh_{Z_X}(X)^{\op} \simeq \Coh_{Z_X}(X).\] 
This implies that $\IndCoh_{Z_X}(X)$
is self-dual over $\QCoh(S)$ via $\bbD(-)$, so that $\Psi$ is an equivalence.
We will now verify that the diagram is commutative and then
that the relevant maps are equivalences.

\hfill

\noindent \textbf{Diagram commutes}: \\
Let us prove that the diagram commutes up to natural equivalence. Since each of $\boxtimes_S$, $\Psi$,
and $\Phi$ preserves colimits, it suffices to give a natural equivalence 
$\Psi_{\F \otimes \calG} = \Phi_{\F \boxtimes_S \calG}$
for $\F \in \Coh(X)$, $\calG \in \Coh(Y)$. Since both functors preserve colimits, we may check
this for $T \in \Coh(X)$:
\begin{align*}
\Phi_{\F \boxtimes_S \calG}(T) 	&:= {\pr_2}_*({\pr_1}^! T \overset{!}{\otimes} (\F \boxtimes_S \calG)) \\
							&\simeq {\pr_2}_*((T \boxtimes_S \omega_Y) \overset{!}{\otimes} (\F \boxtimes_S \calG)) \\
							&\simeq {\pr_2}_*\sHom_{X \times_S Y}(\bbD \F \boxtimes_S \bbD \calG, T \boxtimes_S \omega_Y) \\
							&\simeq {\pr_2}_*\left(\sHom_X(\bbD\F, T) \boxtimes_S \sHom_Y(\bbD\calG, \omega_Y)\right) \\
							&\simeq {\pr_2}_*\left(\sHom_X(\bbD\F, T) \boxtimes_S \sHom_Y(\O_Y, \calG)\right) \\
							&\simeq {\pr_2}_*\pr_1^!\sHom_X(\bbD\F, T) \overset{!}{\otimes} \calG \\
							&\simeq g^!f_*\sHom_X(\bbD\F, T) \overset{!}{\otimes} \calG \\
							&=: f_*\sHom_X(\bbD\F,T) \otimes_S \calG \\
							&\simeq \Hom_S(\bbD\F, T) \otimes_S \calG \\
							&=: \Psi_{\F \otimes_S \calG}(T).
\end{align*}
Here we have implicitly used \cref{lem: shriek1}, \cref{lem: shriek2}, and coherent duality.

\hfill

\noindent \textbf{Equivalences}: \\
Since the diagram commutes and $\Psi$ is an equivalence, it suffices to show that $\boxtimes_S$ is an
equivalence. By \cref{prop: primordial shriek} it preserves compact objects and is fully faithful. It suffices
to show that it is essentially surjective on compact objects. By \cite[Lemma 5.4.5]{PreygelT}, 
the assertion without support conditions is true, since $S$ is a smooth
stack with smooth affine diagonal over the very good stack $\Spec k$. 
Let us show how this implies the general case:
\[\begin{tikzcd}
	{\Coh(Z_X) \otimes_S \Coh(Z_Y)} & {\Coh(Z_X \times_S Z_Y)} \\
	{\Coh_{Z_X}(X) \otimes_S \Coh_{Z_Y}(Y)} & {\Coh_{Z_X \times_S Z_Y}(X \times_S Y)}
	\arrow["\simeq", from=1-1, to=1-2]
	\arrow[from=1-2, to=2-2]
	\arrow[from=1-1, to=2-1]
	\arrow[from=2-1, to=2-2]
\end{tikzcd}\]
We have seen that the bottom horizontal arrow is fully faithful, so since both categories
are stable and idempotent complete it suffices to show that it has dense image. We have
seen that the the top horizontal arrow is an equivalence. The right vertical arrow has dense
image by Lemma 2.2.0.2. Consequently, the bottom horizontal arrow has dense image as
desired.
\eproof

\bcor
\label{cor: BG_m is a very good stack}
$B\G_m$ is a very good stack.
\ecor

\bproof
This follows from \cref{prop: integral kernel theorem}
since $B\G_m$ is a smooth stack with smooth affine diagonal.
\eproof


\section{Singular support}
\label{sec: singular support}
Let $S$ be a Noetherian derived scheme.
By \cite[Proposition 1.6.4]{GaitsgoryIndCoh},
the functor $\Psi_S: \IndCoh(S) \to \QCoh(S)$ is an equivalence
if and only if $S$ is a regular classical scheme,
so the discrepancy between ind-coherent sheaves
and quasicoherent sheaves can be seen as a measure
of the failure of $S$ to be regular.

The theory of singular support for ind-coherent sheaves
developed by Arinkin
and Gaitsgory in \cite{AG15} is a refinement 
of the support of an ind-coherent sheaf
(see \cref{def: support of ind-coherent sheaf})
which quantifies this idea
for a special class of derived schemes: those
which are \emph{quasi-smooth}. The definition of a quasi-smooth
scheme is standard; a derived $k$-scheme is called
quasi-smooth over $k$ if its cotangent complex
is perfect of Tor amplitude $[-1,0]$.
The authors of \cite{AG15} also develop a theory
of singular support for quasi-smooth stacks,
though we will not need it for our purposes.

\subsection{Support theories and $\bbE_2$-algebra actions}
\label{ssec: E2-algebra actions}

Recall that, by definition, an action of an $\bbE_2$-algebra 
$\calA$ on a dg-category $\scrC$ is an 
$\bbE_2$-algebra map, 
$\calA \to \HH^{\bullet}_k(\scrC)$,
where $\HH^{\bullet}_k(\scrC)$ denotes
the $\bbE_2$-algebra of $k$-linear Hochschild cochains
on $\scrC$.
Since $\calA$ has an $\bbE_2$-algebra
structure, the graded classical algebra,
 \[A := \bigoplus_i H^{2i}(\calA),\]
is commutative. From the map of
$\bbE_2$-algebras $\calA \to \HH^{\bullet}_k(\scrC)$
we obtain a graded map of classical commutative
algebras,
\[A \to \bigoplus_{i} \HH^{2i}_k(\scrC).\]
The target of this map is known as the
``graded center of $\scrC$." Given
a graded map into the graded center of $\scrC$---such
as that obtained from an $\bbE_2$-algebra 
action---one may construct a support theory for
objects of $\scrC$ in which the support of
an object is a conic closed subset of $\Spec A$.
The construction and details of this support theory may
be found in \cite[\S 3]{AG15}, but the reader
may treat this as a black box.

\subsection{The definition of singular support}
\label{ssec: the definition of singular support}
\subsubsection{The scheme of singularities}
Let $\Z$ be a quasi-smooth derived scheme,
and let $\bbL_{\Z/k}$ denote its cotangent complex.
Then there is an associated derived scheme,
\[\T^*[-1]\Z := \Spec \left(\Sym_{\O_Z}(\bbL_{\Z/k}^{\vee}[1])\right)\]
whose underlying classical scheme is called
``the scheme of singularities" of $\Z$ by the authors
of \cite{AG15}. 

\brem
The authors of \cite{AG15} denote
this classical scheme by $\Sing(\Z)$. In order to
avoid a clash in notation, we choose to denote it
instead by $T^*[-1]\Z$.
\erem

\subsubsection{Affine quasi-smooth schemes}
Let $\Z$ be an affine quasi-smooth derived scheme.

There is a canonical Lie algebra structure on
$\bbL_{\Z/k}^{\vee}[-1]$ as an object of $\QCoh(\Z)$,
and the underlying $\bbE_1$-algebra of $\HH^{\bullet}_k(\IndCoh(\Z))$
receives a canonical map from $\Gamma(\Z;U_{\O_{\Z}}(\bbL_{\Z/k}^{\vee}[-1]))$,
the global sections of the universal
enveloping algebra of $\bbL_{\Z/k}^{\vee}[-1]$. 
This map is an isomorphism
if $\Z$ is eventually coconnective (as it always is in our applications).
In particular, we obtain a canonical map of
commutative algebras,
\[\Gamma(Z; \O_Z) \to \HH^0_k(\IndCoh(\Z)),\]
hence a map of $\Gamma(Z; \O_Z)$-modules,
\[\Gamma(\Z; H^1(\bbL_{\Z/k}^{\vee})) \to \HH^2_k(\IndCoh(\Z)),\]
where $Z$ denotes the underlying classical
scheme of $\Z$.
It follows that there is a canonical \emph{graded} map
of classical commutative algebras,
\begin{align*}
\Gamma(T^*[-1]\Z; \O_{T^*[-1]\Z}) 		&\simeq \Gamma(Z; \Sym_{\O_Z}(H^1(\bbL_{\Z/k}^{\vee}))) \\ 
								&\to \bigoplus_i \HH^{2i}_k(\IndCoh(\Z)),
\end{align*}
where $\Gamma(Z; H^1(\bbL_{\Z/k}^{\vee}))$
is assigned degree $2$.
Thus we obtain a graded map from the algebra
of regular functions on $T^*[-1]\Z$ to the graded
center of $\IndCoh(\Z)$, and therefore a support
theory assigning to each $\F \in \IndCoh(\Z)$
a closed conical subset of $T^*[-1]\Z$ called
the \emph{singular support of $\F$}, which
we denote by $\singsupp(\F) \subset T^*[-1]\Z$.

\subsubsection{General quasi-smooth schemes}
The singular support of an object $\F \in \IndCoh(\Z)$ for
an arbitrary (i.e. not necessarily affine) quasi-smooth
derived scheme may be defined by choosing an affine
open cover $\{U_{\a}\}$ of $\Z$ and setting,
\[\singsupp(\F) \cap T^*[-1]U_{\a} := \singsupp(\F|_{U_{\a}}),\] 
where $\F|_{U_{\a}} = j_{U_{\a}}^{\IndCoh,*}\F$.
This determines a closed conic subset of $T^*[-1]\Z$, 
called the singular support of $\F$, which
is independent of the choice of affine cover by
\cite[Corollary 4.5.7]{AG15}.

\subsection{Singular support on global complete intersections}
The prototypical example of a quasi-smooth derived
scheme is a fiber of a map between smooth schemes.

\bdef
Let $f: X \to V$ be a map, where $X$ is a smooth
affine scheme and $V$ is a finite dimensional vector
spaces viewed as a scheme, and denote the derived
zero fiber of $f$ by $\Z(f)$. 

A derived scheme
$\Z$ is called a \emph{global complete intersection}
if $\Z \simeq \Z(f)$ for some such $f$.
\edefn

On a global complete intersection $\Z$, the singular support 
for ind-coherent sheaves admits a
description as the support for a certain module structure
on $\IndCoh(\Z)$, as is explained in detail in
\cite[\S 5]{AG15}. We recall the module structure
below for reference in the main text.

\begin{construction}
\label{construction: singular support QCoh action}
Let $\Z$ be a global complete intersection,
and let $f: X \to V$ be a map presenting it
as a derived zero fiber. 
Then $\IndCoh(\Z)$ obtains a natural 
$\QCoh(V^{\vee}/\G_m)$-module
structure as follows.

There is a group structure on the derived
scheme $\Omega_0 V$ that corresponds
heuristically to composition of based loops.
There is an obvious group action of $\Omega_0 V$
on $\Z$\footnotemark 
\footnotetext{This action is the same
one as described in \cite[Construction 3.1.5]{Preygel},
but with $V$ instead of $\Aone$ as the target of
the map $f$.} The action of $\Omega_0 V$ on $\Z$
induces an action of $\IndCoh(\Omega_0 V)$ on
$\IndCoh(\Z)$ under convolution.
Let $\calA_{\Omega_0 V}$ denote the
$\bbE_2$-algebra of endomorphisms
of the unit for the convolution monoidal
structure in $\IndCoh(\Z)$,
also known as the algebra of "Hochschild
cochains on $0$ relative to $V$."\footnotemark 
\footnotetext{See \cite[\S F.4]{AG15} for reference.}
By the functoriality of Hochschild cochains,
this action induces a map of $\bbE_2$-algebras,
\[\calA_{\Omega_0 V} \to \HH^{\bullet}_k(\IndCoh(\Z)),\]
hence an action of $\calA_{\Omega_0 V}^{\op}\mod$ on
$\IndCoh(\Z)$.

Note that the weight $d$ scaling $\G_m$-action on $V$
induces a corresponding $\G_m$-action on the
$\bbE_2$-algebra $\calA_{\Omega_0 V}$.\footnotemark
\footnotetext{i.e. a lift of $\calA_{\Omega_0 V}$ to an
$\bbE_2$-algebra object in $\G_m\mod$.}
For reasons that will become clear below,
we will consider $\calA_{\Omega_0 V}$
with the $\G_m$-action
coming from the the degree $2$ scaling action on $V$.

Let $\Sym(V[-2])$ be the graded commutative dg-algebra
whose $\G_m$-action is also induced by the degree
$2$ scaling action on $V$. There is actually a canonical
commutative (i.e. $\bbE_{\infty}$-algebra) structure on $\calA_{\Omega_0 V}$
that restricts to its initial $\bbE_2$-algebra structure (\cite[Lemma 5.4.2]{AG15}),
so $\calA_{\Omega_0 V}^{\op}$ is equivalent under
Koszul duality to $\Sym(V[-2])$ as graded commutative dg-algebras.
We therefore obtain a $\Sym(V[-2])\mod$-action
on $\IndCoh(\Z)$, hence a $(\Sym(V[-2])\mod)^{\G_m}$-action
on $\IndCoh(\Z)$. The invariance of $\G_m$-invariants
under shearing (see \cref{sec: shearing}
below) gives,
\begin{align*}
(\Sym(V[-2])\mod)^{\G_m} &\simeq (\Sym(V[-2]^{\shear}\mod)^{\G_m} \\
					&\simeq (\Sym(V)\mod)^{\G_m} \\
					&\simeq \QCoh(V^{\vee}/\G_m).
\end{align*}
Hence we obtain a $\QCoh(V^{\vee}/\G_m)$-module structure
on $\IndCoh(\Z)$ as promised.
\end{construction}

The $\QCoh(V^{\vee}/\G_m)$-action on $\IndCoh(\Z)$
constructed above clearly extends to a $\QCoh(X) \otimes \QCoh(V^{\vee}/\G_m)
\simeq \QCoh(X \times V^{\vee}/\G_m)$-action,
where the $\QCoh(X)$-action is given by composing the
canonical $\QCoh(\Z)$-action on $\IndCoh(\Z)$ with
the pullback along $\Z \hook X$.
It is shown in \cite[\S 5.4]{AG15} that the 
support theory obtained from this 
$\QCoh(X \times V^{\vee}/\G_m)$-action
on $\IndCoh(\Z)$ is the singular support theory of
coherent sheaves under the image of the
inclusion $T^*[-1]\Z/\G_m \hook X \times V^{\vee}/\G_m$. 

\brem
The choice to consider the degree $2$ scaling action
on $V$ in \cref{construction: singular support QCoh action}
so that the $\QCoh(V^{\vee}/\G_m$)-action obtained
in that construction was compatible with the
degree assigned to $\Gamma(Z; H^1(\bbL_{\Z/k}^{\vee}))$---and
therefore the scaling action on the fibers of $T^*[-1]\Z$---in the
definition of singular support in \S\labelcref{ssec: the definition
of singular support}.
\erem

\section{Formal groups}
\label{sec: formal groups}
One of the key observations of \cite{PreygelThesis}, which the
author of \textit{op. cit.} attributes to Constantin
Teleman, is that the data of a Landau--Ginzburg pair
$(M,f)$ is equivalent to the data of a certain formal group action on the
category of coherent sheaves on $M$. This formal group action
is the linchpin of the computations in 
\cref{sec: comparing Hochschild invariants} 
comparing the Hochschild
invariants of graded and ungraded matrix factorizations. 
In this appendix, we collect the
facts about formal groups and their actions
relevant to these computations.

\subsection{Formal moduli problems}
\subsubsection{}
Let $\calX \in \PreStk_{\laft}$.
We recall the definition of \emph{formal moduli
problem over $\calX$} given in \cite{GR17II}.

\bdef
Let $\Moduli_{/\calX}$ denote the full subcategory $(\PreStk_{\laft})_{/\calX}$
spanned by inf-schematic nil-isomorphisms $\calY \to \calX$.
Objects of $\Moduli_{/\calX}$ are called formal moduli
problems over $\calX$.
\edefn

\bdef
A \emph{formal group over $\calX$} is a group
object in $\Moduli_{/\calX}$.
\edefn

When $\calX = B\G_m$, we use the terms
``graded formal moduli problem" and ``graded
formal group" instead.

\subsubsection{}
A central principle in deformation theory is that
the formal neighborhood of any point in a moduli
space is controlled by a differential graded Lie algebra.
The relationship between Lie algebras and formal groups
in general is given by the following theorem.

\bthm[{\cite[Theorem 3.1.4]{GR17II}}]
\label{thm: formal group = Lie algebras}
There are adjoint equivalences,
\[\begin{tikzcd}
	{\exp_{\calX}: \LieAlg(\IndCoh(\calX))} && {\Grp(\Moduli_{/\calX}): \Lie_{\calX}}
	\arrow[shift left=1, from=1-1, to=1-3]
	\arrow[shift left=1, from=1-3, to=1-1],
\end{tikzcd}\]
where $\LieAlg(\IndCoh(\calX))$ denotes the category
of Lie algebra objects in $\IndCoh(\calX)$.
\ethm

\subsection{Formal group actions on categories}
\label{ssec: Formal group actions on categories}
Given a formal group $\calG$ over a smooth stack
$\calX$, we would like to formulate a notion of 
$\calG$-action on a small, idempotent complete,
$\Perf \calX$-linear category $\scrC$.

When $\calX = \Spec k$, the theory of such actions
is treated in detail in \cite[\S 5.3]{PreygelThesis}. We
will not be as rigorous in our treatment
of formal groups actions over more general
bases as they are in \textit{loc. cit.},
but the difficulties explained therein
are precisely the same difficulties we
would encounter in trying to make our
presentation more thorough.

Instead, suffice it to make the following
definition, using \cref{thm: formal group = Lie algebras}
for justification.

\bdef
\label{def: modules over a formal group}
Let $\calG$ be a derived formal group over the
smooth stack $\calX$,
and let $L_{\calG}$ denote $\Lie_{\calX}(\calG)$.
Then a $\calG$-category is the data of
of a small $\Perf \calX$-linear category $\scrC$
and a map of Lie algebras,
\[L_{\calG} \to \HH^{\bullet}_{\calX}(\scrC)[+1],\]
where we have used the known fact that
$\HH^{\bullet}_{\calX}(\scrC)[+1] \in \LieAlg(\QCoh(\calX))$.
\edefn

\begin{notn}
We denote the category of all $\calG$-categories
by $\calG\bmod$.
\end{notn}

\subsubsection{}
A formal group action on $\scrC$ may also be
written more directly in terms of formal moduli
problems as follows.

Let ${^{<\infty}}\Sch^{\ft, \aff}_{/\calX}$ denote the
category of affine schematic, finite type prestacks
over $\calX$.

Unwinding the definitions, it is easy to
see that a formal moduli problem over $\calX$ is encoded
by a functor,
\[({^{<\infty}}\Sch^{\ft, \aff}_{/\calX})^{\op} \to \Grpd_{\infty}.\]
Moreover, every such functor admits
a completion to a formal moduli problem.

In fact, by a variant of \cite[Proposition 1.2.2]{GR17II},
a formal moduli problem over $\calX$ is determined by its restriction
to the full subcategory of ${^{<\infty}}\Sch^{\ft, \aff}_{/\calX}$
spanned by the nil-isomorphisms $S \to \calX$. We denote this category
by ${^{<\infty}}\Sch^{\ft, \aff}_{\sim^{\nil}/\calX}$.

\blem
\label{lem: category of test objects is pointed}
Suppose that $\calX$ is a reduced prestack.
Then every $S \to \calX \in {^{<\infty}}\Sch^{\ft, \aff}_{\sim^{\nil}/\calX}$
receives a unique map from $\calX$ such that its composition with $S \to \calX$
is the identity on $\calX$.
\elem  

\bproof
Because the map $S \to \calX$ is a nil-isomorphism, it induces
an isomorphism on reduced prestacks, $S^{\red} \xrightarrow{\simeq} \calX^{\red}$.
On the other hand, there is a canonical map $S^{\red} \to S$ of prestacks
over $\calX$, hence a map $\calX^{\red} \to S$. 
Since $\calX$ is reduced, the canonical map
$\calX^{\red} \xrightarrow{\simeq} \calX$ is an equivalence,
and we therefore obtain a map $\calX \to S$. Now, observe
that $\calX \to S$ is a map of prestacks over $\calX$ by construction,
so composes with $S \to \calX$ to be the identity.
\eproof

As a result of \cref{lem: category of test objects
is pointed}, the category ${^{<\infty}}\Sch^{\ft, \aff}_{\sim^{\nil}/\calX}$
is pointed, and pullback along the unique map $\calX \to S$ furnished by the lemma
gives a symmetric monoidal functor $\Perf S \to \Perf \calX$ 
for every $S \in {^{<\infty}}\Sch^{\ft, \aff}_{\sim^{\nil}/\calX}$.
Using this fact, we formulate the following definition.

\bdef
Given a small $\Perf \calX$-linear
category $\scrC$, we define the formal
moduli problem $(\dgcat_{\Perf \calX}^{\idm})^{\wedge}_{\scrC}$
to the be the formal moduli completion of the
functor,
\[(S \to \calX) \mapsto (\wit{\scrC} \in \dgcat_{\Perf \calY}^{\idm}, 
	\wit{\scrC} \otimes_{\Perf S} \Perf \calX \simeq \scrC),\]
sending $S \in {^{<\infty}}\Sch^{\ft, \aff}_{/\calX}$ to
the $\infty$-groupoid of pairs of a $\Perf S$-linear
category $\wit{\scrC}$ and an equivalence $\wit{\scrC} \otimes_{\Perf S} \Perf \calX \simeq \scrC$.
\edefn

With this definition in hand, one can show that the space of
maps $L_{\calG} \to \HH^{\bullet}_{\calX}(\scrC)$ is
equivalent to the space of maps of formal moduli problems,
\[B\calG \to (\dgcat_{\Perf \calX}^{\idm})^{\wedge}_{\scrC},\]
where $B$ here denotes the delooping functor taken in
the category of formal moduli problems.

\subsubsection{}
Just as in the case of a group stack acting on
a category, there is a functor of $\calG$-invariants,
\[(-)^{\calG}: \calG\bmod \to \Perf \O(\calX)^{\calG}\bmod \simeq \Perf B\calG\bmod,\]
where $\O({\calX})$ denotes the global sections
of the structure sheaf of $\calX$, and $\O(\calX)^{\calG}$ denotes the invariants
for the trivial $\O({\calX})$-linear $\calG$-representation $\O({\calX})$.

\brem
Note that $B\calG$ is obtained by applying the delooping functor \emph{internal}
to the category of formal groups over $\calX$. 
\erem

In the case when $\calX = \Spec k$, 
we note that the algebra $k^{\calG}$ is computed
by the Chevalley--Eilenberg  cochain complex $C^*(\frakg)$ of the
Lie algebra object $\frakg$ associated to $\calG$.

\subsection{Formal groups over $\Spec k$ and $B\G_m$}
In this work, we are only interested in the
cases when $\calX$ is either the point
$\Spec k$, or $B\G_m$. The theory in the
former case is thoroughly developed in \cite[\S 5.3]{PreygelThesis},
as we remarked earlier. In particular, for our purposes
it is important that Preygel shows in \textit{loc. cit.} 
that the Hochschild 
invariants of a $\calG$-module category are
well-defined.

\subsubsection{}
Let $\bbE_2^{\mathrm{calc}}$ denote the $2$-colored operad
governing pairs of an $\bbE_2$-algebra and a module
over it with a $S^1$-action compatible with the circle action
on the $\bbE_2$-operad. Then we have the following proposition.

\bprop[{\cite[Proposition 5.3.3.11]{PreygelThesis}}]
\label{prop: Toly Prop 5.3.3.11}
Let $\calG$ be a formal group over $k$.
An $\calG$-action on a category $\scrC_0$ induces an 
$\calG$-action on its Hochschild invariants. That is,
one has a functor of $\infty$-categories
\[(\HH_{\bullet}, \HH^{\bullet}): \calG\bmod((\Vect_k)^{\sim}) \to \bbE_2^{\mathrm{calc}}\aalg(\calG\mod),\]
where $(\Vect_k)^{\sim}$ denotes the
$\infty$-groupoid given by
discarding non-invertible morphisms of $\Vect_k$.
\eprop

\Cref{prop: Toly Prop 5.3.3.11} is formulated using $(\Vect_k)^{\sim}$
is because Hochschild cohomology is not functorial with respect to
arbitrary maps of dg-categories. Hochschild homology, however, is
functorial with respect to continuous functors of dg-categories,
so we have the following corollary.

\bcor
Let $\calG$ be a formal group over $k$.
An $\calG$-action on a category $\scrC_0$ induces an 
$\calG$-action on its Hochschild homology.
That is, one has a functor,
\[\HH_{\bullet}: \calG\bmod \to S^1\mod(\calG\mod),\]
where $S^1\mod(\calG\mod)$ denotes the category
of $\calG$-modules equipped with an $S^1$-action.
\ecor

\subsubsection{}
One may prove an analogue of \cref{prop: Toly Prop
5.3.3.11} for Hochschild homology of graded
$\calG$-module categories.

\bprop
\label{prop: graded version of Toly Prop 5.3.3.11}
Let $\calG$ be a graded formal group (i.e. a formal
group over $B\G_m$).
A $\calG$-action on a graded category $\scrC_0$ induces an 
$\calG$-action on its Hochschild invariants. That is,
one has a functor of $\infty$-categories
\[(\HH_{\bullet}, \HH^{\bullet}): \calG\bmod((\Rep(\G_m)\mod)^{\sim}) \to \bbE_2^{\mathrm{calc}}\aalg(\calG\mod_{\gr}).\]
\eprop


Similarly, if we consider only Hochschild homology,
we obtain the following corollary.

\bcor
Let $\calG$ be a formal group over $k$.
An $\calG$-action on a category $\scrC_0$ induces an 
$\calG$-action on its Hochschild homology.
That is, one has a functor,
\[\HH_{\bullet}: \calG\bmod(\Rep(\G_m)\mod) \to S^1\mod(\calG\mod_{\gr}),\]
where $S^1\mod(\calG\mod_{\gr})$ denotes the category
of graded $\calG$-modules equipped with a \emph{graded} $S^1$-action.
\ecor

\subsection{Formal groups coming from vector spaces}
The only formal groups we consider in the main
body of the text are those obtained as the delooping
of the formal group given by the completion 
of a vector space at the identity $0$. 
More precisely, let $V$ be a vector space
over $k$, viewed as an abelian group scheme.
Then we are concerned with
\begin{itemize}
\item the formal group $B\wih{V}$ over $k$, where
$\wih{V}$ denotes the formal completion of $V$ at $0 \in V$;
\item the graded formal group $B\wih{V}/\G_m$ where
$B\wih{V}$ is considered with the $\G_m$-action
induced by the weight $d$ scaling action on $V$
for some $d \in \bbN$.
\end{itemize}

\subsubsection{}
The following theorem of Preygel can as an
instance of Cartier duality; see \cite[\S 5.1]{PreygelThesis}.

\bprop[{\cite[Corollary 5.4.0.15]{PreygelThesis}}]
\label{prop: BVhat-actions}
Suppose $\scrC \in \dgcat_k^{\idm}$, and $V$ is a vector space over $k$. 
Let $V[+1]$ be the corresponding abelian 
dg-Lie algebra, and let $B\wih{V}$ 
be the corresponding formal group. Then the following spaces are naturally equivalent
\begin{enumerate}[label=(\alph*), ref=\alph*]
\item $\Map_{\LieAlg}(V[+1], \HH^{\bullet}_k(\scrC)[+1])$;
\item $\Map_{\bbE_2\aalg}(\Sym_k V, \HH^{\bullet}_k(\scrC))$;
\item \{$V[+1]$-actions on $\scrC$\};
\item \{$B\wih{V}$-actions on $\scrC$\}.
\end{enumerate}
Furthermore, regard $\Perf k$ as a $B\wih{V}$-module
via the $\bbE_2$-algebra map $\Sym_k V \to k = \HH^{\bullet}(\Perf k)$
obtained by universal property from the zero map $V \to k$. 
Let $V^{\vee} = \Spec \Sym_k V$ 
as a commutative group scheme. Then 
\begin{itemize}
\item[$-$] there are equivalences
\[\scrC^{B\wih{V}} = \scrC \otimes_{\Sym_k V} k 
	\hspace{10mm} \scrC_{B\wih{V}} = \Fun^{\ex}_{\Sym_k V}(\Perf k, \scrC)\]
of module categories over the symmetric monoidal category $(\Perf k)^{B\wih{V}} = 
\Fun^{\ex}_{\Sym_k V}(\Perf k, \Perf k)$;
\item[$-$] the symmetric monoidal category $(\Perf k)^{B\wih{V}}$ 
can be identified with the convolution category
$(\Coh \Omega_0 V^{\vee}, \circ)$;
\item[$-$] if $V = k$, so that $V^{\vee} = \G_a$, 
then there is a symmetric monoidal
equivalence $(\Coh \Omega_0 \G_a, \circ) \simeq (\Perf k[\b], \otimes)$.
\end{itemize}
\eprop

\subsubsection{}
Predictably, there is an analogue of \cref{prop: BVhat-actions} for
formal groups of the form $B\wih{V}/\G_m$, obtained
by performing the same analysis in the graded context.

\bprop
\label{prop: BVhat/G_m-actions}
Suppose $\scrC \in \Perf(B\G_m)\mod(\dgcat_k^{\idm})$, and 
$V$ is a vector space over $k$, graded by the natural weight $d$ 
$\G_m$-scaling action, for fixed $d \in \bbN$. Let $V[+1]_{\gr}$ be the
Lie algebra in $\Rep(\G_m)$, and let $B\wih{V}/\G_m$ 
be the corresponding graded formal group. Then the following spaces 
are naturally equivalent
\begin{enumerate}[label=(\alph*), ref=\alph*]
\item $\Map_{\LieAlg(\Rep \G_m)}(V[+1]_{\gr}, \HH^{\bullet}_{\Rep \G_m}(\scrC)[+1])$;
\item $\Map_{\bbE_2\aalg(\Rep \G_m)}(\Sym_{\O(B\G_m)} V, \HH^{\bullet}_{\Rep \G_m}(\scrC))$;
\item \{$V[+1]_{\gr}$-actions on $\scrC$\};
\item \{$B\wih{V}/\G_m$-actions on $\scrC$\},
\end{enumerate}
where $V[+1]_{\gr}$ denotes the graded Lie
algebra determined by the $\G_m$-representation
$V$, and $\Sym_{\O(B\G_m)}(-)$ denotes the left
adjoint to the forgetful functor $\bbE_2\aalg(\Rep(\G_m))
\to \Rep(\G_m)$.
 
Furthermore, regard $\Perf B\G_m$ as a $B\wih{V}/\G_m$-module
via the $\bbE_2$-algebra map $\Sym_{\O(B\G_m)} V \to \O(B\G_m) = \HH^{\bullet}_{\Rep(\G_m)}(\Perf B\G_m)$
obtained by universal property from the zero map of graded
$k$-algebras $V \to k$. 
Let $V^{\vee}/\G_m = (\Spec \Sym_k V)/\G_m$ 
as a graded commutative group scheme. Then 
\begin{itemize}
\item[$-$] there are equivalences
\[\scrC^{B\wih{V}/\G_m} = \scrC \otimes_{\Sym_{\O(B\G_m)}V} \O(B\G_m) 
	\hspace{10mm} \scrC_{B\wih{V}/\G_m} = \Fun^{\ex}_{\Sym_{\O(B\G_m)}V}(\Perf B\G_m, \scrC)\]
of module categories over the symmetric monoidal category $(\Perf B\G_m)^{B\wih{V}} = \\
\Fun^{\ex}_{\Sym_k V}(\Perf B\G_m, \Perf B\G_m)$;
\item[$-$] the symmetric monoidal category $(\Perf B\G_m)^{B\wih{V}}$ 
can be identified with the convolution category
$(\Coh \Omega_0 V^{\vee}/\G_m, \circ)$;
\item[$-$] if $V = k$, so that $V^{\vee} = \G_a$, 
then there is a symmetric monoidal
equivalence of graded categories 
$(\Coh \Omega_0 \G_a/\G_m, \circ) \simeq (\Perf k[\b]_{\gr}, \otimes)$,
where $\b$ has $\G_m$-weight $d$.
\end{itemize}
\eprop

\section{Shearing}
\label{sec: shearing}
This appendix is included as a reference
for the basic facts about shearing that we
use throughout the main body of the text.
Our source for all the material appearing in this section
is the upcoming work of Ben-Zvi, Sakellaridis,
and Venkatesh, \cite{BZSV}.

\subsubsection{}
\bdef
\label{def: shearing}
The $n$-shearing functor is an auto-equivalence
of the dg-category $\Rep(\G_m)$, denoted
$(-)^{n\shear}: \Rep(\G_m) \to \Rep(\G_m)$, given on objects by
\[M = \oplus_i M_i \mapsto M^{n\shear} := \oplus_i M_i[ni].\]
We use ``\emph{the} shearing functor" to refer
to this functor for $n=1$, and we denote it simply
by $(-)^{\shear}$.

We denote the inverse to the $n$-shearing
functor by $(-)^{n\,\unshear}$, and refer to it
as ``$n$-unshearing." The $n$-unshearing functor for $n=1$
is again called ``the unshearing functor" and 
denoted simply by $(-)^{\unshear}$.
\edefn

Critically, the shearing functor and its odd variations are 
not \emph{symmetric monoidal} equivalences of
$\Rep(\G_m)$, though its even variants are. In order
to obtain a shearing functor that \emph{is} symmetric
monoidal, we must define the shearing auto-equivalence
on the category of $\G_m$ super-representations.

\subsubsection{}
Recall that a super vector space is a $\bbZ/2$-graded
vector space $V = V_0 \oplus V_1$, or, equivalently, a vector space $V$ 
with an action of the square roots of unity, $\mu_2 \subset \G_m$.

\bdef
The \emph{parity} of an element $v \in V$ is the index
of the summand in the $\bbZ/2$-grading of $V$ in which
it lies. 
\begin{itemize}
\item We say $v$ has \emph{even} parity if $v \in V_0$ or,
equivalently, if it lies in the $1$-eigenspace 
for involution of $V$ given by the action of
$-1 \in \mu_2$. 
\item We say that $v$ has \emph{odd} parity
if $v \in V_1$ or, equivalently, if it lies in the $-1$-eigenspace
for the action of $-1 \in \mu_2$.
\end{itemize}
\edefn

\begin{notation}
\label{notn: parity}
If $v$ has even parity or odd parity, we may occasionally 
also say it has parity $1$ or $-1$, respectively.
\end{notation}

\bdef
Let $\Vect_k^{\super}$ denote the symmetric monoidal
dg-category of complexes of $k$ super-vector spaces.
For $n \in \bbZ$, let $k\langle n \rangle$ 
denote the object of $\Vect_k^{\super}$
whose underlying object of $\Vect_k$ is $k[n]$, and which has
parity $(-1)^n$.
\edefn

If $V \in \Vect_k^{\super}$, we denote by
$V\langle n \rangle$ the object $k\langle n \rangle
\otimes V \in \Vect_k^{\super}$. We call ``$\langle n \rangle$"
the parity-corrected shift of a super-complex. If $V$ is discrete
and even, then $V\langle n \rangle$ is in parity $(-1)^n$\footnotemark
\footnotetext{See \cref{notn: parity}.}
and cohomological degree $-n$.

\bdef
Let $\Rep^{\super}(\G_m)$ denote the symmetric monoidal
dg-category of $\G_m$ super-representations. 
Objects of $\Rep^{\super}(\G_m)$
are complexes of $k$ super-vector spaces equipped with
the structure of a comodule over the coalgebra $\O(\G_m) = k[x^{\pm 1}]$.
\edefn

\bdef
\label{def: super-shearing}
The (super) $n$-shearing functor
is a symmetric monoidal auto-equivalence of
$\Rep(\G_m)^{\super}$, denoted
$(-)^{n\shear}: \Rep^{\super}(\G_m) \to \Rep^{\super}(\G_m)$,
given on objects by,
\[M = \oplus_i M_i \mapsto M^{n\shear} := \oplus_i M_i\langle ni \rangle.\]

All the notational and naming conventions
of \cref{def: shearing} apply \textit{mutatis mutandis}
to the super setting.
\edefn

There is an obvious symmetric monoidal inclusion
$\Rep(\G_m) \subset \Rep^{\super}(\G_m)$
given by the full subcategory of even complexes.
The restriction of the super shearing functor
$(-)^{\shear}$ to this full subcategory induces 
a symmetric monoidal equivalence,
\beqn
\label{eqn: shearing equivalence}
(-)^{\shear}: \Rep(\G_m) \xrightarrow{\simeq} \Rep_{\epsilon}^{\super}(\G_m),
\eeqn
where $\Rep_{\epsilon}^{\super}(\G_m) \subset \Rep^{\super}(\G_m)$
is the full subcategory on objects in which
the $i$-isotypical space, on which $\lambda \in \G_m$ acts
by scalar multiplication by $\lambda^i$,
has parity $(-1)^i$ (i.e. the parity of an element is given by 
the action of $\epsilon = -1 \in \G_m$).

\subsubsection{}
The shearing operations on $\Rep(\G_m)$
and $\Rep^{\super}(\G_m)$ induce corresponding
operations on their categories of modules.

\bdef
\label{def: shearing of categories}
The $n$-shearing functor on graded categories
is an auto-equivalence of $\Rep(\G_m)\mod$,
also denoted $(-)^{n\shear}$, given by
twisting the $\Rep(\G_m)$ action by the $n$-shearing
auto-equivalence of $\Rep(\G_m)$.
Explicitly, it is given on objects by
\[\scrC \mapsto \scrC \otimes_{\Rep(\G_m)} \Rep(\G_m)^{n\shear},\]
where $\Rep(\G_m)^{n\shear}$ denote $\Rep(\G_m)$
a $\Rep(\G_m)$-bimodule, with left action of $Y$ by
given by $Y^{n\shear} \otimes -$ and right action given
by $- \otimes Y$.
\edefn

\bdef
\label{def: super shearing of categories}
The (super) $n$-shearing functor on categories
is a symmetric monoidal auto-equivalence of $\Rep(\G_m)\mod$,
also denoted $(-)^{n\shear}$, given by
twisting the $\Rep^{\super}(\G_m)$ action by the $n$-shearing
auto-equivalence of $\Rep^{\super}(\G_m)$.
Explicitly, it is given on objects by
\[\scrC \mapsto \scrC \otimes_{\Rep^{\super}(\G_m)} \Rep^{\super}(\G_m)^{n\shear},\]
where $\Rep^{\super}(\G_m)^{n\shear}$ denote $\Rep^{\super}(\G_m)$
a $\Rep^{\super}(\G_m)$-bimodule, with left action of $Y$ by
given by $Y^{n\shear} \otimes -$ and right action given
by $- \otimes Y$.
\edefn

\brem
The shearing functor on graded categories
leaves the underlying dg-category unchanged. More
precisely, shearing induces an equivalence
$\scrC \xrightarrow{\simeq} \scrC^{\shear}$.
\erem

It follows from the equivalence \labelcref{eqn: shearing equivalence}
the shearing endomorphism of $\Rep^{\super}(\G_m)\mod$ restricts
to symmetric monoidal equivalence,
\beqn
\label{eqn: shearing equivalence 2}
(-)^{\shear}: \Rep(\G_m)\mod \xrightarrow{\simeq} \Rep_{\epsilon}^{\super}(\G_m)\mod.
\eeqn
As such, given any algebra object $\sfA \in \Rep(\G_m)\mod$,
the shearing functor induces a symmetric monoidal equivalence 
between modules over $\sfA$ and modules 
over the sheared category $\sfA^{\shear}$,
\beqn
\label{eqn: shearing equivalence 3}
(-)^{\shear}: \sfA\mod(\Rep(\G_m)\mod) \xrightarrow{\simeq} \sfA^{\shear}\mod(\Rep_{\epsilon}^{\super}(\G_m)\mod).
\eeqn

\brem
\label{rem: C = C shear}
If $\sfA = A\mod$ for a graded $k$-algebra $A$,
then $\sfA^{\shear} \simeq A^{\shear}\mod$.
\erem

\brem
Because \labelcref{eqn: shearing equivalence 3} is a
\emph{symmetric monoidal equivalence}, the categorical
dimensions of dualizable objects in each category are
exchanged under the equivalence 
$\End(1_{\sfA}) \simeq \End(1_{\sfA^{\shear}})$. 
This is the advantage of working with
$\G_m$ super-representations in this setting.
\erem

\section{Hochschild invariants}
\label{sec: Hochschild invariants}

We assume familiarity with the theory of
traces in homotopical algebra, as outlined,
for example, in \cite{BZNtraces}.

\subsection{Hochschild homology}
\bdef
If $c \in \scrC$ is a dualizable object of the symmetric
monoidal $\infty$-category $\scrC$, we denote the
categorical dimension of $c$ by $\HH_{\bullet}^{\scrC}(c)$,
which we call the \emph{Hochschild homology of $c$}.
\edefn

The Hochschild homology has a canonical $S^1$-action,
where $S^1$ is considered as a group object in $\Grpd_{\infty}$.
More precisely, there is a canonical map
of $\bbE_1$-spaces
 \beqn
\label{eqn: S^1-action on HH}
S^1 \to \Omega_{\HH_{\bullet}^{\scrC}(c)}(\Map(1_{\scrC}, 1_{\scrC})).
\eeqn

When $\scrC$ is enriched over a rigid
monoidal category\footnotemark
\footnotetext{See \cite[\S 9.1.1]{GR17} for definition.} 
$\scrR$, we obtain from the canonical
map \labelcref{eqn: S^1-action on HH} a
map of $\bbE_1$-algebras in $\scrR$,
\[\scrR[S^1] \to \Omega_{\HH_{\bullet}^{\scrC}(c)}(\End_{\scrR}(1_{\scrC})),\]
where $\scrR[S^1]$ denotes the tensor product
$S^1 \otimes 1_{\scrR}$ inside $\scrR$.

\subsection{Tate fixed points}
Without being too precise, let $X$ be an object with
an action of a sufficiently nice topological group $G$. 
The Tate fixed points of this action is, loosely speaking,
the cofiber of a norm map comparing the $G$-coinvariants of $X$
to the $G$-invariants of $X$.

In order to state the Tate fixed point construction
precisely, and in sufficient generality for our purposes, 
we recall the following theorem of
Nikolaus--Scholze.

\bthm[{\cite[Theorem I.4.1]{NikolausScholze}}]
\label{thm: Tate construction}
Let $S$ be a Kan complex, and consider the $\infty$-category $\Sp^S =
\Fun(S, \Sp)$. Let $p : S \to \ast$ be the projection to the point. Then $p^*: \Sp \to \Sp^S$ 
has a left adjoint $p^! : \Sp^S \to \Sp$ given by homology, and a right adjoint $p_*: \Sp^S \to \Sp$
given by cohomology.
\begin{enumerate}[label=(\roman*), ref=\roman*]
\item \label{item: Tate 1} The $\infty$-category $\Sp^S$ is a compactly generated 
presentable stable $\infty$-category. For
every $s \in S$, the functor $s_! : \Sp \to \Sp^S$ takes compact objects to compact objects,
and for varying $s \in S$, the objects $s_! \bbS$ generate $\Sp^S$.

\item \label{item: Tate 2} There is an initial functor $p_*^T: \Sp^S \to \Sp$ with a natural 
transformation $p_* \to p_*^T$ with the property that $p_*^T$ vanishes on compact objects.

\item \label{item: Tate 3} The functor $p_*^T: \Sp^S \to \Sp$ is the unique functor with a natural transformation
$p_* \to p_*^T$ such that $p_*^T$ vanishes on all compact objects, and the fiber of $p_* \to p_*^T$
commutes with all colimits.

\item \label{item: Tate 4} The fiber of $p_* \to p_*^T$ is given by $X \mapsto p_!(D_S \otimes X)$ 
for a unique object $D_S \in \Sp^S$. The object $D_S$ is given as follows. 
Considering $S$ as an $\infty$-category, one has a functor 
$\Map : S \times S \to \Grpd_{\infty}$, sending a pair $(s, t)$ of points in $S$ to the space of paths
between $s$ and $t$. Then $D_S$ is given by the composite 
\[D_S : S \to \Fun(S, \Grpd_{\infty}) \xrightarrow{\Sigma_+^{\infty}} \Fun(S, Sp) = \Sp^S \xrightarrow{p_*} \Sp.\]

\item \label{item: Tate 5} The map $p_!(D_S \otimes -) \to p_*$ is final in the category 
of colimit-preserving functors from $\Sp^S$ to $\Sp$ over $p_*$, 
i.e. it is an assembly map in the sense of Weiss–Williams.

\item \label{item: Tate 6} Assume that for all $s \in S$ and $X \in \Sp$, one has $p_*^T(s_!X) \simeq 0$. 
Then there is a unique lax symmetric monoidal structure on $p_*^T$ which makes 
the natural transformation $p_* \to p_*^T$ lax symmetric monoidal.
\end{enumerate}
\ethm

In the case $S = BG$ for a topological group $G$, the category $\Sp^S$ is precisely
the category of spectra endowed with a $G$-action, and $p_!$ and $p_*$ are
the functors of $G$-coinvariants and $G$-invariants, respectively.

\bdef
Let $S = BG$ for a topological group $G$, and let $X \in \Sp^{BG}$. Then, we define
\[X^{tG} := p_*^T(X)\]
to be the \emph{Tate fixed points} of $X$.
\edefn

If $G$ is such that $S = BG$ satisfies the condition in (\labelcref{item: Tate 6}),
then functor of Tate fixed points,
\[(-)^{tG}: \Sp^{BG} \to \Sp\]
is lax symmetric monoidal by the same. In this case,
$(-)^{tG}$ preserves algebra objects in $\Sp^{BG}$,
as well as modules over such. In particular, if $R$ is a commutative
ring spectrum with trivial $G$-action, we obtain a functor,
\[(-)^{tG}: (R\mod)^{BG} \to (R^{tG})\mod,\] 
where $(R\mod)^{BG}$ denotes the category
$G\mod(R\mod) := R[G]\mod(R\mod)$, which is equivalent
by a free-forgetful adjunction to objects of $R\mod$ whose underlying
spectrum are equipped with a $G$-action.

\brem
\label{rem: Tate invariants}
As a special case, take $R=k$ and $G=S^1$. Then, $k^{tG} \simeq \kb$,
so we obtain the familiar functor,
\[(-)^{tS^1}: S^1\mod(\Vect_k) \to \kb\mod,\]
taking values in $2$-periodic complexes.
\erem

\brem
Note too that $p_*$ is lax symmetric monoidal, so similarly
obtain a functor of $G$-invariants,
\[(-)^G: (R\mod)^{BG} \to R^G\mod.\]
In particular, if $R=k$ and $G =S^1$ as in \cref{rem: Tate invariants},
we obtain a functor,
\[(-)^{S^1}: S^1\mod \to k[\b]\mod.\]
\erem

\subsection{Periodic cyclic homology}
We will define the periodic cyclic homology of module
categories over $R\mod$, for a commutative ring spectrum
$R$.

\bdef
Let $\scrC = (R\mod)\mod$, and let $C \in \scrC$ be a dualizable category.
The Hochschild homology of $C$ is naturally an object
of $R\mod$ equipped with an $S^1$-action. The 
\emph{periodic cyclic homology of $C$} is defined to be
\[\HP_{\bullet}^{\scrC}(C) := \HH_{\bullet}^{\scrC}(C)^{tS^1},\]
where $\HH_{\bullet}^{\scrC}(C)$ is regarded as an object
of $S^1\mod(R\mod)$. 
\edefn

\brem
We will typically use $\HH_{\bullet}^R(C)$ and $\HP_{\bullet}^R(C)$
to denote the Hochschild and periodic cyclic homology, respectively,
when $\scrC = R\mod$.
\erem

The special case of interest to us is the one in which $R$ is
$k$-algebra, specially when $R=k$ or $\kb \simeq C^*(BS^1;k)$, for $\b$
a fixed generator of the cohomology ring of $BS^1$. For the former,
periodic cyclic homology is a functor,
\[\HP_{\bullet}^k(-): \dgcat_k^{\infty} \to \kb\mod,\]
and for the latter, it is a functor,
\[\HP_{\bullet}^{\kb}(-): \kb\mod \to k[\b^{\pm 1}, u^{\pm 1}]\mod,\]
where $u$ is another copy of the fixed generator of $C^*(BS^1;k)$.


\printbibliography

\end{document}